\documentclass[a4paper,11pt]{article}

\usepackage{a4}
\usepackage{latexsym,amssymb}
\usepackage{citesort}
\usepackage{deleq}
\usepackage{ntheorem}
\usepackage{url}
\DeclareFontFamily{OT1}{rsfs}{}
\DeclareFontShape{OT1}{rsfs}{m}{n}{ <-7> rsfs5 <7-10> rsfs7 <10->
rsfs10}{} \DeclareMathAlphabet{\mycal}{OT1}{rsfs}{m}{n}
\newcommand{\mfd}{{\mycal D}}%
\newcommand{\mcV}{{\mycal V}}%
\newcommand{\mcF}{{\mycal F}}%
\newcommand{\tmcf}{\tmcF}%
\newcommand{\tmcF}{\,\,\,\widetilde{\!\!\!\mycal F}}%
\newcommand{\mfdb}{{\mycal D}_{\pM}}%
\newcommand{\mcO}{{\mycal O}}%

\newcommand{\cgg}{c}
\newcommand{\cgamma}{c}

\newcommand{\HH}{\mathbb{H}}
\newcommand{\bM}{\overline{M}}

\newcommand{\ppref}[1]{\ref{#1}, p.~\pageref{#1}}
\newcommand{\bel}[1]{\begin{equation}\label{#1}}
\newcommand{\beal}[1]{\begin{eqnarray}\label{#1}}
\newcommand{\beaa}{\begin{eqnarray*}}
\newcommand{\eeaa}{\end{eqnarray*}}
\newcommand{\eeal}[1]{\label{#1}\end{eqnarray}}
\newcommand{\bed}{\begin{deqarr}}
\newcommand{\eed}{\end{deqarr}}
\newcommand{\eedl}[2]{\arrlabel{#1}\label{#2}\end{deqarr}}
\newcommand{\eea}{\end{eqnarray}}

\newtheorem{Theorem}{\sc Theorem\rm}[section]
\newtheorem{Corollary}[Theorem]{\sc Corollary\rm}
\newtheorem{Lemma}[Theorem]{\sc Lemma\rm}
\newtheorem{Proposition}[Theorem]{\sc Proposition\rm}
\newtheorem{Definition}[Theorem]{\sc Definition\rm}
\newtheorem{hypotheses}[Theorem]{\sc Hypotheses\rm}

\theorembodyfont{\upshape}

\newtheorem{remark}[Theorem]{\sc Remark\rm}
\newtheorem{NRemark}[Theorem]{\sc Remark\rm}

\newcommand{\Proof}{\noindent\textsc{Proof:}\ }
\newcommand{\proof}{\Proof}
\newcommand{\R}{\mathbb{R}}

\newcommand{\ext}{{\mbox{\scriptsize \rm ext}}}

\newcommand{\dvYo}{dv_Y}

\newcommand{\QED}
   {\hfill$\hbox{\vrule height1.3ex width1.3ex depth.1ex}\ $
    \par\medskip}
\newcommand{\qed}{\QED}

\newcommand{\pd}[2]{\frac{\partial #1}{\partial #2}}

\newcommand{\rY}{\big|_{{Y}}}

\newcommand{\YI}{Y\times I}
\newcommand{\go}{\mathring{g}}

\newcommand{\cut}%
{\hbox{%
\,{\vrule height.1ex width.7ex depth.1ex}%
{\vrule height1.3ex width.1ex depth.1ex}\,}%
}

\newcommand{\Norm}[2]{\Vert#1\Vert_{#2}}

\newcommand{\nmi}[1]{\Norm{#1}{H^{1}_{*}(\YI)}}
\newcommand{\nyh}[1]{\Norm{#1}{H^{1/2}_{*}}}
\newcommand{\ran}{{\mathrm{range}\,}}
\newcommand{\loc}{\mathrm{loc}}

\newcommand{\hb}{\mathbb{H}}
\newcommand{\bN}{\mathbb{N}}
\newcommand{\bbC}{\mathbb{C}}

\newcommand{\bbR}{\mathbb{R}}
\newcommand{\bR}{\mathbb{R}}
\newcommand{\bZ}{\mathbb{Z}}
\newcommand{\bbZ}{\mathbb{Z}}
\newcommand{\bbT}{\mathbb{T}}

\newcommand{\cBK}{\mathcal{K}}
\newcommand{\cC}{\mathcal{C}}
\newcommand{\cD}{\mathcal{D}}
\newcommand{\cDb}{{\mathcal{D}_Y}}

\newcommand{\cK}{\mathcal{K}}
\newcommand{\cL}{\mathcal{L}}
\newcommand{\cO}{\mathcal{O}}

\newcommand{\cU}{\mathcal{U}}
\newcommand{\cutut}{\mathcal{U}}
\newcommand{\co}{\mathcal{U}}
\newcommand{\nmet}{h} 
\newcommand{\cN}{{\cal N}}
\newcommand{\tM}{\widetilde{M}}

\newcommand{\nablb}{\bar{\nabla}}
\newcommand{\slashD}{\mathcal{D}\kern-0.5em{/}}
\newcommand{\cl}[1]{c(#1)}  

\newcommand{\Cl}{\mathrm{C}\ell}  

\newcommand{\bref}[1]{(\ref{#1})}
\newcommand{\eq}[1]{(\ref{#1})}
\newcommand{\Eq}[1]{Equation~(\ref{#1})}

\newcommand{\Eqsone}[1]{Equations~(\ref{#1})}

\newcommand{\tfrac}[2]{{\textstyle\frac{#1}{#2}}}
\newcommand{\half}{\tfrac{1}{2}}
\newcommand{\supp}{\mathrm{supp\,}}
\newcommand{\twostar}{\hat{2}}

\newcommand{\be}{\begin{equation}}
\newcommand{\ee}{\end{equation}}
\newcommand{\bp}{\begin{Proposition}}
\newcommand{\ep}{\end{Proposition}}

\newcommand{\Riemn}{pseudo--Riemannian\ }

\newcommand{\gst}{>_{\!\!*}}

\newcommand{\ON}{{\cal ON}_{\!\!M,g}}
\newcommand{\ptwo}{p}
\newcommand{\wkp}{W^{k+1,p}}
\newcommand{\wkpg}{W^{k,p}}

\newcommand{\wlqg}{W^{\ell,q}}
\newcommand{\wlqgl}{\wlqg_\loc}
\newcommand{\wkpl}{\wkp_\loc}
\newcommand{\wkpgl}{\wkpg_\loc}
\newcommand{\divE}{\mathrm{div}\;E}
\newcommand{\divB}{\mathrm{div}\;B}

\newcommand{\hH}{{\hat H}}
\newcommand{\hK}{{\hat K}}
\newcommand{\hA}{{\hat A}}
\newcommand{\Y}{\mathrm{Im}(A)}

\newcommand{\X}{\mathrm{Ker}(A)}

{\catcode `\@=11 \global\let\AddToReset=\@addtoreset}
\AddToReset{equation}{section}
\renewcommand{\theequation}{\thesection.\arabic{equation}}

\newcounter{mnotecount}[section]

\newcommand{\mnote}[1]{}%

\newcommand{\piotr}[1]{}
\newcommand{\rob}[1]{}

\newcommand{\oldmnote}[1]{}
\newcommand{\oldnote}[1]{}

\newcounter{pcheckcount}[section]

\newcommand{\pcheck}[1]{}

\begin{document}

\title{Boundary value problems for Dirac--type equations, with applications}

\author{Robert A.\ Bartnik\thanks{Supported in part by the Australian
    Research Council.
\emph{Email}: \protect\url{bartnik@ise.canberra.edu.au}}\\
School of Mathematics and Statistics\\
University of Canberra\\
ACT 2601 Australia
\\\mbox{~}
\\
Piotr T.\ Chru\'sciel\protect\thanks{%
  Supported in part by grant from the Polish Committee for Scientific
  Research \# 2 P03B 073 15 and by the French Ministry for Foreign
  Affairs.
  \emph{Email}: \protect\url{chrusciel@univ-tours.fr}, URL \protect\url{www.phys.univ-tours.fr\~piotr}} \\
D\'epartement de Math\'ematiques\\
Facult\'e des Sciences\\
Parc de Grandmont\\
F37200 Tours, France } \maketitle
\begin{abstract}
  We prove regularity for a class of boundary value problems for first
  order elliptic systems, with boundary conditions determined by
  spectral decompositions, under coefficient differentiability
  conditions weaker than previously known. We establish Fredholm
  properties for Dirac-type equations with these boundary conditions.
  Our results include sharp solvability criteria, over both compact
  and non-compact manifolds; weighted Poincar\'e and
  Schr\"odinger-Lichnerowicz inequalities provide asymptotic control
  in the non-compact case.  One application yields existence of
  solutions  for the Witten equation with  a spectral boundary condition
  used by Herzlich in his proof of a geometric lower bound for the ADM
  mass of asymptotically flat 3-manifolds.


\end{abstract}


\tableofcontents

\section{Introduction}

Elliptic systems based on the Dirac equation arise frequently in
problems in geometry and analysis. Applications to positive mass and
related conjectures in general relativity motivate this paper, and
involve boundary value problems on compact and non-compact domains
\cite{Witten81,GHHP83,Herzlich97a,Herzlich:mass}.

Previous existence and regularity results
\cite{Morrey,APS75,Hormander85,Bunke,BoosWojc93} are insufficient for
these applications, for various reasons.  The Agmon-Douglas-Nirenberg
approach based on freezing coefficients and explicit kernels for the
constant coefficient inverse operator, leads only to boundary
conditions of Lopatinski-Shapiro type \cite{Morrey66}.  The
pseudo-differential operator approach \cite{Seeley66,Hormander85}
handles non-local boundary conditions such as the spectral projection
condition of Atiyah-Patodi-Singer \cite{APS75}, but the assumptions of
smooth coefficients and product-type boundary metric
\cite{APS75,Hormander85,Bunke,BoosWojc93} are unnatural and, as we
shall show, unnecessary.

In this paper we provide an essentially elementary proof of existence
and regularity for first order elliptic systems with ``Dirac-type''
boundary value conditions.  These encompass both pointwise
(Lopatinski-Shapiro) and non-local (spectral) boundary conditions, and
do not require product metric structures on the boundary.  We obtain
explicit necessary and sufficient conditions which ensure the
solvability of natural inhomogeneous boundary value problems, over
both compact and non-compact manifolds with compact boundary.

The coefficient regularity conditions, for both the elliptic system
and the boundary conditions, are rather general.  For example, they
are weaker than those in the pseudo-differential operator approach of
Marschall \cite{Marschall88}.  It seems likely that the boundary
conditions will admit some generalizations; the boundary data is
$H^{1/2}$ whereas there are recent results for a certain constant
coefficient Dirac equation with $L^2$ boundary values on a Lipschitz
hypersurface \cite{AxelssonEtal00}.

Note that there is an extensive literature on applications of Dirac
operators to index problems on compact and non-compact manifolds
\cite{BoosWojc93,Carron96} which we do not address, although many
aspects of our results are no doubt relevant to such applications; the
results here are focussed on applications to energy theorems in general
relativity.

The motivating example of the Dirac (Atiyah-Singer) operator is
described in some detail in \S\ref{secA}, where the
Schr\"odinger-Lichnerowicz identity with suitable boundary conditions
combines with a Lax-Milgram argument to reduce the existence question
to that of showing that a weak ($L^{2}$) solution of an adjoint
problem is in fact a strong ($H^{1}$) solution.  This
\emph{weak-strong} regularity property turns out to be the key
technical step, and the focus of much of the paper.  The difficult
case is regularity at the boundary; interior regularity is established
in \S\ref{secIR} using standard Fourier techniques, for general first
order elliptic systems.

\S\ref{secSC} reviews conditions under which a symmetric operator
has a complete set of eigenfunctions; these are used to to control
the boundary operator in later sections. In \S\ref{secB} we prove
regularity results at the boundary, for a class of operators much
broader than Dirac equations, with weak assumptions on the
continuity/regularity of the operator coefficients.  The main
technical tools are the $H^{1}$ identity \bref{B1:est}, and some
basic spectral theory.  The boundary conditions of \S\ref{secB}
follow from the requirements of the arguments of the regularity
theorem, and some additional work is required to apply them to
first order systems. This is carried out in \S\ref{secC}, for
equations of Dirac-type near the boundary, for which the boundary
operator is self-adjoint.  The resulting boundary conditions are
naturally presented in terms of graphs over the space of negative
eigenfunctions of the boundary operator.


The boundary value problems considered have a Fredholm property, and
admit an explicit solvability criteria involving solutions of the
homogeneous adjoint problem.  These properties are established for
compact manifolds with boundary in \S\ref{secDI}, and for a large
class of non-compact manifolds with boundary in \S\ref{secD}.  The
analysis of the non-compact case relies on two \emph{a priori}
inequalities: a weighted Poincar\'e inequality, and a
Schr\"odinger-Lichnerowicz inequality.  These inequalities imply the
manifold is \emph{non-parabolic at infinity} in the sense of
\cite{Carron96}. The weighted Poincar\'e inequality is established in
\S\ref{secF} in a number of cases, including the important cases of
manifolds with asymptotically flat or hyperbolic ends.  The
Schr\"odinger-Lichnerowicz inequality follows in applications from an
$H^1$ estimate derived from an identity of Schr\"odinger-Lichnerowicz
type.

In section~\ref{secE} we show that common pointwise and spectral
boundary conditions for the Dirac equation are elliptic in the
sense of our conditions.  These calculations form the basis for
\S\ref{Section:pmt}, which verifies several positive mass theorems
\cite{Witten81,GHHP83,Herzlich97a}.  Appendix \ref{AppAGC}
collects some relevant properties of tensor and spinor fields on
manifolds with $W^{k+1,p}$ differentiable structure and $W^{k,p}$
metric, $k>n/p$.

\section{The model problem}
\label{secA}
In this section we use the Riemannian Dirac equation to illustrate and
motivate the existence and regularity results of the following
sections.

Consider an oriented manifold $M$ with Riemannian metric $g$ and a
representation $c : \Cl(TM)\to \textrm{End}(S)$ of the Clifford
algebra $\Cl(TM)$ on some bundle $S$; with our conventions,
$$\cl{v}\cl{w}+\cl{w}\cl{v}=-2g(v,w)\;.$$   Clifford
representations are discussed in detail in
\cite{ABS64,LawsonMichelsohn89}.
$S$ carries an invariant inner product, $\langle \cl{v}\psi,
\cl{v}\psi\rangle = |v|^2\,\langle\psi,\psi\rangle =
|v|^2\,|\psi|^{2}$, with respect to which $\cl{v}$ is skew-symmetric,
for all vectors $v$.

A \emph{Dirac} connection \cite{LawsonMichelsohn89,BoosWojc93} is a
connection on the space of sections of $S$ which satisfies the
compatibility relation
\begin{equation}
  d\langle \phi,\cl{v}\psi\rangle = \langle
  \nabla\phi,\cl{v}\psi\rangle +\langle
  \phi,\cl{v}\nabla\psi\rangle+\langle \phi,\cl{\nabla v}\psi\rangle\;,
\label{A:compat}
\end{equation}
where $\nabla$ also denotes  the Levi-Civita connection on vector
fields.

Spin manifolds provide the fundamental example, with $S$ a bundle of
spinors associated with a Spin principal bundle which double covers
the Riemannian orthonormal frame bundle.  In this case there is a
covariant derivative $\nabla$ defined in terms of a local orthonormal
frame $e_k$, $k=1,\dots,n$, with Riemannian connection matrix
$\omega_{ij}(e_{k})=g(e_i,\nabla_{e_k}e_j)$, by
\begin{equation}
  \nabla_{e_{k}}\psi = D_{e_{k}}\psi^{I}\,\phi_{I} -
  {\textstyle{\frac{1}{4}}} \psi^{I}\,\omega_{ij}(e_{k})
  \cl{e^{i}e^{j}}\phi_{I}\;,
\label{A:gpsi}
\end{equation}
where $\psi=\psi^{I}\phi_{I}$ and $\phi_I$, $I=1,\dots,\dim S$, is a
choice of spin frame associated with the orthonormal frame $e_k$.  The
expression \bref{A:gpsi} may be abbreviated to
$\nabla=d-\frac{1}{4}\omega_{ij}e^{i}e^{j}$.  Note  that there
are other examples of Dirac bundles and connections, eg.~\cite[example
II.5.8]{LawsonMichelsohn89}.

The \emph{Dirac operator} of a Dirac connection $\nabla $ is
\begin{equation}
  \cD\psi = \cl{e^{i}}\nabla_{e_{i}}\psi\;;
\label{A:dirac}
\end{equation}
in the spin case this is sometimes called the \emph{Atiyah-Singer
  operator}.  When the spinor representation is
irreducible\footnote{Reducible representations lead to interesting
  formulas with $\tfrac{1}{4}R(g)$ replaced by more complicated
  curvature endomorphisms, \emph{c.f.}~\S\ref{Section:pmt}}, a
classical and very important computation \cite{Schrodinger32} shows
that
\begin{equation}
\label{A-1}
\cD^{2}\psi = \nabla^{*}\nabla\psi + {\textstyle{\frac{1}{4}}}R(g)\psi\;,
\end{equation}
where $R(g)$ is the (Ricci) scalar curvature of $g$.  This leads to
the Schr\"odinger-Lichnerowicz identity
\cite{Lichnerowicz63,Schrodinger32}
\begin{equation}
\label{A-2}
\left(|\nabla\psi|^{2} + {\textstyle{\frac{1}{4}}}R(g)|\psi|^{2}
  -|\cD\psi|^{2}\right) *\!1 =
d\left(\langle\psi,(\cl{e_{i}e_{j}}+g_{ij})\nabla^{j}\psi\rangle
  *\!e^{i}\right)\;,
\end{equation}
which when integrated over the compact manifold $M$ with
boundary\footnote{Throughout this paper we use the geometer's
  convention, that a manifold with boundary contains its boundary as a
  point set.}  $Y$ becomes
\begin{equation}
  \int_{M}\left(|\nabla\psi|^{2} +
    {\textstyle{\frac{1}{4}}}R(g)|\psi|^{2} -|\cD\psi|^{2}\right) =
  \oint_{Y}\langle\psi,\cl{ne^{A}}\nabla_{A}\psi\rangle\;.
\label{A:Lich}
\end{equation}
Here $n$ is the outer normal vector at $Y=\partial M$ and $\{e_{A}\}$
is a compatible orthonormal frame on $Y$.  The boundary term
may be simplified by introducing the boundary covariant derivative
\[
\nablb = d-\tfrac{1}{4}\omega_{AB}\cl{e^{A}e^{B}}\;,
\]
and the \emph{boundary} Dirac operator\footnote{Both $e^A\to
\cl{e^A}$ and $e^A\to\cl{ne^A}$ give representations of the Clifford
algebra of the boundary tangent space; the choice of $\cl{ne^A}$ is
made here for convenience \cite{Herzlich97a}.  }
\begin{equation}
\label{A:bdirac}
\cDb\psi=\cl{ne^{A}}\nablb_{A}\psi\;.
\end{equation}
Denoting the mean curvature by $H=H_{Y}=g(n,\nabla_{e_{A}}e^{A})$
gives
\begin{equation}
  \oint_{Y}\langle\psi,\cl{ne^{A}}\nabla_{A}\psi\rangle =
  \oint_{Y} \langle\psi,\cDb\psi + \half H \psi\rangle\;.
\label{A:Lichb}
\end{equation}
We use conventions which give $H=2/r>0$ for $M=\bR^{3}-B(0,r)$, the
exterior of a ball of radius $r$, with the outer normal
$n=-\partial_{r}$.  If $x$ is a Gaussian boundary coordinate ($x\ge0$
in $M$, $x=0$ on $Y$ and $\partial_x = -n$), then near the boundary we
have
\begin{equation}
\cD \psi =  -\cl{n}(\partial_x + \cDb + \half H)\psi\;.
\end{equation}

We now seek boundary conditions for which the equation $\cD\psi=f$
is solvable, following a well-known argument
\cite{ParkerTaubes82,GerochPerng94,Herzlich97a}. Suppose $M$ is a
compact\footnote{The case of an $(M,g)$ which is asymptotically
flat with compact interior, which is relevant to the positive mass
theorem, is discussed along similar lines in \S\ref{secD} and
\ref{Section:pmt}.}  manifold with non-negative scalar curvature,
$R(g)\ge0$, and $\cBK:H^{1/2}(Y)\to H^{1/2}(Y)$ is a bounded
linear operator such that
\begin{equation}
  \oint_{Y} \langle\psi,\cDb\psi + \half H \psi\rangle \le 0
  \qquad \textrm{whenever}\ \ \cBK\psi=0\;.
\label{A:bMle0}
\end{equation}
Suppose further that $M$ admits no parallel spinors.  Define the space
 $H^{1}_{\cBK}(M)$ as the completion of the smooth spinor fields with
 compact support (in $M\cup Y$) which satisfy the boundary condition
 $\cBK\psi=0$, in the norm
\begin{equation}
  \Norm{\psi}{H^{1}_{\cBK}(M)}^{2} :=
  \int_{M}(|\nabla\psi|^{2}+\tfrac{1}{4}R(g)|\psi|^{2})\;.
\label{A:normh1}
\end{equation}
The boundary condition \bref{A:bMle0} combined
with the Lichnerowicz identity \bref{A:Lich} and the curvature condition
$R(g)\ge0$ now ensures that the
bilinear form
\[
a(\psi,\phi) = \int_{M}\langle \cD\psi,\cD\phi\rangle, \quad
\phi,\psi\in H^{1}_{\cBK}(M)\;,
\]
is strictly coercive,
$a(\psi,\psi)\ge\Norm{\psi}{H^{1}_{\cBK}(M)}^{2}$.  For any spinor
field $f\in L^{2}(M)$, the linear functional $\phi\mapsto
\int_{M}\langle f,\cD\phi\rangle$ is bounded on $H^{1}_{\cBK}(M)$.
Coercivity and the Lax-Milgram lemma show there is a unique $\psi\in
H^{1}_{\cBK}(M)$ such that
\[
\int_{M}\langle \cD\psi-f,\cD\phi\rangle =0\qquad \forall \phi\in
H^{1}_{\cBK}(M)\;,
\]
and we would like to deduce that $\cD\psi=f$.  Now $\Psi:=\cD\psi-f\in
L^2(M)$ is a \emph{weak solution} of the Dirac equation; that is,
\begin{equation}
        \int_{M}\langle\Psi,\cD\phi\rangle = 0 \qquad \forall\ \phi\in
        H^{1}_{\cBK}(M)\;.
\label{A:weak}
\end{equation}
\underline{If} we could show that $\Psi$
is in fact a strong solution, that is, $\Psi\in H^{1}(M)$, then we
could integrate by parts to conclude
\[
\int_{M}\langle\cD\Psi,\phi\rangle + \oint_{Y} \langle
\Psi,\cl{n}\phi\rangle = 0 \qquad \forall \phi\in H^{1}_{\cBK}(M)\;,
\]
and thus $\cD\Psi=0$ and $\oint_{Y}\langle
\Psi,\cl{n}\phi\rangle=0$ for all $\phi\in H^{1/2}(Y)$ such that
$\cBK\phi=0$.  This would give the boundary condition $\Psi\rY \in
\cl{n}(\ker\cBK)^\perp$, which we suppose may be re-expressed as
$\widetilde\cBK\Psi=0$, for some ``adjoint'' boundary operator
$\widetilde\cBK$.  This would give $\Psi\in H^1_{\widetilde\cBK}$,
so if finally we suppose that $\widetilde\cBK$ also satisfies the
boundary positivity condition \bref{A:bMle0}, then we could
conclude from $a(\Psi,\Psi) =0$ and the coercivity of
$a(\cdot,\cdot)$ with respect to the norm
$\Norm{\cdot}{H^1_{\widetilde\cBK}}$, that $\Psi=0$ as desired.

The key technical difficulty in this classical argument lies in
establishing the ``Weak-Strong'' property, that weak ($L^2$) solutions
lie in $H^1$.  In the following sections we will prove this property
for a large class of elliptic systems, under rather general boundary
conditions; see \S\ref{secB} and \S\ref{secC}.

Two model boundary operators illustrate the possibilities for
achieving the required conditions. The \emph{APS} (or \emph{spectral
projection} \cite{APS75}) condition arose in Herzlich's work
\cite{Herzlich97a}:
\begin{equation}
\label{A:herzlich}
\label{APS}
\cBK  =P_{+}\;,
\end{equation}
where $P_{+}$ is the $L^{2}(M)$-orthogonal projection onto the
positive spectrum eigenspace of the boundary Dirac operator $\cDb$.
Using the relation $\cl{n} \cDb = - \cDb \cl{n}$, which shows that the
spectrum of $\cDb$ is symmetric about $0\in\bbR$, we find that
$\widetilde\cBK = \cBK$, provided there are no zero eigenvalues.

The eigenvalue estimate for $Y\simeq S^2$ of Hijazi and B\"ar
\cite{Bar92,Hijazi91}
\begin{equation}\label{Hbest}
|\lambda(\cDb)|\ge \sqrt{4\pi/\mathrm{area}(Y)},
\end{equation}
shows that in this case there are no zero eigenvalues.  In addition,
if we have the mean curvature condition
\begin{equation}
H_{Y}\le \sqrt{16\pi/\textrm{area}(Y)},
\label{A:BH}
\end{equation}
then $\cBK$ (and $\widetilde\cBK$) will satisfy the boundary
positivity condition \bref{A:bMle0}.  In conclusion, if $Y=\partial
M\simeq S^2$ satisfies \bref{A:BH}, then (assuming the Weak-Strong
property can be established) the above argument shows $\Psi=0$ and
thus the equation $\cD \psi =f$ with boundary condition $P_+\psi = 0$
is uniquely solvable, for any $f\in L^2(M)$.

The \emph{chirality} condition was used
in \cite{GHHP83,DouganMason91}.  For a slightly simplified version of
\cite{GHHP83}, suppose $M$ is a totally geodesic hypersurface in a
Lorentz spacetime, with future unit normal vector $e_0$, and consider
the connection on spacetime spinors, restricted to $M$.  Along $Y
=\partial M$ we define
\begin{equation}\label{A:eps}
\epsilon=\cl{e_{0}n},
\end{equation}
which satisfies the \emph{chiral conditions}
\begin{equation}
\label{A:chiral}
\begin{array}{rclcrcl}
\epsilon^{2} &\ =\ & 1\;,  &\quad&
\epsilon \cl{n} +\cl{n}\epsilon &\ =\ & 0 \;,\\[2pt]
\langle\phi,\epsilon\psi\rangle &=& \langle\epsilon\phi,\psi\rangle\;,
&\quad&
\epsilon\cDb + \cDb\epsilon & = & 0\;,
\end{array}
\end{equation}
and then the boundary operators
\begin{equation}\label{GHHP}
  \cBK_{\pm} = \half(1\pm\epsilon) \;.
\end{equation}
Assuming either of the two conditions $\cBK_\pm\psi = 0$ gives
$        \epsilon\psi=\mp \psi $ which implies
\begin{eqnarray*}
        \langle\psi,\cDb\psi\rangle
        &\ =\ \mp \langle\psi,\cDb\epsilon\psi\rangle
        &\ =\ \pm\langle\psi,\epsilon\cDb\psi\rangle
        \\[2pt]
        &\ =\   \pm\langle\epsilon\psi,\cDb\psi\rangle
        &\ =\   -\langle\psi,\cDb\psi\rangle=0\;.
\end{eqnarray*}
If we further assume that $H_{Y}\le0$ then \bref{A:bMle0} follows
directly.  In general relativity the condition $H_{Y}\le0$ is the
defining property for $Y$ to be a \emph{trapped surface}.

Since the $\cBK_\pm$'s are complementary orthogonal projections, we
have $(\ker\cBK_\pm)^\perp = \ker \cBK_\mp$, so $\psi\in \cl{n}(\ker
\cBK_\pm)^\perp$ exactly when $\cl{n}\psi\in \ker \cBK_\mp$, which
gives $\psi\in\ker\cBK_\pm$, and $\widetilde\cBK_\pm = \cBK_\pm$. In
this case we conclude (still assuming the Weak-Strong property can be
established) that if $H_Y\le0$ then $\cD \psi=f$, with either of the
boundary conditions $\epsilon \psi = \pm\psi$, is uniquely solvable.
These examples are discussed further in \S\ref{secE}.

\section{Interior Regularity}
\label{secIR}

In this section we establish regularity away from the boundary for
weak ($L^2$) solutions of first order elliptic systems.  We
consider equations of the form
\begin{equation}
\label{IR-1}
  \cL u := a^j \partial_j u + b u =f\;,
\end{equation}
where $u,f$ are sections respectively of $N$-dimensional real vector
bundles $E, F$, both over an $n$-dimensional manifold $M$ without
boundary\footnote{If $M$ has boundary $\partial M\ne\emptyset$, then
  the interior $\mathring{M}=M-\partial M$ is a (noncompact) manifold
  without boundary, to which the results of this section will apply.},
and $a^j,b$, $j=1,\ldots,n$, are sections of the bundle of
endomorphisms of $E$ to $F$.  We assume that $E,F$ are equipped with
fixed smooth inner products, denoted by $\langle\cdot,\cdot\rangle$.
The length determined by $\langle\cdot,\cdot\rangle$ will be denoted
invariably by $|u|^2 = \langle u,u\rangle$.  To simplify notation, the
respective bundles usually will be understood, and thus $L^2(M)$ will
generally mean $L^2\Gamma(E)$, the space of $L^2$ sections of $E$, or
$L^2\Gamma(F)$, depending on context.

\begin{remark} There is no loss of generality in
  considering real bundles, since complex and quaternionic bundles may
  be viewed simply as real bundles with additional algebraic
  structure.  For example, a Hermitean vector space of dimension $n$
  is equivalent to a real vector space of dimension $2n$ with a skew
  endomorphism $J$ satisfying $J^2=-1$, with the Hermitean inner
  product $(\ ,\ )$ and real inner product $\langle\ ,\ \rangle$
  related by
$
(u,v) = \langle u,v\rangle -i\langle u,Jv\rangle
$.
\end{remark}

Define the indices $\twostar=\twostar(n)$, $n^*=n^*(n)$ by
\begin{equation} \label{IR-1b}
\begin{array}{lll}
   \twostar =\frac{ 2n}{n-2},\ \  & n^* = n & \textrm{for }\ n\ge 3 ,
   \\[3pt]
 \twostar = 10^6, &  n^* = \frac{2}{1-2/\twostar}\ \ & \textrm{for }\ n = 2,
       \\
   \twostar = \infty\;, &  n^* = 2  & \textrm{for }\ n = 1\;,
\end{array}
\end{equation}
where $10^6$ represents any large constant.
Note that if $M$ admits a Sobolev inequality with constant $C_S$
\[
  \Norm{u}{L^{\twostar}} \le C_S \Norm{u}{H^1},
\]
then we also have
\begin{equation}
\label{IR-1c}
  \Norm{fu}{L^2} \le C_S \Norm{f}{L^{n^*}}\Norm{u}{H^1}\;.
\end{equation}
 Another basic fact is
the inequality
\[
   \Norm{fg}{W^{1,p}} \le C\Big(\Norm{f}{L^\infty}\Norm{g}{W^{1,p}}
                         +\Norm{g}{L^\infty}\Norm{f}{W^{1,p}}\Big),
\]
which shows that $W^{1,n^*}\cap C^0$ (in particular) forms a ring
under addition and multiplication of functions.  For $n=1,2$ the $C^0$
is superfluous here, of course.

With one exception, it suffices to assume throughout that the
underlying manifold has a $C^\infty$ differentiable structure.  The
exceptional point arises in \S\ref{secE} in the construction of
approximately Gaussian coordinates in a neighbourhood of the boundary,
when the metric has low regularity.  The description in Appendix
\ref{AppAGC} of $W^{k+1,p}$ differential structures, $k>n/p$,
establishes the necessary consistency conditions in this case.

We will assume that $a^j$, $b$ satisfy the regularity conditions
\begin{equation}
\label{IR-2}
\begin{array}{rcl}
a^j& \in& W^{1,n^*}_{\loc}(M)\cap C^0(M)\;,
\\[2pt]
b &\in &L^{n^*}_{\loc}(M)\;.
\end{array}
\end{equation}

The conditions \bref{IR-2} are preserved by bundle frame changes in
$W^{1,n^*}\cap C^0$, by the above ring property. In particular, even
if the bundle metrics $\langle\cdot,\cdot\rangle$ on $E, F$ are only
in $W^{1,n^*}\cap C^0$, by the Gram-Schmidt process we may construct
$W^{1,n^*}\cap C^0$ frame changes which make the metric coefficients
constant.  Since this changes the operator coefficients $a^j$, $b$
respectively by $W^{1,n^*}\cap C^0$, $W^{n^*}$ affine linear
transformations, there is no loss of generality in assuming the
metrics on $E, F$ to be locally constant.

The conditions \bref{IR-2} mean that $M$ can be covered by open
neighbourhoods $\mcO_\alpha$ with $W^{1,n^*}\cap C^0$ bundle
transition functions, such that the local coefficients $a^j$, $b$
satisfy the stated regularity.  Frame changes satisfying Sobolev
conditions are also discussed in detail in Appendix~\ref{AppAGC}.

We require that $a^j$ satisfy the ellipticity condition, that for each
$p\in M$ there is a coordinate neighbourhood $p\in U\subset M$ and a
constant $\eta>0$ such that
\begin{equation}
\label{IR-3}
\eta^2 |\xi|^2 |V|^2 \le |\xi_j a^j(x) V|^2 \le \eta^{-2} |\xi|^2 |V|^2\;,
\end{equation}
for all $x\in U$, $\xi\in T_x^*M$ and $V\in E_x$, where $|\xi|^2$ is
measured by a fixed background metric $\go$, which we may assume to be
$C^\infty$.  Note that \bref{IR-3} implies the fibres of $E,F$ must be
of the same dimension.

A \emph{weak solution} of \bref{IR-1} is $u\in L^2_{\loc}(E)$ such
that
\begin{equation}
\label{IR-4}\int_M \langle \cL^\dagger\phi,u\rangle \,dv_M =
\int_M\langle\phi,f\rangle\, dv_M\;,
\end{equation}
for all $\phi\in C_c^\infty(M)$, where $dv_M = \gamma dx$, $\gamma>0$,
is a coordinate-invariant volume measure on $M$ with $\gamma \in
W^{1,n^*}(U)\cap C^0(U)$ and $dx$ is coordinate Lebesgue measure, in
any local coordinate neighbourhood $U$. Here the formal 
adjoint $\cL^\dagger$ is defined with respect to $dv_M$ and the inner
products on $E, F$.  Thus in local coordinates, 
\begin{equation}
\label{IR-5}
\cL^\dagger\phi = -{}^ta^j \partial_j\phi + ({}^tb -
\gamma^{-1}\partial_j( {}^ta^j\gamma ))\phi\;,
\end{equation}
where the transposes ${}^ta^j$ are defined with respect to the local
framing forms of the inner products of $E, F$.

The proof proceeds by establishing various special cases, starting
with a constant coefficient operator acting on sections of a trivial
bundle $E$ over the torus $\bbT^n$.  This type of argument is very
standard.

\begin{Proposition}
\label{PropIR1}
   Suppose $u\in
   L^2(\bbT^n)$ is a weak solution of $\cL_0u=f$ where $f\in
   L^2(\bbT^n)$, $\cL_0 = a_0^j\partial_j$ with $a_0^j$ constant and
   satisfying the ellipticity condition \bref{IR-3}.  Then $u\in
   H^1(\bbT^n)$.
\end{Proposition}

\Proof We regard $\bbT^n = \bR^n/\bbZ^n$.  Fix a mollifier
$\phi_\epsilon = \epsilon^{-n} \phi((x-y)/\epsilon) \in C^\infty(\bR^n)$ with
$\phi(-x)=\phi(x)$, and set $u_\epsilon = \phi_\epsilon*u \in
C^\infty(\bbT^n)$.  Then
\[
   \cL_0 u_\epsilon = \int_{\bbT^n} a_0^j \frac{\partial}{\partial x^j}
   \phi_\epsilon(x-y)u(y)\, dy
\]
and thus the definition of weak solution gives
\begin{eqnarray*}
   \int_{\bbT^n} \langle \psi, \cL_0u_\epsilon\rangle dx
   &=& \int_{\bbT^n}\int_{\bbT^n} \phi_\epsilon(x-y)
       \langle \cL_0^* \psi(x), u(y)\rangle \,dy\,dx
\\
   &=& \int_{\bbT^n}\int_{\bbT^n} \langle -{}^ta_0^j
        \frac{\partial}{\partial y^j}
        \phi_\epsilon(y-x)\psi(x),u(y)\rangle \,dy\,dx
\\
   &=& \int_{\bbT^n} \langle \cL_0^* (\phi_\epsilon * \psi)(y),u(y)\rangle \,dy
\\
   &=& \int_{\bbT^n} \langle \phi_\epsilon * \psi(y),f(y)\rangle \,dy
\\
   &=& \int_{\bbT^n} \langle \psi(y),\phi_\epsilon*f(y)\rangle \,dy\;.
\end{eqnarray*}
Thus $\cL_0 u_\epsilon = f_\epsilon = \phi_\epsilon*f$, and we note that
$f_\epsilon\to f$ strongly in $L^2$.
Now the ellipticity condition \bref{IR-3} and the Plancherel theorem ensure
that for all $v\in H^1(\bbT^n)$,
\begin{eqnarray*}
   \int_{\bbT^n} |\partial v|^2 \,dx
   &=& \int_{\bbT^n} |\xi|^2|\hat{v}|^2\, d\xi
\\
   &\le& \eta^{-1} \int_{\bbT^n} |a^j_0 \xi_j\hat{v}|^2\,d\xi
\\
   &=& \eta^{-1} \int_{\bbT^n} |\cL_0v|^2\,dx\;,
\end{eqnarray*}
and thus
\[
   \int_{\bbT^n} |\partial u_\epsilon|^2 \,dx
   \le \eta^{-1} \int_{\bbT^n} |f_\epsilon|^2\, dx\;.
\]
Since $u,f\in L^2$, it follows that $u_\epsilon\to u$ strongly in $H^1$.
\QED

\begin{Proposition}
\label{PropIR2}
Under the conditions of Proposition~\ref{PropIR1}, the map
$\cL_0+\lambda:H^1(\bbT^n)\to L^2(\bbT^n)$ where $\lambda=\pi\eta$, is
uniquely invertible, and for all $u\in H^1(\bbT^n)$,
\begin{equation}
\label{IR-7}
   \Norm{u}{H^1} \le \sqrt{5}/\eta \Norm{(\cL_0+\lambda)u}{L^2}\;.
\end{equation}
\end{Proposition}

\Proof Write $u=\sum_{k\in\bbZ^n} u_k e^{2\pi i k\cdot x}$, where the
coefficients $u_k = \int_{\bbT^n}u(x)e^{-2\pi i k\cdot x}dx$ are
valued in $\bbC^N$, the complexification of the real vector space
modelling the fibres of $E$. We then have
\[
  \int_{\bbT^n} |(\cL_0+\lambda)u|^2 \,dx = \sum_{k\in\bbZ^n}
   |(2\pi i k_j a^j_0 +\lambda)u_k|^2,
\]
and using the vector length inequality $|a+b|^2 \ge \chi |a|^2 -
\chi/(1-\chi)|b|^2$ with $\chi=\half$, we find
\[
  |(2\pi i k_j a^j_0 +\lambda)u_k|^2 \ge \half  |2\pi i k_j a^j_0 u_k|^2
  - \lambda^2 |u_k|^2.
\]
For $k\ne0$ this is greater than $ \pi^2\eta^2 |u_k|^2$, whilst for $k=0$
we have $ |(2\pi i k_j a^j_0 +\lambda)u_k|^2 = \lambda^2 |u_k|^2 =
\pi^2\eta^2|u_k|^2$, hence
\[
  \int_{\bbT^n} |(\cL_0+\lambda)u|^2 \,dx \ge \sum_{k\in\bZ^n}
  \pi^2\eta^2|u_k|^2 = \pi^2\eta^2 \int_{\bbT^n} |u|^2 \,dx\;,
\]
which shows $\cL_0+\lambda$ has trivial kernel. Choosing
$\chi=1-1/(2|k|)$ shows in fact that
\[
  |(2\pi i k_ja^j_0 + \lambda)u_k|^2 \ge \pi^2\eta^2 (2|k|-1)^2|u_k|^2,
\]
for all $k\in\bZ^n$ and all $u_k\in\bbC^N$.  Since $2|k|-1 \ge 1$ for
all $k\in\bZ^n$, we obtain \bref{IR-7}.  Moreover, this shows
also that the $N\times N$ complex matrices $2\pi i k_ja^j_0 + \lambda$
are invertible for any $k\in\bZ^n$, which gives a direct construction
of the inverse of the operator $\cL_0+\lambda$.
\QED

\begin{Theorem}
\label{ThmIR3}
Suppose $u\in L^2(\bbT^n)$ is a weak solution of
\begin{equation}
\label{IR-8}
  \cL_0 u + B_0u + B_1u = f
\end{equation}
where $f\in L^2$ and $\cL_0=a_0^j\partial_j$ is a constant coefficient
first order operator satisfying the conditions of Proposition
\ref{PropIR2} with ellipticity constant $\eta$, where
$B_1:L^2\to L^2$ is bounded, and where
$B_0:H^1\to L^2$ is a linear map satisfying
\begin{equation}
\label{IR-8b}
  \Norm{B_0}{H^1\to L^2}\le \eta/3,\quad \Norm{B_0^\dagger}{H^1\to
  L^2} \le \eta/3\;,
\end{equation}
where $B_0^\dagger$ is the $L^2(\bbT^n)$-adjoint of $B_0$.  Then $u\in
H^1(\bbT^n)$ is a strong solution of \bref{IR-8}, and there is a
constant $C$, depending only on $\eta$ and $\Norm{B_1}{L^2\to L^2}$,
such that
\begin{equation}
\label{IR-8c}
\Norm{u}{H^1} \le C\, (\Norm{f}{L^2} + \Norm{u}{L^2})\;.
\end{equation}
\end{Theorem}

\Proof Construct the iteration sequence $w^{(k)}\in H^1$,
$k=0,1,\dots$ by defining $w^{(k+1)}$ to be the solution of
\begin{equation}
\label{IR-9}
  (\cL_0+\lambda)w^{(k+1)} = -B_0w^{(k)} + \tilde{f}\;,
\end{equation}
with $w^{(0)}=0$, where $\tilde{f}=f+\lambda u - B_1u\in L^2$ by the
assumptions.  This equation with $\lambda = \pi\eta$ is uniquely
solvable by Proposition~\ref{PropIR2}.  The difference $v^{(k+1)} =
w^{(k+1)}-w^{(k)}$ satisfies $(\cL_0+\lambda)v^{(k+1)} = -B_0v^{(k)}$ and
the estimate \bref{IR-7} shows that
\[
\frac{\eta}{\sqrt{5}}\Norm{v^{(k+1)}}{H^1}
\le \Norm{B_0v^{(k)}}{L^2} \le \frac{\eta}{3}\Norm{v^{(k)}}{H^1}.
\]
The iteration is thus a contraction and converges in $H^1$, to $w\in
H^1$ satisfying $(\cL_0+B_0+\lambda)w = \tilde{f}$, and then $v=u-w\in L^2$
is a weak solution of $(\cL_0+B_0+\lambda)v=0$.  Now $\cL_0^\dagger$ is also
elliptic with the same ellipticity constant $\eta$, so there is $z\in
H^1$ satisfying $(\cL_0^\dagger + B_0^\dagger+\lambda)z=v$.  Since $v$ is a
weak solution,
\[
   \int_{\bbT^n} \langle (\cL_0^\dagger + B_0^\dagger+\lambda)\phi,v\rangle
  \,dx=0 \quad \forall \phi\in H^1(\bbT^n),
\]
we may test with $\phi=z$ to see that $\int |v|^2 =0$ and $v=0$.  Thus
$u=v+w=w\in H^1$ as required.
By Proposition~\ref{PropIR2} and \bref{IR-8b}, we have
\begin{eqnarray*}
\frac{\eta}{\sqrt{5}} \Norm{u}{H^1}
&\le& \Norm{(\cL_0+B_0+B_1)u}{L^2} + \Norm{B_0u}{L^2} + \Norm{B_1u}{L^2}
       +\Norm{\lambda u}{L^2}
\\
&\le& \Norm{f}{L^2} + \frac{\eta}{3}\Norm{u}{H^1}
       + (\Norm{B_1}{L^2\to L^2}+\eta\pi)\Norm{u}{L^2}\;.
\end{eqnarray*}
Since $\sqrt{5}<3$, the estimate \bref{IR-8c} follows.
\QED

Next we consider operators with non-constant coefficients.  Let $C_S$
be the $\bbT^n$ Sobolev constant
\begin{equation}
\label{IR-10}
  \Norm{u}{L^{\twostar}(\bbT^n)} \le C_S \Norm{u}{H^1(\bbT^n)}\;,
\end{equation}
where $\twostar$ is defined in \bref{IR-1b}.

\begin{Proposition}
\label{PropIR4}
Suppose $u\in L^2(\bbT^n)$ is a weak solution of the system of equations
\begin{equation}
\label{IR-11}
   \cL u := a^j\partial_j u + bu = f
\end{equation}
where $f\in L^2$ and the coefficients $a^j\in W^{1,{n^*}}\cap C^0$,
$b\in L^{n^*}$ satisfy
\begin{eqnarray}
\label{IR-12}
   \Norm{a^j-a_0^j}{L^\infty} &\le& \frac{\eta}{10}
\\
\label{IR-13}
   \Norm{\partial_j a^j}{L^{n^*}} &\le& \frac{\eta}{10 C_S}
\end{eqnarray}
where $a_0^j$, $j=1,\dots,n$, are constant matrices with ellipticity
constant $\eta$.  Then $u\in H^1(\bbT^n)$ is a strong solution of
\bref{IR-11} .
\end{Proposition}

\Proof It will suffice to show that $\cL$ admits a decomposition
satisfying the conditions of Theorem \ref{ThmIR3}.  Since $L^\infty$
is dense in $L^{n^*}$, for any $\epsilon>0$ we may find $b_0\in L^{n^*}$,
$b_1\in L^\infty$, such that $b=b_0+b_1$ and
$\Norm{b_0}{L^{n^*}}<\epsilon$.  We choose $\epsilon = \eta/(10C_S)$.
Then $B_0u:=b_0u +(a^j-a_0^j)\partial_ju$ satisfies
\[
  \Norm{B_0u}{L^2} \le \Norm{b_0}{L^{n^*}}\Norm{u}{L^{\twostar}} +
  \Norm{a^j-a_0^j}{L^\infty} \Norm{\partial u}{L^2}.
\]
Using \bref{IR-13} and the Sobolev inequality \bref{IR-10} gives
\[
   \Norm{B_0u}{L^2} \le \frac{\eta}{10}(\Norm{u}{H^1}+ \Norm{\partial
   u}{L^2}) \le \frac{\eta}{3} \Norm{u}{H^1},
\]
so $\Norm{B_0}{H^1\to L^2}\le \eta/3$.  Clearly $B_1u:= b_\infty u$ is
bounded on $L^2$, and it remains to verify the $H^1\to L^2$ bound on
the adjoint operator
\[
  B_0^\dagger w := {}^tb_0 w - ({}^ta^j-{}^ta_0^j)\partial_jw
  - \partial_j({}^ta^j-{}^ta_0^j) w .
\]
Again using the Sobolev inequality and the conditions
\bref{IR-12},\bref{IR-13} we find
\begin{eqnarray*}
  \Norm{B_0^\dagger w}{L^2} &\le& C_S \Norm{b_0}{L^{n^*}}\Norm{w}{H^1} +
  \Norm{{}^ta^j-{}^ta_0^j}{L^\infty} \Norm{\partial w}{L^2} +
  C_S \Norm{\partial_ja^j}{L^{n^*}}\Norm{w}{H^1}
\\
  &\le& \frac{\eta}{3} \Norm{w}{H^1}\;,
\end{eqnarray*}
so the conditions of Theorem \ref{ThmIR3} are met and the result
follows.  \QED

On a general compact manifold we define the Sobolev space $H^1(M)$ by
the norm
\begin{equation}
\label{IR-14b}
\Norm{u}{H^1(M)}^2 = \int_M (|\nabla u|^2 +  |u|^2)\,dv_M\;,
\end{equation}
where the lengths $|u|^2$, $|\nabla u|^2$ are measured using the metric
$\langle\ ,\ \rangle$ on sections of $E$ and a fixed smooth background
metric $\go$ on $TM$, and where $\nabla$ is a (covariant) derivative
defined in local coordinates on $M$ and a local framing on $E$ by
\begin{equation}
\label{IR-14c}
\nabla_i = \partial_i - \Gamma_i\;.
\end{equation}
We assume the charts on $E,M$ are such that
\begin{equation}
\label{IR-14d}
 \Gamma_i\in L^{n^*}_\loc \;.
\end{equation}
Note we do not require that $\nabla$ be compatible with the metric on
$E$.
If $M$ is compact then the space $H^1(M)$ is independent of the choice
of covariant derivative:
\begin{Lemma}\label{LemmaIR5a}
Suppose $M$ is compact and $\nabla, \hat{\nabla}$ are covariant
derivatives satisfying \bref{IR-14d}.  Then there is $C>0$ such that
for all $u\in H^1(M)$,
\begin{equation}\label{IR-14e}
C^{-1} \int_M (|\nabla u|^2+|u|^2)dv_M
\le  \int_M (|\hat{\nabla} u|^2+|u|^2)dv_M
\le C  \int_M (|\nabla u|^2+|u|^2)dv_M\;.
\end{equation}
Moreover, there is a constant $C_S$, depending on $M,\nabla$, such that
\begin{equation}\label{IR-14f}
\left(\int_M |u|^{\hat{2}} dv_M\right)^{2/\hat{2}}
\le C_S  \int_M (|\nabla u|^2+|u|^2)dv_M\;.
\end{equation}
\end{Lemma}

\Proof
There is a finite covering of $M$ by charts
$U_\alpha$ with a corresponding partition of unity $\phi_\alpha$.
Using the Sobolev inequality for $U_\alpha\subset \bR^n$, in  each chart we may
estimate the localisation $u_\alpha=\phi_\alpha u$ by
\[
\int_{U_\alpha} |\partial u_\alpha|^2\,dx \le C\int_M (|\nabla
u_\alpha|^2+|u_\alpha|^2)\,dv_M \;,
\]
where $C$ depends also on the decomposition $\Gamma=\Gamma^\infty +
\Gamma^{n^*}\in L^\infty+L^{n^*}$, with $\Gamma^{n^*}$ small.
Again using the $\bR^n$ Sobolev inequality and $\Gamma,\hat{\Gamma}\in
L^{n^*}$ we have
\begin{eqnarray*}
  \int_M (|\hat{\nabla}u|^2 + |u|^2)\,dv_M
&\le& C\sum_\alpha \int_{U_\alpha} (|\hat{\nabla}u_\alpha|^2 +
|u_\alpha|^2)\,dx
\\
&\le& C\sum_\alpha \int_{U_\alpha} (|\partial u_\alpha|^2 + |u_\alpha|^2)\,dx
\\
&\le& C\sum_\alpha \int_{U_\alpha} (|\nabla u_\alpha|^2 + |u_\alpha|^2)\,dv_M\;,
\end{eqnarray*}
from which the equivalence of the norms follows easily.  The Sobolev
inequality follows from very similar arguments.
\QED

We may now complete the proof of interior regularity.

\begin{Theorem} \label{ThmIR5}
Suppose $M$ is a $C^\infty$ $n$-dimensional manifold without boundary,
and $E,F$ are real vector bundles over $M$, each with fibres modelled
on $\bR^N$.  Suppose $u\in L^2_{\loc}(M)$ is a weak solution of
$\cL u=f$, where $\cL$ is a first order operator satisfying the conditions
(\ref{IR-2},\ref{IR-3}).  Then $u\in H^1_{\loc}(M)$ and $u$ is a strong
solution of $\cL u=f$.  Moreover, if $M$ is compact there is a constant
$C>0$, depending on $a^j,b$ and $\Gamma$, such that for all $u\in H^1(M)$,
\begin{equation}
\label{IR-14}
   \Norm{u}{H^1(M)} \le C(\Norm{\cL u}{L^2(M)} + \Norm{u}{L^2(M)})\;.
\end{equation}
\end{Theorem}

\Proof Since $\cL$ is locally of the form $\cL u =a^j\partial_j u+bu$ with
$a^j\in W^{1,{n^*}}_{\loc}\cap C^0$, $b\in L^{n^*}_{\loc}$, for each $p\in M$
there is a coordinate neighbourhood $U$ and a constant $\eta>0$ such
that $\eta$ is the ellipticity constant of $a^j_0 = a^j(p)$, and with
respect to the local trivialisation of $E|_U\simeq U\times \bR^N$ we
have the bounds
\begin{eqnarray*}
   \Norm{a^j-a_0^j}{L^\infty(U)} &\le& \frac{\eta}{10},
\\
   \Norm{\partial_j a^j}{L^{n^*}(U)} &\le& \frac{\eta}{10 C_S}\;,
\end{eqnarray*}
where we assume without loss of generality that $U=Q_R=(0,R)^n$ is a
cube of side length $R\le1$.  By paracompactness there is a locally
finite countable covering $\{p_\alpha,U_\alpha\}_{\alpha\in\bZ}$ of
$M$ by such charts, with a subordinate $C^\infty$ partition of unity
$\{\phi_\alpha\}_{\alpha\in\bZ}$.  Noting that $\supp(\phi_\alpha
u)\Subset Q_R$ and that $\phi_\alpha u$ satisfies
\[
   \cL(\phi_\alpha u) = \phi_\alpha f + a^j
   \partial_j\phi_\alpha u
\]
weakly, we see that it suffices to consider the case where $\supp\,u
\Subset Q_R$.  Assuming this, rescaling by $y=x/R, x\in Q_R$ and defining
$\tilde{u}(y)=u(x)$, $\tilde{f}(y)=Rf(x)$, $\tilde{a}^j(y)=a^j(x)$ and
$\tilde{b}(y)=Rb(x)$, it follows that $\tilde{u}\in L^2(\bbT^n)$ is a weak
solution of
\[
   \tilde{a}^j \frac{\partial}{\partial y^j} \tilde{u}(y)
  + \tilde{b}(y) \tilde{u}(y) = \tilde{f}(y).
\]
In particular we have $\tilde{b}\in L^{n^*}(\bbT^n)$ and
\begin{eqnarray*}
   \Norm{\tilde{a}^j-a_0^j}{L^\infty(\bbT^n)}
 &=& \Norm{{a}^j-a_0^j}{L^\infty(Q_R)} \le \eta/10,
\\
   \Norm{\partial_{y^j}\tilde{a}^j}{L^{n^*}(\bbT^n)}
 &\le& \Norm{\partial_j{a}^j}{L^n(Q_R)} \le \eta/10\;.
\end{eqnarray*}
The conditions of Proposition~\ref{PropIR4} are satisfied, so
$\tilde{u}\in H^1(\bbT^n)$ and thus $u\in H^1_{\loc}(M)$.

When $M$ is compact there is a finite covering by charts
$\{p_\alpha,U_\alpha\}$, and by Theorem \ref{ThmIR3}, in each chart we may
estimate the localisation $u_\alpha=\phi_\alpha u$ by
\[
\Norm{u_\alpha}{H^1(U_\alpha)} \le C_\alpha
       (\Norm{\cL u_\alpha}{L^2(M)}+\Norm{u_\alpha}{L^2(M)}) ,
\]
where the $H^1(U_\alpha)$ norm uses the coordinate partial derivatives
$\partial_i$ in $U_\alpha$.  To estimate $\int_{U_\alpha}|\nabla
u_\alpha|^2$, note that the Sobolev inequality \bref{IR-1c}
in $U_\alpha$ gives
\[
   \int_{U_\alpha}|\Gamma u_\alpha|^2 \, dv_M
\le C \Norm{\Gamma}{L^{n^*}(U_\alpha)}^2 \int_{U_\alpha}
    (|\partial u_\alpha|^2 + |u_\alpha|^2)\,dv_M\;,
\]
so by the $H^1(U_\alpha)$ estimate we have
\[
\Norm{u_\alpha}{H^1(M)} \le C_\alpha
       (\Norm{\cL u_\alpha}{L^2(M)}+\Norm{u_\alpha}{L^2(M)}) ,
\]
for some constant $C_\alpha$ depending also on
$\Norm{\Gamma}{L^{n^*}(U_\alpha)}$.  Since $u=\sum u_\alpha$ and $\cL
u_\alpha = \phi_\alpha \cL u + \partial_j(\phi_\alpha) a^j u$, with
$|\partial\phi_\alpha|\le c$ and $|\phi_\alpha| \le 1$, the estimate
\bref{IR-14} follows easily.  \QED

The constant $C$ of \bref{IR-14} can be controlled by
$\Norm{a^j}{W^{1,p}}$, $\Norm{b}{L^{p}}$ for any $p>n^*$, or by otherwise
controlling the decompositions $\partial_j a^j,b\in L^\infty+L^{n^*}$.

Higher regularity follows easily from Theorem \ref{ThmIR5} by a
standard bootstrap argument:
\begin{Theorem}\label{ThmIR6}
Suppose $u\in L^2_{\loc}$ is a weak solution of $\cL u=f$ in the situation of
Theorem \ref{ThmIR5}, where the coefficients of $\cL u=f$ satisfy the
regularity conditions
\begin{equation}
\label{IR-15}
  a^j\in W^{k,n^*}_{\loc}\cap C^0 ,\ \  b \in
  W^{k,{n^*}}_{\loc},\quad\textrm{and }\ f\in H^k_{\loc}\;,
\end{equation}
for some integer $k\ge 1$.  Then $u\in H^{k+1}_{\loc}$.  If $M$ is a
compact manifold without boundary
then there is a constant $C=C(k,\cL)$, depending on $k$ and
$\Norm{a^j}{W^{k,{n^*}}}$, $\Norm{b^j}{W^{k,{n^*}}}$, such that
\begin{equation}
 \label{IR-16}
   \Norm{u}{H^{k+1}(M)} \le C(\Norm{f}{H^k(M)} + \Norm{u}{H^k(M)})\;.
\end{equation}
Thus for any $u\in L^2(M)$ such that $\cL u$ (defined weakly) satisfies
$\cL u\in H^k(M)$, we have
\begin{equation}
 \label{IR-16b}
   \Norm{u}{H^{k+1}(M)} \le C(\Norm{\cL u}{H^k(M)} + \Norm{u}{L^2(M)})\;.
\end{equation}
\end{Theorem}

\Proof
For simplicity we first treat the case $k=1$.  Theorem
\ref{ThmIR5} shows $u\in H^1_{\loc}$, so the vector of first derivatives
$\partial u\in L^2_{\loc}$ itself is a weak solution of the system of
equations
\begin{equation}
\label{IR-17}
 \cL\partial u+\partial(a^j)\, \partial_j u
 = \partial f - \partial(b)\, u\;.
\end{equation}
Since $\partial a^j\in L^{{n^*}}$ and
\[
\Norm{\partial(b)\,u}{L^2} \le C \Norm{b}{W^{1,n^*}} \Norm{u}{H^1},
\]
so $\Norm{\partial(b)\,u}{L^2}$ is bounded, this system satisfies the
conditions of Theorem \ref{ThmIR5}, hence $u\in H^2_{\loc}$.  The
general induction step applies a similar argument: if the result is
established $\forall\ k\le K-1$, and if $\cL u=f$ with coefficient
conditions \bref{IR-15} with $k=K$, then $\partial u$ satisfies an
elliptic system \bref{IR-17} of the same form with coefficient
conditions \bref{IR-15} with $k=K-1$, so by induction $\partial u\in
H^{K}(M)$ and thus $u\in H^{K+1}(M)$ as required.  The estimates
\bref{IR-16}, \bref{IR-16b} follow easily by a similar argument and
Theorem \ref{ThmIR5}.  \QED

The coefficient conditions in Theorem \ref{ThmIR6}
are not optimal in most cases.  For
example, if $n=3$ then $b\in W^{1,2}$ suffices to show $u\in H^2$
(rather than $b\in W^{1,3}$). This follows by interpolation,
\begin{eqnarray*}
\Norm{\partial b\,u}{L^2} &\le& \Norm{\partial b}{L^2}\Norm{u}{L^\infty}
\\
&\le& \epsilon \Norm{u}{H^2}
+ C(\epsilon,\Norm{\partial b}{L^2})  \Norm{u}{L^2} \;,
\end{eqnarray*}
which shows that $\partial b\,u$  may be thought of as the sum of a small
second order operator, and a large bounded operator on $L^2$.  The
small operator term may be absorbed as a perturbation of $\cL$,
and the remainder contributes to the right hand side source term.

\section{Spectral Condition}
\label{secSC}

In this section we review conditions under which an operator will have
a complete set of eigenfunctions.  These conditions will be used in
\S\ref{secB} to analyse boundary conditions, and thus the case
of most interest concerns operators on a compact manifold without
boundary, and in particular the first order elliptic systems
considered in \S\ref{secIR}.  However, the main result, Theorem
\ref{ThmSC1}, is stated in slightly more generality, which could
be used to extend the eigenfunction representation to operators on
manifolds with boundary.


%


Let $H$ be a closed subspace of $W^{1,2}(Y)$, with the induced norm,
where $Y$ is a compact manifold perhaps with boundary, and as in
\S\ref{secIR}, it is understood that these spaces refer to
sections of a (real) vector bundle $E$ over $Y$.


The abstract spectral theorem for the  map $A:H\to L^{2}(Y)$
uses the following conditions:
\begin{enumerate}
\item[($\cC0$)] $A:H\to L^{2}(Y)$ is linear and bounded 
in the  $W^{1,2}$ topology on $H$.
\item[($\cC1$)] The G\aa rding inequality holds: there exists a
  constant $C$ such that for all $ \psi \in H $ we have
  \begin{eqnarray}
    \Norm{\psi}{H}^2 & \le & C \int_Y \left(\langle A\psi,
    A\psi\rangle  +  \langle \psi,
    \psi\rangle \right)\, \dvYo \ .
\label{eq:garding}
  \end{eqnarray}
\item[($\cC2$)] Weak solutions are strong solutions (``elliptic
  regularity''): If $\phi\in
  L^2(Y)$  satisfies
  \begin{equation}
    \label{eq:weak}
     \int_Y \langle A\psi,\phi \rangle \, \dvYo
     =0 \ ,\qquad \forall\, \psi \in H ,
  \end{equation}
then $\phi\in H$.
\item[($\cC3$)] $A$ is symmetric:
  \begin{equation}
    \label{eq:symmetry}
    \forall \ \phi,\psi \in H \qquad \int_Y \langle A\phi,
    \psi\rangle \, \dvYo =\int_Y \langle \phi,
    A\psi\rangle \, \dvYo \ .
  \end{equation}
\item[($\cC4$)] density:
  \begin{equation}
    \label{eq:1.40}
    \mbox{$H$ is dense in $L^2(Y)$. }
  \end{equation}
\end{enumerate}
Note that in the case $\partial Y\ne \emptyset$,
the space $H$ must incorporate boundary conditions, and these will
play an important role in verifying $(\cC2)$, as will be seen in
\S\ref{secB}.

The main result of this section is the following:
\begin{Theorem}
  \label{ThmSC1} Under the conditions $(\cC0)$--$(\cC4)$, there exists
  a countable orthonormal basis of $L^2$ consisting of eigenfunctions
  of $A$, with eigenvalues all real and having no accumulation point
  in $\bR$.
\end{Theorem}
\proof Let $\X \subset H$ be the kernel of $A$; it is a standard fact
that $\X $ is finite dimensional when the G\aa rding inequality holds
--- we give the proof for completeness. Let $\{\psi_i\}_{i=1}^I$,
$I\le \infty$, be an $L^2$--orthonormal basis of $\X $, the
equation $A\psi_i = 0$ together with \eq{eq:garding} shows that
$\{\psi_i\}_{i=1}^I$ is bounded in $W^{1,2}$. The
Rellich theorem \cite[Theorem
7.22]{GilbargTrudinger77} implies that from the sequence $\psi_i$
we can extract a subsequence $\psi_{i_j}$ converging strongly in
$L^2$, weakly in $W^{1,2}$.
The G\aa rding inequality \eq{eq:garding} with $\psi$
replaced by $\psi_{i_j}-\psi_{i_k}$ shows that $\psi_{i_j}$ is Cauchy
in $W^{1,2}$, hence converges in norm to some $\psi\in W^{1,2}$. By
continuity of $A$, condition $(\cC0)$, we have $A\psi=0$, by
continuity of $L^2$ norm on $W^{1,2}$ it holds that
$\|\psi\|_{L^2}=1$, and it easily follows that
$\psi\in\{\psi_i\}_{i=1}^I$. We have thus shown that
$\{\psi_i\}_{i=1}^I$ is compact, which yields $I< \infty$, as desired.

Let now
$$\hH = \{\psi \in H:\forall\,  \phi \in \X \quad \int_ Y \langle \phi,
    \psi\rangle \, \dvYo =0 \}\ .
    $$
For $\phi\in L^2$ the map $H\ni \psi\to \int_ Y \langle \phi,
\psi\rangle \, \dvYo\in \bR$ is continuous in the $L^2$ topology (and
therefore also in the $W^{1,2}$ topology), thus $\hH$ is closed (being
an intersection of closed spaces), and hence a Banach space. We note
the following:
    \begin{Lemma}
      \label{LemSC2}There exists a constant $C$ such that
      \begin{equation}
        \label{eq:in1}
        \forall\, \psi \in \hH \qquad \|\psi\|_{L^2}\le C
        \|A\psi\|_{L^2}\ .
      \end{equation}
    \end{Lemma}
\proof
Suppose that this is not the case, then there exists a sequence
$\psi_n\in\hH$ such that
\begin{equation}
  \label{eq:in2}
  \|\psi_n\|_{L^2}\ge n
        \|A\psi_n\|_{L^2}\ \;.
\end{equation}
Rescaling $\psi_n$ if necessary we can without loss of generality
assume that $\|\psi_n\|_{L^2}=1$. The inequality \eq{eq:garding} shows
that $\psi_n$ is bounded in $W^{1,2}$ norm. By the Rellich
theorem \cite[Theorem 7.22]{GilbargTrudinger77} we can extract a
subsequence, still denoted $\psi_n$, converging to a $\psi_\infty\in
\hH$, weakly in $W^{1,2}$ and strongly in ${L^2}$. Equation
\eq{eq:in2} shows that the sequence $A\psi_n$ converges to zero in
$L^2$, and \eq{eq:garding} with $\psi$ replaced with $\psi_n-\psi_m$
shows that $\psi_n$ is Cauchy in the $W^{1,2}$ norm. Continuity of $A$
and Equation
\eq{eq:in2} imply that $A\psi_\infty = 0$, and since $A$ has no kernel
on $\hH$ we obtain $\psi_\infty = 0$, which contradicts
$\|\psi_\infty\|_{L^2}=1$, and the lemma follows. \qed

Returning to the proof of Theorem \ref{ThmSC1}, define $\Y$ to be the
image of $\hH$ under $A$. Then $\Y$ is a closed subspace of $L^2$,
which can be seen as follows: Let $\psi_i$ be any sequence in $\hH$
such that the sequence $\chi_i\equiv A\psi_i$ converges in $L^2$ to
$\chi_\infty \in L^2$. The inequality \eq{eq:in1} shows that $\psi_i$
is Cauchy in $L^2$, which together with the G\aa rding inequality
shows that $\psi_i$ is Cauchy in the $W^{1,2}$ norm. As $\hH$ is
closed, it follows that there exists $\psi_\infty\in \hH$ such that
$\psi_i$ converges to $\psi_\infty$ in the $W^{1,2}$ norm, and the
equality $\chi_\infty= A\psi_\infty$ follows from continuity of $A$.

Let $\phi\in L^2 $ be any element of $\Y^\perp$, the $L^2$ orthogonal of
$\Y$; by definition we have
$$
 \forall\, \psi \in H \qquad \int_Y \langle \phi,
    A\psi\rangle \, \dvYo =0 \ .
$$
The hypothesis $(\cC2)$ of elliptic regularity implies that $\phi\in H
$, so we can use the symmetry of $A$ to conclude
$$
 \forall\, \psi \in H \qquad \int_Y \langle A\phi,
    \psi\rangle \, \dvYo =0 \ .
$$
Density of $H$ in $L^2$ implies $A\phi=0$, thus
\begin{equation}
  \label{eq:perp}
\Y^\perp = \X  
\ \;.
\end{equation}
Define
$\hA:\hH\to \Y$ by $\hA\psi = A\psi$. By the definition of all the
objects involved the map $\hA$ is continuous, surjective and
injective, hence bijective.  Let $\hK:\Y\to \hH$ denote its inverse,
then $\hK$ is continuous by the open mapping theorem. Let $i$ be the
embedding of $ W^{1,2}(Y)$ into $L^2(Y)$; we have $i(\hH) \subset
\X ^\perp$ which coincides with $\Y $ by \eq{eq:perp}. 
It follows that for all $\chi\in \Y$ we have $i\circ \hK (\chi)\in
\Y$, so that $i\circ \hK $ defines a map of $\Y$ into $\Y$, which we
will denote by $K$. Now $\hK$ is continuous and $i$ compact, which
implies compactness of $K$.

We note that $\Y$ is a closed subset of the Hilbert space $L^2$, hence
a Hilbert space with respect to the induced scalar product. The
operator $K$ is self--adjoint with respect to this scalar product,
which can be seen as follows: let $\psi_a=K\phi_a$, $\phi_a\in \Y$,
$a=1,2$, thus $\psi_a\in H$ and $A\psi_a=\phi_a$. We then have
\[
\int_Y \langle \phi_1,
    K\phi_2\rangle \, \dvYo =\int_Y \langle A\psi_1,
    \psi_2\rangle \, \dvYo =\int_Y \langle \psi_1,
    A\psi_2\rangle \, \dvYo=
\int_Y \langle K\phi_1,
    \phi_2\rangle \, \dvYo
\ ,
\]
as desired.
By the spectral theorem for compact self adjoint operators
\cite{weidmann} there  exists a countable $L^2$--orthonormal basis of $\Y$
consisting of eigenfunctions of $K$:
\[
K\phi_\alpha=\mu_\alpha\phi_\alpha\ ,
\]
with eigenvalues $\mu_{\alpha}$ accumulating only at $0$.  Since $K$
is invertible we have $\mu_\alpha\ne 0$, hence
\[
 A\phi_\alpha=\lambda_\alpha
\phi_\alpha\ , \qquad \lambda_\alpha=\mu_\alpha^{-1}\ .
\]
The required basis of $L^2$ is obtained by completing $\{\phi_\alpha\}$
with any $L^2$--orthonormal basis of the finite dimensional kernel of
$A$. \QED

\begin{Definition}
\label{DefSC}
$A$ is said to satisfy the \emph{spectral condition} if $A$ is an
operator on $C^\infty$ sections of $E$ over $Y$ which is symmetric
with respect to the $L^2$ integration pairing with measure $dv_Y$ and
inner product $\langle\cdot,\cdot\rangle$, and there is a countable
orthonormal basis $\{\phi_\alpha\}_{\alpha\in\Lambda}$ of
$L^2\Gamma(E)$ consisting of eigenfunctions,
\begin{equation}
\label{SC-8}
A\phi_\alpha = \lambda_\alpha\phi_\alpha,\quad \alpha\in \Lambda\;,
\end{equation}
such that the eigenvalues $\lambda_\alpha\in\bR$, counted as always
with multiplicity, have no accumulation point in $\bR$.
\end{Definition}

\begin{Corollary}
\label{CorSC3}
Suppose $Y$ is a compact manifold without boundary
and $A:H^1(Y)\to L^2(Y)$ is an elliptic system between sections of the
bundles $E$, $F$, which satisfies the conditions
(\ref{IR-2},\ref{IR-3}) of Theorem \ref{ThmIR5}.
If $A=A^\dagger$ is formally
self-adjoint (see \bref{IR-5}), then $A$ satisfies the spectral
condition, Definition \ref{DefSC}.
\end{Corollary}

\Proof Take $H=H^1(Y)$.  Condition ($\cC0$) follows from the
coefficient bounds \bref{IR-2} and the inequality \bref{IR-1c}, and
condition ($\cC1$) is conclusion \bref{IR-14} of Theorem \ref{ThmIR5},
which also provides condition ($\cC2$).  Finally, ($\cC3$) follows
from the definition \bref{IR-5} of the $L^2$-adjoint $A^\dagger$,
since integration by parts is permitted in $H$, and ($\cC4$) is
standard.  The conclusions now follow from Theorem \ref{ThmSC1}.  \QED

\begin{Corollary}
\label{CorSC4}
Suppose $Y$ is a compact manifold without boundary and $A:H^1(Y)\to
L^2(Y)$ is an elliptic system between sections of the bundles $E$,
$F$, which satisfies the conditions of Theorem \ref{ThmIR5}.  There
are bases $\phi_\alpha\in L^2(E)$, $\psi_\alpha\in L^2(F)$, $\alpha\in
\Lambda$, with real numbers $\lambda_\alpha$ having no accumulation
point in $\bR$, which satisfy
\begin{equation}
\label{SC-9}
A\phi_\alpha  =\lambda_\alpha \psi_\alpha,\quad
A^\dagger\psi_\alpha  =\lambda_\alpha \phi_\alpha\;.
\end{equation}
The fields $\phi_\alpha$, $\psi_\alpha$ are all $H^1(Y)$.
\end{Corollary}
\Proof
This follows directly by applying Corollary \ref{CorSC3} to the
formally self-adjoint operator
\begin{equation}
\label{SC-10}
\mathbb{A} = \left[\begin{array}{cc} 0 & A^\dagger \\
                                     A & 0 \end{array} \right]\;,
\end{equation}
which acts between sections of the bundle $E\oplus F$.
\QED

\section{Boundary Regularity}
\label{secB}

In this section we introduce a broad class of boundary conditions which are
elliptic in the sense that the Weak-Strong property $(\cC2)$ can be
established, at least for solutions supported near the boundary.  When
combined with the interior regularity results of \S\ref{secIR}, this
will give the Weak-Strong property for compact manifolds with boundary
(\S\ref{secC}), and for a large class of noncompact manifolds with
compact boundary (\S\ref{secD}).  The main result is the boundary
regularity Theorem \ref{ThmB5}, and the primary ingredient in the arguments
is the energy identity \bref{B1:est} cf.~\cite{APS75} and \bref{A:Lich}.


We consider operators which may be written abstractly in the form
\begin{equation}
\label{B:L}
L = L_{0}+B = \partial_{x} + A + B\;,
\end{equation}
acting on sections of a (real) vector bundle $E$ over $\YI$, where $Y$
is a compact manifold without boundary\footnote{Although many of the
arguments of this section may be extended to allow $\partial
Y\ne\emptyset$, this would introduce technical complications which are
not relevant to the applications we have in mind.}, and for some
constant $\delta>0$,
\[  I=[0,\delta]. \]

Let $E\rY = i^*E$ be the pullback bundle over $Y$, where $i:Y\to\YI$,
$y\mapsto (y,0)$.  We assume that $A$ is an operator on sections of
$E\rY$ which is formally self-adjoint with respect to the pairing
defined by integration over $Y$ with the measure $dv_Y$ and the real
inner product $\langle\cdot,\cdot\rangle$ on the fibres of $E\rY$.
The operator $A$ and the inner product extend naturally to act on
sections of $E$ over $\YI$, and we likewise extend the definition of
the integration pairing by using the product measure $dv_{Y}dx$ on
$\YI$. Thus, $A$ is $x$--independent, but we allow $B$ to depend upon $x$.

We assume also that $A$ satisfies the spectral condition, Definition
\ref{DefSC}, so there is a countable index set $\Lambda$ and an
orthonormal basis $\left\{\phi_{\alpha}\right\}_{\alpha\in\Lambda}$ of
$L^{2}(E\rY)$ consisting of eigenfunctions
\begin{equation}
  A\phi_{\alpha} = \lambda_{\alpha}\phi_{\alpha},\qquad
   \alpha\in\Lambda,
  \label{B:fa}
\end{equation}
with eigenvalues $\lambda_{\alpha}\in\bbR$ having no accumulation
points in $\bbR$.

A formally self-adjoint first order elliptic operator with the
coefficient conditions of Theorem \ref{ThmIR5}, will satisfy these
conditions, by Corollary \ref{CorSC3}.


Although we have in mind primarily the case where $A,B$ are first
order differential (Dirac-type) operators, the results here will be
presented in an abstract form, because they could be applied more
widely.  For example, $A=-\Delta_{Y}$ will also satisfy the spectral
conditions, so the boundary regularity result Theorem \ref{ThmB5} may
also be applied to the heat equation.

We fix an eigenvalue cutoff parameter
$\kappa>0$, which is used to partition the index set $\Lambda$ into
\begin{eqnarray}
\label{B:Lam+}
\Lambda^{+} & = & \{\alpha\in\Lambda, \lambda_{\alpha}\ge\kappa\}, \\
\label{B:Lam-}
\Lambda^{-} & = & \{\alpha\in\Lambda, \lambda_{\alpha}\le -\kappa\} ,\\
\label{B:Lam0}
\Lambda^{0} & = & \{\alpha\in\Lambda, |\lambda_{\alpha}| < \kappa \}\;,
\end{eqnarray}
and we set $\Lambda' = \Lambda^{+}\cup\Lambda^{-}$.
It will be useful also to introduce a scale parameter
\begin{equation}
   \theta_{0} = \kappa^{-1} \max_{\alpha\in\Lambda^{0}}|\lambda_{\alpha}| < 1,
\label{B:kap0}
\end{equation}
which measures the size of the ``small'' eigenvalues. For example, we
could choose $\kappa$ to be the smallest nonzero eigenvalue, $\kappa =
\inf_{\lambda_\alpha\ne 0} |\lambda_{\alpha}|$, in which case
$\theta_0=0$.  Choosing $\kappa$ appropriately will lead to estimates
for $L$ which are uniform under perturbations of $A$ which create or
destroy small and zero eigenvalues.

The eigenfunction expansion $u=\sum_{\alpha\in\Lambda}
u_{\alpha}\phi_{\alpha}$ of $ u\in L^{2}(Y)$, where
\[  u_{\alpha}:= \oint_{Y} \langle  u,\phi_{\alpha}\rangle \, dv_{Y}, \]
leads to projection operators $P_{+},P_{-},P_{0},P'$, defined by
\begin{eqnarray}
\label{B:Ppm}
   P_{\pm} u &=& \sum_{\alpha\in\Lambda^{\pm}} u_{\alpha}\phi_{\alpha}, \\
\label{B:P0}
   P_{0} u &=& \sum_{\alpha\in\Lambda^{0}} u_{\alpha}\phi_{\alpha}\;,
\end{eqnarray}
and $P'=1-P_{0}=P_{+}+P_{-}$.

For $s\ge 0$, the Sobolev-type space $H^{s}_{*}(Y)$ is defined as the
completion of the space of smooth sections $C^{\infty}(Y)$, with
respect to the norm
\begin{equation}
\label{B:Hs}
   \Norm{u}{H^{s}_{*}(Y)}^{2} = \sum_{\alpha\in\Lambda'}
   |\lambda_{\alpha}|^{2s} |u_{\alpha}|^{2} +
   \kappa^{2s}\sum_{\alpha\in\Lambda^{0}} |u_{\alpha}|^{2}\;.
\end{equation}
The space $H^{1}_{*}(\YI)$ is likewise the completion of
$C^{\infty}(\YI)$ with respect to the norm
\begin{eqnarray}
  \nmi{u}^{2} &=& \int_{0}^{\delta} \left(
    \sum_{\alpha\in\Lambda}|u_{\alpha}'|^{2} +
    \sum_{\alpha\in\Lambda'} |\lambda_{\alpha}|^{2} |u_{\alpha}|^{2}
    + \kappa^{2}\sum_{\alpha\in\Lambda^{0}} |u_{\alpha}|^{2}
  \right)\, dx \nonumber \\
\label{B:H1}
&=& \int_{0}^{\delta}\oint_{Y} \left( |\partial_{x}u|^{2} +
  |Au|^{2} + \kappa^{2}|P_{0}u|^{2} \right)\,dv_{Y} dx\;,
\end{eqnarray}
where $'=\tfrac{d}{dx}$.  Of course in the typical case where $A$ is a
first order elliptic operator, these norms will be equivalent to the
usual Sobolev norms, defined using the Fourier transform.
Note that the normalization \bref{B:Hs} ensures that
the $L^{2}$ norm is controlled by the Sobolev norm:
\[
\Norm{u}{L^{2}(Y)} \le \kappa^{-s}\Norm{u}{H^{s}_{*}(Y)},\quad s\ge0.
\]

In addition, this formulation leads simply to a useful trace lemma,
stated in terms of the parameter $\ell$,
\begin{equation}
        \ell = \kappa\delta\;,
\label{B:ell}
\end{equation}
which measures the thickness of the boundary layer $Y\times[0,\delta]$
in units of $\kappa^{-1}$.

\begin{Lemma}\label{LemB0}
  The restriction map $r_{Y}:u\mapsto r_{Y}u = u(0,\cdot)$ from
  $C^{\infty}(Y)$ to $C^\infty(\YI)$, $I=[0,\delta]$, extends to
  a bounded linear map $r_{Y}:H^{1}_{*}(\YI)\to H^{1/2}_{*}(Y)$ satisfying
  \begin{equation}
          \nyh{r_{Y}u}^{2}  \le c_{1} \nmi{u}^{2}\;,
  \label{B0.1}
  \end{equation}
  where $c_{1}=c_{1}(\ell) = \ell^{-1}(1+\sqrt{1+\ell^{2}})$.
  The map $x\mapsto r_{Y,x}u$ (where $r_{Y,x}u = u(x,\cdot)$ is the
  restriction to $Y\times\{x\}$), is likewise bounded and continuous in $x$
  from $I$ to $H^{1/2}_{*}(Y)$.
  Moreover, $r_{Y}$ is surjective: there is an extension map
  $e_{Y}:H^{1/2}_{*}(Y)\to H^{1}_{*}(\YI)$ such that $r_{Y}e_{Y}u = u$
  for all $u \in H^{1/2}_{*}(Y)$ and $e_{Y}$ satisfies
  $r_{Y,\delta}e_Y(u) = 0$ and
\begin{equation}
  \nmi{e_{Y}u}^{2} \le \frac{2}{\sqrt{3}} \nyh{u}^{2}\ .
\label{B0.2}
\end{equation}
\end{Lemma}

\Proof
For $x\in [0,\infty)$ set $\chi(x)=\max(0,1-x)$.
For any $u\in
C^{\infty}(Y\times[0,\delta])$ with $u_{\alpha}(x)$,
$\alpha\in\Lambda$ denoting the spectral coefficients, and with
$\tilde{\chi}(x) = \chi(x/\delta)$, we find that
\begin{eqnarray*}
        |u_{\alpha}(0)|^{2} & = & -2\int_{0}^{\delta}
        \langle\tilde{\chi}u_{\alpha},
        \tfrac{d}{dx}(\tilde{\chi}u_{\alpha})\rangle \, dx
\\
         & \le & \int_{0}^{\delta}
         \left( (2\tilde{\chi}|\tilde{\chi}'| + \eta \tilde{\chi}^{2})
                   \, |u_{\alpha}|^{2}
               + \eta^{-1}\tilde{\chi}^{2}|u_{\alpha}'|^{2} \right) \, dx
\end{eqnarray*}
for any $\eta>0$.  For $\alpha\in\Lambda'$ we take $\eta=a
|\lambda_{\alpha}|$ and with an appropriate choice of $\chi$ we find
\[
   |\lambda_{\alpha}| |u_{\alpha}(0)|^{2}
   \le \int_{0}^{\delta} \left( a^{-1}|u_{\alpha}'|^{2}
     + |\lambda_{\alpha}|^{2}(a+2/\ell) |u_{\alpha}|^{2}\right)\, dx\;.
\]
Likewise for $\alpha\in \Lambda^{0}$, setting $\eta=a\kappa$ gives
\[
   \kappa |u_{\alpha}(0)|^{2}
   \le \int_{0}^{\delta} \left( a^{-1}|u_{\alpha}'|^{2}
     + \kappa^{2}(a+2/\ell) |u_{\alpha}|^{2}\right)\, dx\;.
\]
Choosing $a=\ell^{-1}(\sqrt{1+\ell^{2}}-1)$ and combining the two estimates
gives \bref{B0.1} for all $u\in C^{\infty}(Y\times[0,\delta])$.
\pcheck{14.III.98}
But this space is dense in $H^{1}_{*}(\YI)$ by definition, and it
follows easily that $r_{Y}u$ is defined and \bref{B0.1} is valid for
all $u\in H^{1}_{*}(\YI)$.

A very similar argument shows that for $x\in[0,\delta]$,
\[
   \nyh{r_{Y,x}u}^{2}  \le \ell^{-1}(2+\sqrt{4+\ell^{2}}) \nmi{u}^{2}.
\]

To establish continuity of $x\mapsto r_{Y,x}u$ as a map
$[0,\delta]\to H^{1/2}_{*}(Y)$,
note first that for any $v\in H^{1}_{*}(Y\times[x_{0},x_{1}])$, $x_{0}<x_{1}$,
the spectral coefficients $v_{\alpha}$ lie in $H^{1}([x_{0},x_{1}])$ and we
may compute:
\begin{eqnarray*}
        \left| |v_{\alpha}(x_{1})|^{2}- |v_{\alpha}(x_{0})|^{2}\right| & \le &
        2 \int_{x_{0}}^{x_{1}} \left| \langle v_{\alpha},v'_{\alpha}\rangle
        \right|\,dx
\\
         & \le & \eta_{\alpha}^{-1} \int_{x_{0}}^{x_{1}} \left( |v'_{\alpha}|^{2}+
         \eta_{\alpha}^{2} |v_{\alpha}|^{2}\right)\,dx\;.
\end{eqnarray*}
Choosing $\eta_{\alpha}=|\lambda_{\alpha}|$ for $\alpha\in\Lambda'$ and
$\eta_{\alpha}=\kappa$ for $\alpha\in\Lambda^{0}$ and summing gives
\begin{eqnarray}
        \left| \nyh{r_{Y,x_{1}}v}^{2} -  \nyh{r_{Y,x_{0}}v}^{2}\right|& \le &
        \int_{x_{0}}^{x_{1}} \biggl(\sum_{\alpha\in\Lambda}|u'_{\alpha}|^{2}
        +\sum_{\alpha\in\Lambda'}|\lambda_{\alpha}|^{2}|u_{\alpha}|^{2}
        +\sum_{\alpha\in\Lambda^{0}} \kappa^{2}|u_{\alpha}|^{2} \biggr)\,dx
\nonumber
\\
\label{B0.3}
         & \le & \Norm{v}{H^{1}_{*}(Y\times[x_{0},x_{1}])}\;.
\end{eqnarray}
Given $u\in H^{1}_{*}(\YI)$ and $\bar{x}\in[0,\delta)$,
$\epsilon\in(0,(\delta-\bar{x})/2)$, we set
\[
  v(x) = u(\bar{x}+\epsilon+x) - u(\bar{x}+\epsilon-x),
\]
where $x\in[-\epsilon,\epsilon]$ and the $Y$-dependence of $u,v$
is understood. Applying \bref{B0.3} with $x_{0}=0$,
$x_{1}=\epsilon$, gives $v(x_{0})=0$, $v(x_{1}) =
u(\bar{x}+2\epsilon)-u(\bar{x})$ and
\begin{eqnarray*}
     \nyh{r_{Y,\bar{x}+2\epsilon}u - r_{Y,\bar{x}}u}^{2}
         & \le & \Norm{v}{H^{1}_{*}(Y\times[0,\epsilon])}^{2}
\\
         & \le & 2 \Norm{u}{H^{1}_{*}(Y\times[\bar{x},\bar{x}+2\epsilon])}^{2}
\\
         & = & o(1)\ \ \mathrm{as}\ \ \epsilon\searrow0\;.
\end{eqnarray*}
This establishes continuity from the right, and left continuity follows
similarly.

To see that $r_{Y}$ is surjective, we construct an extension map
$e_{Y}:H^{1/2}_*(Y)\to H^1_*(\YI)$, such that $r_Y\circ e_Y = Id$.
For any $u\in H^{1/2}_{*}(Y)$ let $\{u^{k}\}$ be an approximating
Cauchy sequence of smooth fields with spectral coefficients
$u^{k}_{\alpha}$ and consider the sequence $\{\tilde{u}^{k}\}$ defined
by
\[
\tilde{u}^{k}(x,y) =  \sum_{\alpha\in\Lambda}
u^{k}_{\alpha}\phi_{\alpha}(y) \chi(x/\eta_{\alpha}),
\]
where $\eta_{\alpha} = \sqrt{3}/|\lambda_{\alpha}|$ for
$\alpha\in\Lambda'$ and $\sqrt{3}/\kappa$ for $\alpha\in
\Lambda^{0}$, so
\begin{eqnarray*}
  \nmi{\tilde{u}^{k}}^{2} & \le &
  \int_{0}^{\delta}\sum_{\alpha\in\Lambda'}|\lambda_{\alpha}|^{2}
  |u^{k}_{\alpha}|^{2}
    \Bigl( \chi^{2}(x/\eta_{\alpha}) +
    \tfrac{1}{3}\chi'^{2}(x/\eta_{\alpha})\Bigr) \, dx
  \\ && {}
  + \int_{0}^{\delta}\sum_{\alpha\in\Lambda^{0}}\kappa^{2}|u^{k}_{\alpha}|^{2}
  \Bigl(\chi^{2}(x/\eta_{\alpha}) +
    \tfrac{1}{3}\chi'^{2}(x/\eta_{\alpha})\Bigr)\, dx\;.
\end{eqnarray*}
Using the bounds $\int_{0}^{\infty}\chi^{2}(x)dx\le 1/3$,
$\int_{0}^{\infty}\chi'^{2}(x)dx\le1$, and noting that
\[
  \int_{0}^{\delta} \psi^{2}(x/\eta)dx \le \eta\int_{0}^{\infty}\psi^{2}(x)dx
\]
for any $\psi$, we have
\begin{eqnarray*}
  \nmi{\tilde{u}^{k}}^{2} & \le &  \frac{2}{\sqrt{3}}\Big\{
  \sum_{\alpha\in\Lambda'}|\lambda_{\alpha}| |u^{k}_{\alpha}|^{2}
    +\sum_{\alpha\in\Lambda_{0}}\kappa |u^{k}_{\alpha}|^{2} \Big\}
  \\ &
  \le & \frac{2}{\sqrt{3}} \nyh{\tilde{u}^{k}}^{2}  \;.
\end{eqnarray*}
Hence the sequence $\{\tilde{u}^{k}\}$ is uniformly bounded, and a
similar argument shows that it is also Cauchy, with limit $\tilde{u}
=e_{Y}u \in H^{1}_{*}(\YI)$.  It follows easily that the sequence has
boundary values converging to $u \in H^{1/2}_{*}(Y)$, and $u,e_{Y}u$
satisfy the bound \bref{B0.2}.\pcheck{14.III.98} \QED


The next result relates $H^{1}$ estimates to boundary conditions, and
is the key to understanding the nature of ellipticity for boundary
data.  It may be considered as a generalization either of the
integration by parts formula for the Dirac operator \bref{A:Lich}, or
of the estimate underlying the analysis in \cite{APS75}.  We use
$u(0), u(\delta)$ to denote the restrictions $r_Yu =
u\big|_{Y\times\{0\}}$, $r_{Y,\delta}u = u\big|_{Y\times\{\delta\}}$
respectively.

\begin{Lemma} \label{LemB1}
  Suppose $f\in L^{2}(\YI)$ and $u\in H^{1}_{*}(\YI)$ satisfies
  $L_{0}u=f$, then
\begin{eqnarray}
  \nmi{u}^{2} &\le&  \Norm{(1-P_{0})f}{L^{2}(\YI)}^{2} +
  (1+\theta_{0}^{2}) \Norm{P_{0}f}{L^{2}(\YI)}^{2}
\nonumber \\[3pt] &&{}
  + 3 \kappa^{2}\Norm{P_{0}u}{L^{2}(\YI)}^{2}
\nonumber \\ &&{} +
  \sum_{\alpha\in\Lambda^{+}}\lambda_{\alpha}
  \left(|u_{\alpha}(0)|^{2}-|u_{\alpha}(\delta)|^{2}\right)
\nonumber \\
  &&{} + \sum_{\alpha\in\Lambda^{-}}|\lambda_{\alpha}|
  \left(|u_{\alpha}(\delta)|^{2}-|u_{\alpha}(0)|^{2}\right)
\label{B1:est}
\end{eqnarray}
\end{Lemma}

%


\Proof The coefficient functions $u_{\alpha}(x)$ are measurable and,
by Fubini's theorem, square-integrable over $[0,\delta]$.  Testing the
weak formulation with $\phi(x,y)=\chi(x)\phi_{\alpha}(y)$ where
$\chi\in C^{\infty}_{c}((0,\delta))$ shows that $u_{\alpha}$ satisfies
\begin{equation} \label{ode}
    \int_{0}^{\delta} \left( -u_{\alpha}\chi' + u_{\alpha}
  \lambda_{\alpha}\chi       - f_{\alpha}\chi\right) \,dx =0
\end{equation}
for all $\chi\in C^{\infty}_{c}((0,\delta))$.  Because $u\in
H^{1}_{*}(\YI)$, the spectral coefficient $u_{\alpha}(x)$ is
differentiable for $a.e.\ x\in[0,\delta]$, with $u_{\alpha}'$
square-integrable, and \eq{ode} shows that it satisfies the ordinary
differential equation $u_{\alpha}'(x) + \lambda_{\alpha}u_{\alpha}(x)
= f_{\alpha}(x)$.  The trace lemma also shows that the restrictions
$u_{\alpha}(0)$, $u_{\alpha}(\delta)$ are well defined.  {}From the
ODE we derive the fundamental identity
\begin{eqnarray}
\nonumber
  \int_{0}^{\delta}|f_{\alpha}|^{2}\,dx & = & \int_{0}^{\delta} |u_{\alpha}' +
  \lambda_{\alpha}u_{\alpha}|^{2}\,dx
\nonumber
\\ & = &
  \int_{0}^{\delta}  \left(|u_{\alpha}'|^{2}
    +\lambda_{\alpha}^{2}|u_{\alpha}|^{2}\right)dx
\nonumber
\\ &  &
     {}+\lambda_{\alpha}(|u_{\alpha}(\delta)|^{2} - |u_{\alpha}(0)|^{2}).
\label{B:est}
\end{eqnarray}
Summing over $\alpha\in\Lambda^{+}\cup\Lambda^{-}$ and noting that the boundary
restrictions $u(0)$, $u(\delta)$ are in $H^{1/2}_{*}(Y)$ by Lemma
\ref{LemB0} since $u\in H^{1}_{*}(\YI)$ by assumption, we find
\begin{eqnarray}
  \nmi{P'u}^{2} &=&
  \int_{0}^{\delta}  \sum_{\alpha\in\Lambda^{+}\cup\Lambda^{-}}
  (|u'_{a}|^{2}+|\lambda_{\alpha}|^{2}|u_{\alpha}|^{2})\,dx
\nonumber \\
  & = & \int_{0}^{\delta}\oint_{Y}|P'f|^{2}dv_{Y}dx
\nonumber \\ &&{} +
  \sum_{\alpha\in\Lambda^{+}\cup\Lambda^{-}}\lambda_{\alpha}
  \left(|u_{\alpha}(0)|^{2}-|u_{\alpha}(\delta)|^{2}\right) \;.
\label{B:Hest2}
\end{eqnarray}
For $\alpha\in\Lambda^{0}$ we use $u'_{\alpha} =
f_{\alpha}-\lambda_{\alpha}u_{\alpha}$ to estimate
\begin{eqnarray}
  \nmi{P_{0}u}^{2} &=&
  \int_{0}^{\delta}   \sum_{\alpha\in\Lambda^{0}}
  (|u'_{a}|^{2}+\kappa^{2}|u_{\alpha}|^{2})\,dx
\nonumber \\
  & \le & \int_{0}^{\delta} \sum_{\alpha\in\Lambda^{0}}
     \left(  (1+\varepsilon) |f_{\alpha}|^{2} +
     (\kappa^{2}+(1+\varepsilon^{-1})\kappa_{0}^{2})|u_{\alpha}|^{2} \right)
\nonumber \\
  &\le&   3\kappa^{2}\Norm{P_{0}u}{L^{2}(\YI)}^{2} +
  (1+\theta_{0}^{2}) \Norm{P_{0}f}{L^{2}(\YI)}^{2}\;,
\label{B:Hest3}
\end{eqnarray}
having chosen $\varepsilon = \theta_0^2$.  Combining \bref{B:Hest2}
and \bref{B:Hest3} gives \bref{B1:est}. \pcheck{14.III.98} \QED


\begin{remark} The above proof could be generalized to allow $u\in
L^{2}$ and to show then that $u$ is in $H^{1}_{\loc}$, but this
refinement is unnecessary as we soon will show a more general
regularity theorem.  Working with $u\in H^{1}$ allows us to use the
boundary terms with impunity---a freedom that is not possible with
weak solutions at this stage.
\end{remark}

The fundamental estimate \bref{B1:est} shows that in order to obtain a
useful \emph{a priori} elliptic estimate for a general solution of
$L_{0}u=f$, it is necessary to impose boundary conditions which
control $P_{+}u(0)$ (and $P_{-}u(\delta)$).  Motivated by the examples
of the spectral and pointwise boundary conditions for the Dirac
equation (see \S\ref{secA}), we introduce a class of boundary
conditions which allow us to exploit the ``good'' terms in $P_{-}u(0)$
in \bref{B1:est} to provide the required control.  The effect of the
parameterization below is to describe the class of admissible boundary
data as graphs over the complementary subspace of ``good'' data
$(1-P_{+})H^{1/2}_{*}(Y)$.  The first justification of this approach
is the following existence result and its corresponding elliptic
estimate \bref{B2:est}.

\begin{Lemma}\label{LemB2}
  Let $P=P_{\Lambda^{+}\cup\hat{\Lambda}}$ be the spectral projection
  determined by $\Lambda^{+}$ and some subset $\hat{\Lambda}\subset
  \Lambda^{0}$ of the set of small eigenvalues.  Let $\sigma\in
  PH^{1/2}_{*}(Y)$ and $f\in L^{2}(\YI)$ be given, and suppose
  $K:(1-P)H^{1/2}_{*}(Y)\to PH^{1/2}_{*}(Y)$ is a continuous linear
  operator, so there is a constant $k\ge0$ such that for all $w\in
  H^{1/2}_{*}(Y)$,
  \begin{equation}
  \label{B2:K}
    \nyh{K(1-P)w} \le k\nyh{(1-P)w}.
  \end{equation}
  Then there exists a
  solution $u\in H^{1}_{*}(\YI)$ to the boundary value problem
  \begin{eqnarray} \label{B2:a}
      L_{0}u & = & f
  \\ \label{B2:b} Pu(0) & = & \sigma + K(1-P)u(0)
  \\ \label{B2:c} (1-P)u(\delta) & = & 0.
  \end{eqnarray}
  Moreover, the solution $u$ satisfies the estimate
  \begin{equation}
     \nmi{u}^{2} \le c_{4} (\Norm{f}{L^{2}(\YI)}^{2} +
     \nyh{\sigma}^{2})
  \label{B2:est}
  \end{equation}
  where the constant
  $c_{4}$ depends on $\ell, \theta_{0}$ and $k$.
\end{Lemma}

%

\Proof The solution to an ordinary differential equation
$u'(x) +\lambda u = f$ may be written in
either of the two forms
\begin{equation}
  u(x) = \left\{
  \begin{array}{c}
  \displaystyle
        e^{-\lambda x}u(0) + \int_{0}^{x}e^{\lambda(s-x)}f(s)ds  \\
  \displaystyle
        e^{\lambda(\delta-x)}u(\delta) -
        \int_{x}^{\delta}e^{\lambda(s-x)}f(s)ds .
  \end{array}  \right.
\label{B2.0}
\end{equation}
Consider first the spectral coefficients $u_\alpha(x)$ for $\alpha\in
\Lambda^{-}\cup\hat{\Lambda}'$, where $\hat{\Lambda}' =
\Lambda^{0}\backslash \hat{\Lambda}$.  The boundary condition
\bref{B2:c} is achieved for  $u_{\alpha}$
by letting $u_{\alpha}(\delta)=0$, so we define
\begin{equation}
        u_{\alpha}(x) = -\int
        _{x}^{\delta}e^{\lambda_{\alpha}(s-x)}f_{\alpha}(s)\,ds,
        \qquad \alpha\in \Lambda^{-}\cup\hat{\Lambda}',
\label{B2.1}
\end{equation}
where $f_{\alpha}(x)$ is the spectral coefficient of $f$.
Note that $f_{\alpha}\in L^{2}(I)$, so the integral in \bref{B2.1} is
well-defined.
The identity \bref{B:est} and $u_{\alpha}(\delta)=0$ shows that
\begin{equation}
        \int_{0}^{\delta} (|u'_{\alpha}|^{2} +
        \lambda_{\alpha}^{2}|u_{\alpha}|^{2})\,dx =
        \lambda_{\alpha} |u_{\alpha}(0)|^{2}
        + \int_{0}^{\delta} |f_{\alpha}|^{2}\,dx\;, \quad \forall \alpha\in
        \Lambda^{-}\cup\hat{\Lambda}'\ .
\label{B2.2}
\end{equation}
It follows that for all $\alpha\in\Lambda^{-}$,
\begin{equation}
        |\lambda_{\alpha}| |u_{\alpha}(0)|^{2}
        \le \int_{0}^{\delta} |f_{\alpha}|^{2} dx
        - \int_{0}^{\delta}(|u'_{\alpha}|^{2} +
        \lambda_{\alpha}^{2}|u_{\alpha}|^{2})\,dx .
\label{B2.2b}
\end{equation}

{}To control the small eigenvalues $\alpha\in \hat{\Lambda}'$ we use an
elementary lemma, the proof of which is an exercise:
\begin{Lemma}\label{LemB2a}
  For any $f\in L^{2}([0,\delta])$ and $\lambda\in\bbR$,
  \begin{equation}
 \int_{0}^{\delta}\left(\int_{x}^{\delta}e^{\lambda(s-x)}f(s)ds\right)^{2}dx
  \le \left\{
  \begin{array}{ll}
        \half \delta^{2}e^{2\lambda\delta}\int_{0}^{\delta}f^{2}(x)dx &
        \lambda>0 \\[5pt]
        \half \delta^{2}\int_{0}^{\delta}f^{2}(x)dx & \lambda \le 0.
  \end{array} \right.
  \label{B2a.1}
  \end{equation}
\end{Lemma}

{}From \bref{B2.1} and Lemma \ref{LemB2a} it follows that for
$\alpha\in \hat{\Lambda}'$,
\[
  \int_{0}^{\delta} \kappa^{2}|u_{\alpha}|^{2}dx
  \le \half \ell^{2}e^{2\ell\theta_{0}}
  \int_{0}^{\delta}|f_{\alpha}|^{2}\,dx\;.
\]
Using $u'_\alpha=f_\alpha - \lambda_\alpha u_\alpha$ as in \bref{B:Hest3}
we obtain
\begin{equation}
        \int_{0}^{\delta} (|u'_{\alpha}|^{2}+\kappa^{2}|u_{\alpha}|^{2})dx \le
        (1+c_{2})\,\int^{\delta}_{0}|f_{\alpha}|^{2}dx\;,
\label{B2.3}
\end{equation}
where
\begin{equation}
        c_{2}=c_{2}(\ell,\theta_0) =
        \theta_0^2+\tfrac{3}{2}\ell^{2}e^{2\ell\theta_0}.
\label{B:c2}
\end{equation}
Combining \bref{B2.2} and \bref{B2.3} shows that
\[
   u^{(-)}(x,y) := \sum_{\alpha\in \Lambda^{-}\cup\hat{\Lambda}'}
   u_{\alpha}(x) \phi_{\alpha}(y)
\]
is a sum converging in $H^{1}_{*}(\YI)$, and $u^{(-)}$ satisfies
\begin{eqnarray}
        \nmi{u^{(-)}}^{2} & \le & \Norm{(1-P)f}{L^{2}(\YI)}^{2}
        + c_{2}\Norm{P_{\hat{\Lambda}'}f}{L^{2}(\YI)}^{2}
\nonumber \\ &&{}
        - \sum_{\alpha\in\Lambda^{-}}|\lambda_{\alpha}|\,|u_{\alpha}(0)|^{2}\;,
\label{B2.4}
\end{eqnarray}
where $1-P=P_{-}+P_{\hat{\Lambda}'}$.  For $\alpha\in\hat{\Lambda}'$
such that $\lambda_\alpha\ne0$, using \bref{B2.1} we find that
\begin{eqnarray}
\nonumber
   |u_{\alpha}(0)|^{2} & \le &
   \left(\int_{0}^{\delta}e^{\lambda_{\alpha}s}f_{\alpha}(s)ds\right)^{2}
\\
\nonumber &\le& \frac{1}{2\lambda_{\alpha}}(e^{2\lambda_{\alpha}\delta}-1)
   \int_{0}^{\delta} |f_{\alpha}|^{2} dx\;,
\end{eqnarray}
while if $\lambda_\alpha = 0$ then $|u_{\alpha}(0)|^{2} \le \delta
\int_{0}^{\delta} |f_{\alpha}|^{2} dx$.  Defining
\begin{equation}
        c_{3}\  =\  \ell e^{2\ell\theta_0},
\label{B:c3}
\end{equation}
it follows that
\[
\kappa |u_{\alpha}(0)|^{2}  \le  c_{3}\int_{0}^{\delta}
|f_{\alpha}|^{2} dx\;, \quad \forall \ \alpha\in \hat{\Lambda}'.
\]
Thus, combining with \bref{B2.4} we have
\begin{eqnarray}
   \lefteqn{ \nmi{u^{(-)}}^{2} + \nyh{u^{(-)}(0)}^{2} \hfill \qquad {}}
\nonumber \\ &\le&
    \Norm{(1-P)f}{L^{2}(\YI)}^{2}
    + (c_{2}+c_{3})\Norm{P_{\hat{\Lambda}'}f}{L^{2}(\YI)}^{2} .
\label{B2.5}
\end{eqnarray}

For $\alpha\in\Lambda^{+}\cup\hat{\Lambda}$ we use \bref{B2.0} to
define $u_{\alpha}$ by
\begin{equation}
   u_{\alpha}(x) = e^{-\lambda_{\alpha}x}(\sigma_{\alpha}+
   K_{\alpha}u^{(-)}(0) ) + \int_{0}^{x}
   e^{\lambda_{\alpha}(s-x)}f_{\alpha}(s)\,ds ,
\label{B2.6}
\end{equation}
where $\sigma_{\alpha}$, $K_{\alpha}u^{(-)}(0)$ denote the
$\phi_{\alpha}$ coefficients of $\sigma$ and $Ku^{(-)}(0)$
respectively.  Note in particular that \bref{B2.5} shows that
$u^{(-)}(0)\in H^{1/2}_{*}(Y)$, so $Ku^{(-)}(0)\in H^{1/2}_{*}(Y)$ by
the hypothesis \bref{B2:K}, hence the coefficients
$K_{\alpha}u^{(-)}(0)$ are well defined.


For $\alpha\in\Lambda^{+}$ we estimate using 
\eq{B:est} and \bref{B2:b}:
\begin{eqnarray}
\nonumber \int_{0}^{\delta} (|u'_{\alpha}|^{2}+\lambda_{\alpha}^{2})dx
        & \le &
        \lambda_{\alpha}|\sigma_{\alpha}+K_{\alpha}u^{(-)}(0)|^{2} +
        \int_{0}^{\delta}|f_{\alpha}|^{2} dx \\
\label{B2.7}
         & \le &
    2\lambda_{\alpha}(|\sigma_{\alpha}|^{2}+|K_{\alpha}u^{(-)}(0)|^{2})
    + \int_{0}^{\delta}|f_{\alpha}|^{2} dx  \;.
\end{eqnarray}
For $\alpha\in \hat{\Lambda}$ we estimate directly from \bref{B2.6}:
\begin{eqnarray}
\nonumber
\kappa^{2}      \int_{0}^{\delta}|u_{\alpha}|^{2}dx & \le &
  3\kappa^{2}(|\sigma_{\alpha}|^{2} + |K_{\alpha}u^{(-)}(0)|^{2}) \int_{0}^{\delta}
          e^{-2\lambda_{\alpha}x}dx
\\
\nonumber
         &  & {}+ 3 \kappa^{2}\int_{0}^{\delta}\left(
         \int_{0}^{x}e^{\lambda_{\alpha}(s-x)}f_{\alpha}(s)ds\right)^{2}dx
\\
\label{B2.8}
     & \le & 3 c_{3} \kappa (|\sigma_{\alpha}|^{2} + |K_{\alpha}u^{(-)}(0)|^{2})
         +\tfrac{3}{2}\ell c_{3} \int_{0}^{\delta}|f_{\alpha}|^{2}dx \;,
\end{eqnarray}
\pcheck{14.III.98}
where Lemma \ref{LemB2a} has been used to control the final term.
Using $u_{\alpha}' = f_{\alpha}-\lambda_{\alpha}u_{\alpha}$ to estimate
$|u_{\alpha}'|^{2}\le (1+\varepsilon)|f_{\alpha}|^{2}
                     +(1+\varepsilon^{-1})\lambda_{\alpha}^{2}|u_{\alpha}|^{2}$
with $\varepsilon = \theta_0^2$, \bref{B2.8} gives
for $\alpha\in\hat{\Lambda}$,
\begin{eqnarray}
\nonumber
        \int_{0}^{\delta}\left(
        |u'_{\alpha}|^{2}+\kappa^{2}|u_{\alpha}|^{2}\right) dx & \le &
           9c_{3}\kappa  (|\sigma_{\alpha}|^{2}
           + |K_{\alpha}u^{(-)}(0)|^{2})
\\
\label{B2.9}
         &  & {}+ (1+3c_{2}) \int_{0}^{\delta}|f_{\alpha}|^{2} dx \;.
\end{eqnarray}
Combining \bref{B2.7} and \bref{B2.9} we have
(setting $u^{(+)}=\sum_{\Lambda^{+}\cup\hat{\Lambda}}u_{\alpha}\phi_{\alpha}$)
\begin{eqnarray}
\nonumber
        \nmi{u^{(+)}}^{2} & \le & \Norm{Pf}{L^{2}(\YI)}^{2}
            + 3c_{2} \Norm{P_{\hat{\Lambda}}f}{L^{2}(\YI)}^{2}
\\
\nonumber
         &  & {} + 2 (\nyh{P_{+}\sigma}^{2} + \nyh{P_{+}Ku^{(-)}(0)}^{2})
\\
\label{B2.10}
         &  & {} + 9c_{3} (\nyh{P_{\hat{\Lambda}}\sigma}^{2}
               + \nyh{P_{\hat{\Lambda}}Ku^{(-)}(0)}^{2}) \;,
\end{eqnarray}
\pcheck{14.III.98}
where $P_{+}=P_{\Lambda^{+}}$, $P=P_{+}+P_{\hat{\Lambda}}$.  Since we
have already shown that $K_{\alpha}u^{(-)}(0)\in H^{1/2}_{*}(Y)$, all
terms on the right hand side of \bref{B2.10} are bounded, which shows
that $u^{(+)}$ is well-defined in $H^{1}_{*}(\YI)$.

With $u=u^{(+)}+u^{(-)}$, we add an appropriate multiple of
\bref{B2.5} to \bref{B2.10} to control the bad terms in $Ku^{(-)}(0)$
with the good term $u^{(-)}(0)$ of \bref{B2.5}.  This gives the
elliptic estimate \bref{B2:est}:
\begin{eqnarray*}
        \nmi{u}^{2} & \le & \nmi{u^{(+)}}^{2} + \max(1,2k^{2}\;,9c_{3}k^{2})
        \nmi{u^{(-)}}^{2}
\\
         & \le & c_{4}^2(\Norm{f}{L^{2}(\YI)}^{2} + \nyh{\sigma}^{2})\;,
\end{eqnarray*}
where $c_{4}=c_{4}(\ell,\theta_0,k)$ as required.  The
definitions \bref{B2.1}, \bref{B2.6} ensure $u$ is a solution
satisfying the boundary conditions \bref{B2:b}, \bref{B2:c}.  \QED

Explicitly, we may take $c_4^2 = 3c_2 +(k^2+1)(2+9c_3)$ in general, and
$c_4^2 = 2\max(1,k^2)$ if $\theta_0=0$.

%
%

The next result is the key to handling operators with coefficients
depending on $x$.  Recall that the operator norm $\Vert B\Vert_{op}$
of a linear map $B:X_{1}\to X_{2}$ between Banach spaces is the
smallest constant such that
\begin{equation}
  \Vert Bu\Vert_{X_{2}} \le \Vert B\Vert_{op} \Vert u\Vert_{X_{1}},
  \quad\forall u\in X_{1}\; .
\label{B:Bop}
\end{equation}

\begin{Lemma}\label{LemB3}
  Suppose $L_{0},A,f,\sigma,K$ are as in
  Lemma \ref{LemB2}, and suppose $B:H^{1}_{*}(\YI)\to L^{2}(\YI)$
  is a linear map satisfying
   \begin{equation} \label{B3:bnd}
        {c_{4}}\,\Norm{B}{op} < 1.
   \end{equation}
   Then there exists $u\in H^{1}_{*}(\YI)$ satisfying
   \[
   (L_{0}+B)u = f
   \]
   and the boundary conditions \bref{B2:b},\bref{B2:c}, such that
   \begin{equation}
         \nmi{u}^{2} \le
         \frac{c_{4}}{1-c_4\Norm{B}{op}}
                 (\Norm{f}{L^{2}(\YI)} + \nyh{\sigma}).
   \label{B3:est}
   \end{equation}
\end{Lemma}

\Proof Let $u^{(0)}\in H^{1}_{*}(\YI)$ be any function satisfying the
boundary conditions \bref{B2:b},\bref{B2:c}; the trace lemma \ref{LemB0}
ensures the existence of a suitable $u^{(0)}$.  Construct a
sequence $\{u^{(k)}\}\subset H^{1}_{*}(\YI)$ by solving the problems
\begin{eqnarray}
\label{B3:it}
   L_{0}u^{(k)} & = & f - Bu^{(k-1)} \\
\label{B3:b}
   Pu^{(k)}(0) &=& \sigma + K(1-P)u^{(k)}(0) \\
\label{B3:c}
   (1-P)u^{(k)}(\delta)& =&0,\quad n=1,2,\ldots \;.
\end{eqnarray}
Lemma \ref{LemB2} ensures this problem has a solution for every
$n\ge1$, and the difference $w^{(k)}=u^{(k)}-u^{(k-1)}$ satisfies
\begin{eqnarray*}
L_{0}w^{(k)}         & = &  - Bw^{(k-1)} \\
Pw^{(k)}(0)          & = & K(1-P)w^{(k)}(0) \\
(1-P)w^{(k)}(\delta) & = & 0\;.
\end{eqnarray*}
The elliptic estimate \bref{B2:est} gives
\[
\nmi{w^{(k)}}\le c_{4} \Vert
Bw^{(k-1)}\Vert_{L^{2}(\YI)} .
\]
If $\Vert B\Vert_{op}<1/c_{4}$ then the iteration is a contraction and
thus the sequence $u^{(k)}$ is Cauchy, converging to
${u}=\lim_{n\to\infty}u^{(k)}$ strongly in $H^{1}_{*}(\YI)$.  Taking
the limit of \bref{B3:it} shows that $(L_{0}+B){u}=f$, and boundedness
of the trace operator $r_{Y}$ shows that ${u}$ satisfies the boundary
conditions (\ref{B2:b})--(\ref{B2:c}).  The elliptic estimates
\bref{B2:est} satisfied by $u^{(k)}$ are preserved in the limit, so
$u$ satisfies
\[
        \nmi{u} \le  c_{4}(\Norm{f-Bu}{L^{2}(\YI)} + \nyh{\sigma}),
\]
from which \bref{B3:est} follows easily.
\pcheck{14.III.98}
\QED

Observe that the proof of Lemma \ref{LemB3} relies on just two
properties of the operator $L_0$; namely, the solvability of the
problem \bref{B3:it} with boundary conditions \bref{B3:b},
\bref{B3:c}, and the elliptic estimate \bref{B2:est}, which provides
the size bound \bref{B3:bnd} for the perturbation $B$. This suggest
that it should be possible to extend this existence result to more
general operators $L=L_0+B$, for which a strictly coercive estimate
such as \bref{B2:est} can be established.

Consider, for example, the case where $E$ has a complex structure
$J:E\to E$, $J^2=-1$, and $A$ is a \emph{normal} operator ($[A,A^*]=0$), so
$A=A_0+JA_1$ where $A_0,A_1$ are self-adjoint and commuting, $A_0$
satisfies the spectral condition, and both commute with $J$.  Then $A$
admits an eigenfunction basis
\[
A \phi_\alpha = (\lambda_\alpha+\mu_\alpha J)\phi_\alpha,
\ \ \lambda_\alpha, \mu_\alpha \in \bR, \ \forall \alpha\in\Lambda,
\]
and the results of this section extend with only minor modifications,
provided the eigenvalues $\lambda_\alpha, \mu_\alpha$ satisfy the
sectorial condition
\[
\sup_{\alpha\in\Lambda} |\mu_\alpha|/(1+|\lambda_\alpha|) < \infty.
\]

Before stating the main uniqueness theorem, a definition of
weak solution with boundary conditions is required.  Note that although
the definition is consistent with just $L^{2}$ boundary data, the
regularity theorem \ref{ThmB5} will require data in $H^{1/2}_{*}$.

\begin{Definition} \label{Def:weak}
  Suppose $L=\partial_{x}+A+B$ where $A$ satisfies the
  spectral conditions (Definition \ref{DefSC}) and
  $B:H^{1}_{*}(\YI)\to L^{2}(\YI)$ is a bounded linear operator
  for which there exists an $L^2$-adjoint $B^{\dagger}:  H^{1}_{*}(\YI)\to
  L^2(\YI)$ such that:
  \begin{equation}
    \int_{\YI} \langle Bu,v\rangle \,dv_Ydx
    = \int_{\YI} \langle u,B^{\dagger}v\rangle \,dv_Ydx\;,
    \quad \forall\ u,v\in H^{1}_{*}(\YI).
  \label{B:Bdag}
  \end{equation}
  Suppose further that $P=P_{+}+P_{\hat{\Lambda}}$
  (as in  Lemma \ref{LemB2}),
  that $K:(1-P)L^{2}(Y)\to PL^{2}(Y)$ is a bounded
  linear map with $L^2$-adjoint
  $K^{\dagger}:PL^{2}(Y)\to (1-P)L^{2}(Y)$,
  and let $\sigma \in PL^{2}(Y)$, $f\in L^{2}(\YI)$ be given.
  A \emph{weak solution} of the boundary value problem
\oldnote{$Pu(\delta)$ changed to $(1-P)u(\delta)$ in \eq{B:c}}
\begin{eqnarray}
  \label{B:Luf}
        Lu & = & f \\
  \label{B:b}
        Pu(0) & = & \sigma + K(1-P) u(0)  \\
   \label{B:c}
        (1-P)u(\delta) & = & 0
   \end{eqnarray}
   is a field $u\in L^{2}(\YI)$ satisfying
   (with $L^{\dagger}=-\partial_{x}+A+B^{\dagger}$)
   \begin{equation}
        \int_{\YI}\langle u, L^{\dag}\phi\rangle \, dv_Ydx=
        \int_{\YI}\langle f,\phi\rangle\, dv_Ydx
        + \oint_{Y} \langle \sigma,\phi(0)\rangle\, dv_Y\;,
   \label{B:wLuf}
   \end{equation}
        for all $\phi\in H^{1}_{*}(\YI)$ satisfying the
       \emph{adjoint boundary conditions}
   \begin{eqnarray}
   \label{B:adb}
        (1-P + K^{\dag}P)\phi(0) &=& 0 \\
   \label{B:adc}
        P\phi(\delta) &=& 0\;.
\end{eqnarray}
\end{Definition}
\oldnote{check $K$ bndd $\Leftrightarrow$ $K^{\dag}$ bndd? {\bf piotr}:
         this is correct}


The boundary values $\phi(0), \phi(1)$ both lie in $ H^{1/2}_{*}(Y)$ by the
trace lemma, so the adjoint boundary conditions are well-defined on the
space of test fields.  Since $C^{\infty}$ fields are dense in
$H^{1}_{*}(\YI)$ and in $H^{1/2}_{*}(Y)$, to verify the weak equation
\bref{B:wLuf} it suffices to test just with $C^{\infty}$ fields $\phi$;
however the uniqueness argument of Lemma \ref{LemB4} requires the use of an
$H^{1}_{*}$ test field.


The structure of the adjoint boundary condition \bref{B:adb} is
explained by the next lemma, which is applied with $v=u(0)-\sigma$ and
$H=L^2(Y)$.

\begin{Lemma} \label{LemB4a}
  If $H$ is a Hilbert space, $P:H\to H$ is an orthogonal projection
  and $K:\ker P\to \ran P$ is bounded, and if $v\in H$ satisfies
  \begin{equation}
          \langle v, \phi\rangle_{H} =0\quad \forall\ \phi\in\ker(1-P+K^{\dag}P)
  \label{B4a}
  \end{equation}
  (where $K^{\dag}$ is the adjoint of $K$ in $H$),
  then $v\in\ker(P-K(1-P))$.
\end{Lemma}
\Proof Since $\ker P=\ran(1-P) \perp \ran P$, it follows that $P$ is
  self-adjoint and there is an orthogonal decomposition $H=
  (1-P)H\oplus PH$.  Setting $\phi_{1}=(1-P)\phi$, $\phi_{2}=P\phi$,
  the condition $\phi=\phi_{1}+\phi_{2}\in\ker(1-P+K^{\dag}P)$ is
  equivalent to $\phi_{1}=-K^{\dag}\phi_{2}$, which exhibits
  $\ker(1-P+K^{\dag}P)$ as a graph over $PH$.  Similarly decomposing
  $v=v_{1}+v_{2}$, the condition $\langle v,\phi\rangle=0$ is
  equivalent to $\langle v_{2}-Kv_{1},\phi_{2}\rangle=0$.  Since this
  holds for all $\phi_{2}\in PH$, it follows that $v_{2}=Kv_{1}$, or
  equivalently, $v\in\ker(P-K(1-P))$.
\QED

In other words, $H=\ker(P-K(1-P)) \oplus \ker(1-P+K^\dag P)$ is an
orthogonal splitting of $H$, where $P-K(1-P)$, $1-P+K^\dag P$
are projections, which are not orthogonal in general.

Using Lemma \ref{LemB4a} we next show that an $H^1$ weak solution of the
boundary value problem (\ref{B:Luf})-(\ref{B:c}), in fact satisfies the
equation \bref{B:Luf} and the boundary conditions (\ref{B:b},\ref{B:c}) in
the strong sense:

\begin{Lemma}  \label{LemB4b}
  If $u\in H^1_*(\YI)$ is a weak solution of \bref{B:Luf} with the boundary
  conditions \bref{B:b}, \bref{B:c}, then $u$ satisfies the equation $Lu=f$
  in the sense of strong ($H^1$) derivatives, and the restrictions
  $u(0)=r_Y(u)$, $u(\delta) = r_{Y,\delta}u$ satisfy the boundary
  conditions \bref{B:b},\bref{B:c} in $L^{2}(Y)$.  Conversely, if $u\in
  H^1_*(\YI)$ is a strong solution of (\ref{B:Luf},\ref{B:b},\ref{B:c}),
  then $u$ is also a weak solution.
\end{Lemma}
\Proof
Integration by parts gives
\[
  \int_{\YI} \langle u,L^{\dag}\phi\rangle =
  \int_{\YI}\langle Lu,\phi\rangle
  + \oint_{Y}\langle u(0),\phi(0)\rangle
  - \oint_{Y}\langle u(\delta),\phi(\delta)\rangle
\]
for any $u,\phi\in H^{1}_{*}(\YI)$.  Testing $u$ with arbitrary
$\phi\in C^{\infty}_c(\YI)$ shows that a $H^{1}_{*}$ weak solution
satisfies $Lu=f$ in the sense of strong derivatives.
Comparing this formula with \bref{B:wLuf} shows also that
\[
        \oint_{Y}\langle u(0)-\sigma,\phi(0)\rangle  =0,\qquad
        \oint_{Y}\langle u(\delta),\phi(\delta)\rangle  =0,
\]
for all $\phi(0),\phi(\delta)$ satisfying the adjoint boundary conditions
\bref{B:adb}, \bref{B:adc}.  Since $H^{1/2}_*(Y)$ is dense in
$\ker(1-P+K^\dagger P)\subset L^2(Y)$, Lemma \ref{LemB4a} may be applied
with $v=u(0)-\sigma$ to show the boundary condition \bref{B:b} holds in
$L^2(Y)$, and \bref{B:c} follows similarly.

To show the converse,  integration by parts again  gives
\begin{eqnarray*}
\lefteqn{\int_{\YI} \langle u,L^{\dag}\phi\rangle - \langle Lu,\phi\rangle
\hspace*{2cm}}
\\
&=&
   \oint_{Y}\langle u(0),\phi(0)\rangle
  - \langle u(\delta),\phi(\delta)\rangle
\\
&=&
   \oint_Y \langle \sigma +(1+K)(1-P)u(0),\phi(0)\rangle
  -  \langle Pu(\delta),\phi(\delta)\rangle
\\
&=&
   \oint_Y \langle \sigma,\phi(0)\rangle
         + \langle u(0),(1-P+K^\dagger P)\phi(0)\rangle
         - \langle u(\delta),P\phi(\delta)\rangle \;,
\end{eqnarray*}
and the final two terms vanish by the adjoint boundary conditions
(\ref{B:adb},\ref{B:adc}).
\QED

By solving an adjoint problem, we now show that weak solutions of
(\ref{B:Luf})-(\ref{B:c}) are unique.

\begin{Lemma}  \label{LemB4}
  Let $u\in L^{2}(\YI)$ be a weak solution of the boundary value
  problem (\ref{B:Luf})-(\ref{B:c}), with $\sigma\in
  L^{2}(Y)$ and $f\in L^{2}(\YI)$.  Suppose that the operator
  $L=\partial_{x}+A+B$ satisfies the conditions of Lemma \ref{LemB3}
  and Definition \ref{Def:weak}, and the $L^2(\YI)$-adjoint
  $B^\dagger:H^1_*(\YI)\to  L^2(\YI)$ satisfies
  \begin{equation}  \label{B:Bdbnd}
   c_4  \Norm{B^\dagger}{op} < 1.
  \end{equation}
Suppose also that the boundary operators $K,K^\dagger$ of
  Definition \ref{Def:weak} satisfy
  \begin{eqnarray} \label{B:K1}
     \nyh{K(1-P)w} &\le&  k\,\nyh{(1-P)w}\ ,
  \\ \label{B:K2}
     \nyh{K^{\dagger}Pw} &\le& k\,\nyh{Pw} .
  \end{eqnarray}
  for some constant $k\ge0$ and all $w\in H^{1/2}_{*}(Y)$.
  Then $u$ is unique.
\end{Lemma}

\Proof It will suffice to show that any weak solution $\tilde{u}$ of
\bref{B:Luf}--\bref{B:c} with $\sigma=0$, $f=0$, must vanish.
Consider the adjoint problem $L^{\dag}\phi= \tilde{u}$ with boundary
conditions \bref{B:adb},\bref{B:adc}; writing $L^{\dag}\phi=
\tilde{u}$ as $(\partial_x-A-B^{\dag})\phi=-\tilde{u}$, we see that
$L^{\dag}$, $K^{\dag}$ satisfy the conditions required by Lemma
\ref{LemB3}, since interchanging $A\leftrightarrow -A$ means replacing
$P$ by $1-P$, and perhaps changing a finite number of eigenfunctions
in $\hat\Lambda$ (without modifying $\Lambda_0$). The elliptic
estimate \bref{B3:est} does not depend on $\hat\Lambda$. Thus, by
Lemma \ref{LemB3} there exists a solution $\phi\in H^{1}_{*}(\YI)$ of
this boundary value problem. By construction, $\phi$ satisfies the
boundary conditions required of test functions in \bref{B:wLuf}, so
testing $\tilde{u}$ in \bref{B:wLuf} with $\phi$ gives
\[
\int_{\YI}|\tilde{u}|^{2} = \int_{\YI} \langle
\tilde{u},L^{\dag}\phi\rangle = 0
\]
and thus $\tilde{u}=0$.  \QED

It is easy to check that \bref{B:K2} is equivalent to requiring
that $K:(1-P)H^{-1/2}_*(Y) \to PH^{-1/2}_*(Y)$ is bounded, with
constant $k$.

We now obtain the main result on boundary regularity of weak
solutions.  Note that although the definition of weak solution assumes
boundary data $\sigma\in L^{2}(Y)$ only, and uniqueness of weak
solutions holds also in this generality, this condition is
incompatible with regularity $u\in H^{1}_{*}(\YI)$, which would imply
(by simple restriction) that $\sigma\in H^{1/2}_{*}(Y)$.
However, some results for $L^2$ boundary conditions on domains with
uniformly Lipschitz boundary are known \cite{McIntoshEtal96,AxelssonEtal00},
so it is plausible that the results here could be extended.

\begin{Theorem} \label{ThmB5}
  Suppose $u\in L^{2}(\YI)$ is a weak solution of the boundary value
  problem (\ref{B:Luf})-(\ref{B:c}) with operator
  $L=\partial_{x}+A+B_0+B_1$, where $A$ satisfies the spectral conditions
  (Definition \ref{DefSC}), $B_0$ satisfies the size condition
  \bref{B3:bnd} with $L^2$-adjoint $B_0^\dagger$ satisfying
  \bref{B:Bdbnd}, and $B_1:L^{2}(\YI)\to L^{2}(\YI)$ is bounded.
  Further suppose the boundary operators $K,K^\dagger$ satisfy
  (\ref{B:K1},\ref{B:K2}), and $\sigma\in PH^{1/2}_{*}(Y)$.
  Then $u \in H^{1}_{*}(\YI)$
  (so $u$ is a strong solution) and $u$ satisfies the \emph{a priori}
  estimate
  \begin{equation}
  \label{B:apriori}
        \nmi{u} \le
          \frac{c_{4}}{1-c_4\Norm{B_0}{op}}
                ( \Norm{f}{L^{2}(\YI)} + \nyh{\sigma}
         + \Norm{B_1}{L^2\to L^2} \Norm{u}{L^{2}(\YI)} ).
  \end{equation}
\end{Theorem}

\Proof Since $\Norm{B_1u}{L^{2}(\YI)} \le \Norm{B_1}{op}
\Norm{u}{L^{2}(\YI)}$, $u$ satisfies $ \tilde{L}u :=
(\partial_{x}+A+B_0)u=\tilde{f} $ where $\tilde{f}=f-B_1u\in L^{2}(\YI)$.
Lemma \ref{LemB3} constructs a solution $\bar{u}\in H^{1}_{*}(\YI)$ of
$\tilde{L}\bar{u}=\tilde{f}$ satisfying the same boundary conditions, and
it follows that $\bar{u}$ is also a weak solution.  By the Uniqueness Lemma
\ref{LemB4} we have $u=\bar{u}$ and thus $u\in H^{1}_{*}(\YI)$, as
required.  The estimate \bref{B3:est} of Lemma \ref{LemB3} (with $f$
replaced by $\tilde f$) leads directly to \bref{B:apriori}.  \QED

\section{Boundary regularity for first order systems}
\label{secC}


In this section we determine conditions under which the boundary
regularity results of \S\ref{secB} apply to a first order equation
of Dirac type (see \bref{C-4f}, \bref{C-6}) at the boundary,
with suitable boundary operator, to show $H^1$ regularity of an $L^2$
weak solution.

We assume $M$ is a smooth manifold with compact boundary $Y$,
$E\to M$ and $F\to M$ are real vector bundles over $M$
with scalar products, and $\cL$ is a first order
elliptic operator on sections of $E$ to sections of $F$, which in
local coordinates $x^j$, $j=1,\ldots,n$ takes the form
\begin{equation}
\label{C-1}
\cL u = a^j\partial_j u + b u\;,
\end{equation}
where $a^j$, $b$ are homomorphisms of $E$ to $F$ as before.  Note that we
do not assume $Y$ to be connected.

To apply the preceding results, the coefficients must satisfy the
interior regularity conditions \bref{IR-2}, and the boundary
restrictions must be defined and satisfy the corresponding conditions
in dimension $n-1$,
\begin{equation}
\label{C-2}
\begin{array}{rcl}
a^j\rY &\in&   W^{1,(n-1)^*}(Y)\cap C^0,\ j=1,\ldots,n, \\
b\rY &\in&  L^{(n-1)^*}(Y)\;.
\end{array}
\end{equation}
Conditions (\ref{IR-2},\ref{IR-3}) and \bref{C-2} will be assumed
throughout this section.

\begin{NRemark}\label{Rreg}
The $H^s$ conditions
\begin{equation}
\label{C-3}
\begin{array}{rcl}
a^j & \in& W^{s,2}_{\loc}(M) \cap C^0(M)\;,
\\
b & \in& W^{s-1,2}_{\loc}(M)\;,
\end{array}
\end{equation}
where
\begin{equation}
\label{C-4}
 \begin{array}{ccl}
     s =  n/2 & \textrm{\ for\ } & n\ge 4\;, \\
     s >  3/2 & \textrm{\ for\ } & n = 3\;, \\
     s =  3/2 & \textrm{\ for\ } & n = 2\;,
 \end{array}
\end{equation}
imply the interior \bref{IR-2} and boundary \bref{C-2} coefficient
regularity conditions, through the Sobolev embedding and trace
theorems \cite{Taylor96}.
The $C^0$ condition in \bref{C-3} is superfluous for $n=2,3$.
\end{NRemark}

Let $x=x^n$ be a boundary coordinate, defining a tubular neighbourhood
$Y\times[0,1] \subset M$ of $Y$ with local coordinates $(y^i,x)\in
Y\times[0,1]$ (where we identify $Y$ with $Y\times\{0\}$).  Let $dv_M$
be a volume measure on $M$, and define $dv_Y =
(-1)^{n-1}\partial_x\cut dv_M\rY$ on $Y\times\{0\}$.  The local
coordinate integration factor $\gamma$ is defined in $Y\times[0,1]$ by
$dv_M = \gamma \, dy\,dx$, where $dy\,dx$ is coordinate Lebesgue
measure.  We assume the local coordinate condition\footnote{As in
  \S\ref{secIR}, this means that we can cover $Y$ by a finite number
  of coordinate charts ${\mycal O}_\alpha $ so that $\gamma$ has the stated
  regularity in the local coordinates on ${\mycal O }_\alpha\times [0,1]$.}
\begin{equation}
\label{C-4d}
\gamma \in  (W^{1,n^*}\cap C^0)(Y\times[0,1])\;.
\end{equation}
In order to directly apply  the results of the previous section, we assume that
 $\gamma=\gamma(0)$ is independent of $x$ in $Y\times[0,1]$, so
 $dv_M =dv_Y\,dx$.  This involves no loss of generality, as the $x$-dependence
of $\gamma$ in the integral form \bref{C-8a} of the equation can
be absorbed through a rescaling of the coefficients $a^i,b$. Since
the restrictions to $Y\times \{0\}$ are unchanged, this does not affect the
boundary operator $A$.

To minimize confusion
with the outer unit normal $n=-\partial_x$ we set
\begin{equation}
\label{C-4e}
\nu := {}^ta^n\;,
\end{equation} and for simplicity we assume
\begin{equation}
\label{C-4f}
{}^t\nu = \nu^{-1}\;.
\end{equation}
This is satisfied by the Dirac operator.  \Eq{C-4f} can also be
achieved in several other situations of interest by pre-multiplying
$\cL$ by a suitable homomorphism, or by making a frame change in $F$
(this can be done {\em e.g.}~when $F$ is trivial).
 Note that such a pre-multiplication does not affect the values
of $\tilde{a}^i$, $\tilde{b}$, where we define
\begin{equation}
\label{C-5}
\tilde{a}^i := \nu a^i\rY,\ i=1,\dots,n-1 ,\quad
\tilde{b} :=  \nu b\rY\;.
\end{equation}
By extending independent of $x$ we regard $\tilde{a}^i,\tilde{b}$
as defined on $Y\times[0,1]$.  We assume the important boundary
symmetry condition
\begin{equation}
\label{C-6}
\tilde{a}^i = {}-{} {}^t\tilde{a}^i,\ \ i=1,\dots,n-1 \;,
\end{equation}
where the transpose ${}^t\tilde{a}^i$ is taken with respect to the
inner product on $E$.
We wrap up these conditions into a definition.
\begin{Definition}
A first order system \bref{C-1} is of \emph{boundary Dirac type} if
the coefficients $a^j$ satisfy the conditions (\ref{C-4f},\ref{C-6})
in some neighbourhood of the boundary.
\end{Definition}
%
As discussed in \S\ref{secE}, this class includes the examples of
\S\ref{secA} and \S\ref{Section:pmt}.  These conditions will be
assumed henceforth.

 Using \bref{C-6} we define the
boundary operator
\begin{equation}
\label{C-7}
Au := \sum_{i=1}^{n-1}\tilde{a}^i\partial_iu + \tilde{a}_0u + \tilde{b}_0u\;,
\end{equation}
where
\begin{equation}
\tilde{a}_0 =\half \sum_{i=1}^{n-1}\left(  
    \partial_i\tilde{a}^i +
    \tilde{a}^i\partial_i\log\gamma \right) \;,
\end{equation}
and $\tilde{b}_0 = b_0\rY$ for some symmetric endomorphism $b_0$ of
$E$. We require that $b_0\in L^{n^*}(Y\times I))$ and $\tilde{b}_0 \in
L^{(n-1)^*}(Y)$, compare \bref{IR-2}, \bref{C-2}. Then $A$ is formally
self-adjoint ($A^\dagger=A$) on the bundle $E\rY$ over $Y$, with
respect to the measure $dv_Y$.

Note that the choice of zero-order term $b_0$ gives some freedom in
the definition of the boundary operator $A$, which is thus not
uniquely determined by $\cL$.  Near the boundary, $\cL u =
{}^t\nu(\partial_x+A+B)$, which may be expressed as
\begin{equation}
\label{C-7b}
   \tilde{L}u := (\partial_x + A + B) u = \tilde{f}\;,
\end{equation}
where $\tilde{L}u=\nu \cL u$, $\tilde{f}=\nu f$, and $B$
denotes the difference
\begin{equation}
\label{C-7c}
  Bu = \sum_{i=1}^{n-1}(\nu a^i-\tilde{a}^i)\partial_iu +
  (\nu b -\tilde{a}_0 -  \tilde{b}_0)u\;.
\end{equation}

By Corollary \ref{CorSC3}, under the above conditions, $A$ will satisfy the
spectral condition \ref{DefSC}.  Denote the eigenvalue index set by
$\Lambda$ and fix a cutoff $\kappa$, as in \S\ref{secB}.  Similarly let
$\Lambda^+=\{\alpha\in\Lambda : \lambda_\alpha\ge\kappa\}$, fix some subset
$\hat\Lambda\subset \{\alpha\in\Lambda : |\lambda_\alpha| <\kappa\}$ and
let $P=P_{\Lambda^+}+P_{\hat{\Lambda}}$ be the associated spectral
projection operator.  Note that Theorem \ref{ThmIR5} (see \bref{IR-14})
shows that the $H^s_*$ norms defined using $A$ will be equivalent to the
corresponding $H^s$ norms defined on $Y$ and $\YI$, at least for $s=0,1$;
it seems likely that this will hold, by interpolation,
for all $s\in[0,1]$.
If the coefficient regularity allows an $H^{k+1}$ elliptic estimate
(\emph{cf.\/}~Theorem \ref{ThmIR6}) then this should extend to $s\in[0,
k+1]$.

The definition of weak solution on a manifold with boundary which we
are about to give is slightly simpler than the tubular neighbourhood
Definition \ref{Def:weak} of \S\ref{secB}, since the conditions
\bref{B:c}, \bref{B:adc} may be
imposed by localising with a boundary cutoff function; see the proof
of Theorem \ref{ThmC1}.
  Let $H^1_c(M)$ denote the dense subspace of
$H^1(M)$ consisting of functions of compact support, where we recall
that because $M$ is a manifold with boundary $Y=\partial M$,
$H^1_c(M)$ includes functions which are non-zero on $Y$.  As in
\S\ref{secB}, the boundary condition is expressed using a positive
spectrum projection $P:L^2(Y)\to L^2(Y)$ and a bounded linear map
$K:(1-P)L^2(Y)\to PL^2(Y)$ and its $L^2(Y)$ adjoint $K^\dag$.

\begin{Definition}\label{defC1}
Let $f\in L^2_\loc(M)$ and $\sigma\in PL^2(Y)$ be given.  A section
 $u\in L^2_\loc(M)$ is a
\emph{weak solution of $\cL u=f$ with boundary condition}
\begin{equation}
\label{C-8o}
Pu_0 = \sigma + K(1-P)u_0\;,
\end{equation}
if
\begin{equation}
\label{C-8a}
        \int_{M}\langle u,{\cL}^{\dag}\phi\rangle \,dv_M =
        \int_{M}\langle {f},\phi\rangle \,dv_M
        +\oint_{Y} \langle \sigma,\nu\phi_0\rangle \,dv_Y\;,
\end{equation}
for all $\phi\in H^1_c(M)$ satisfying the boundary condition
\begin{equation}
\label{C-8b}
  (1-P+K^\dagger P)(\nu\phi_0) = 0\;,
\end{equation}
where $\cL^\dag$ is the $L^2$ adjoint given by \bref{IR-5}.
\end{Definition}
Note we are using the notation $u_0,\phi_0$, {\em etc.}, to denote
the restriction (trace) on the boundary $Y$.  The additional term
$\nu$ in (\ref{C-8a},\ref{C-8b}) (cf.~(\ref{B:wLuf},\ref{B:adb}))
arises from the relation $\cL = {}^t\nu(\partial_x+A+B)$ between
$\cL$ and the boundary form $\partial_x+A+B$ used in \S\ref{secB}.

The boundary condition \bref{C-8o} restricts $u_0=u\rY$ to lie in the
affine subspace of $L^2(Y)$ given by the graph of $x\mapsto \sigma+Kx$
over the negative spectrum subspace $x\in (1-P)L^2(Y)$.  It will be useful
to re-express \bref{C-8o} as
\begin{equation}
\label{C-8c}
\cK u_0 = \sigma\;,
\end{equation}
where we have introduced the operator $\cK$ on $L^2(Y)$,
\begin{equation}
\label{C-8d}
\cK := P - K(1-P)\;,
\end{equation}
and likewise to re-express the ``adjoint'' boundary condition \bref{C-8b} as
\begin{equation}
\label{C-8e}
\cK^\dag \phi_0 = 0\;,
\end{equation}
where we define
\begin{equation}
\label{C-8f}
\cK^\dag := \nu^{-1}(1-P+K^\dag P)\nu\;.
\end{equation}

The next result generalizes the interior weak-strong Theorem \ref{ThmIR5}
to boundary value problems.

\begin{Theorem}
\label{ThmC1}
Suppose $\cL$ and $A$ satisfy the conditions
(\ref{IR-2},\ref{IR-3},\ref{C-2}, \ref{C-4d},\ref{C-4f},\ref{C-6}),
and suppose $\sigma\in PH^{1/2}_{*}(Y)$, $f\in L^{2}_\loc(M)$.
Further suppose $K:(1-P)L^2(Y)\to PL^2(Y)$ is bounded linear and
satisfies \bref{B:K1}, with $L^2$ adjoint $K^\dagger$ satisfying
\bref{B:K2}.  Assume $u\in L^2_\loc(M)$ is a weak solution of $\cL
u=f$ with the boundary condition $\cK u_0 = \sigma$ \bref{C-8c}.  Then
$u\in H^1_\loc(M)$ and $u$ is a strong solution. Moreover, there are
constants $\delta\in(0,1]$, $c_5$, depending only on $\kappa,k$ and
$a^i,b$, and intervals $I'=[0,\delta/2]$, $I=[0,\delta]$ such that
\begin{equation}
\label{C-9b}
\Norm{u}{H^1(Y\times I')}
\le c_5 ( \Norm{f}{L^{2}(\YI)} + \Norm{\sigma} {H^{1/2}_*(Y)}
         + \Norm{u}{L^{2}(\YI)})\;.
\end{equation}
\end{Theorem}

\Proof Theorem \ref{ThmIR5} ensures $u\in H^1_\loc(\mathring{M})$,
where $\mathring{M}$ is the interior of $M$, so it suffices to
consider $u$ compactly supported in $Y\times[0,\delta)$, for any
choice of $\delta\in(0,1)$.  In particular, because $u$ then vanishes
near $Y\times\{\delta\}$, it follows from \bref{C-8o} that $u$ is a
weak solution with boundary conditions, in the sense of Definition
\ref{Def:weak}.

It will suffice to show, for a sufficiently small choice of
$\delta>0$, that we may decompose $B=B_0+B_1$ into pieces satisfying
the size conditions of Theorem \ref{ThmB5}.
Write $B=\beta^i\partial_i + \beta$ where
\[
\beta^i(y,x) = (a^n(y,x))^{-1} a^i(y,x)- (a^n(y,0))^{-1} a^i(y,0),
\quad i=1,\ldots,n-1,
\]
so $\beta^i\in W^{1,n^*}\cap C^0$ and $\beta^i(y,0)=0$.  Since the
constant $c_4$ of Lemma \ref{LemB2} depends only on $\theta_0<1$,
$\ell=\kappa\delta$ and the constant $k$ of (\ref{B:K1},\ref{B:K2}),
it is bounded uniformly in $\delta\le1$.  Consequently for any
$\epsilon>0$ there is $\delta_0>0$ such that
$c_4\Norm{\beta^i}{L^\infty(Y\times[0,\delta])} < \epsilon$ for all
$\delta\le \delta_0$.

Likewise, since $\gamma\in W^{1,n^*}$ \bref{C-4d}, we have $\beta\in
L^{n^*}$ and there is a decomposition $\beta=\beta_0+\beta_1$ with $ c_4\,
C_S\, \Norm{\beta_0}{L^{n^*}(Y\times[0,1])} \le \epsilon$, where $C_S$ is
the Sobolev constant on $Y\times[0,1]$, and $\beta_1\in L^\infty$.  Then
$B_0=\beta^i\partial_i+\beta_0$ satisfies \bref{B3:bnd}, as does
$B_0^\dagger$ (possibly after decreasing $\delta$), and $B_1=\beta_1$ is
bounded on $L^2$.  Theorem \ref{ThmB5} now applies and shows $u\in
H^1(Y\times[0,\delta])$, since $H^1(\YI)=H^1_*(\YI)$ as remarked above.  The
elliptic estimate \bref{C-9b} follows by applying \bref{B:apriori} to
$\tilde{u}=\chi u$, where $\chi=\chi(x)$ is a cutoff function, $\chi(x)=1$
for $0\le x\le\delta/2$, $\chi(x)=0$ for $x\ge\tfrac{3}{4}\delta$.
\QED

\begin{Corollary}
\label{CorC1}
Suppose $\cL,A,K,\cK$ satisfy the conditions of Theorem \ref{ThmC1}. Then
for all $u\in H^1_{\loc}(M)$ we have the boundary estimate
\begin{equation}
\label{C-9c}
\Norm{u}{H^1(Y\times I')}
\le c_5 ( \Norm{\cL u}{L^{2}(\YI)} + \Norm{\cK u_0} {H^{1/2}_*(Y)}
         + \Norm{u}{L^{2}(\YI)})\;.
\end{equation}
\end{Corollary}

\Proof
If $u\in H^1_{\loc}(M)$ then $f:=\cL u\in L^2_{\loc}(M)$,
$\sigma:=\cK u_0\in H^{1/2}_*(Y)$, and $u$ is a strong solution.  Since
integration by parts may be applied to show $u$ satisfies the weak equation
\bref{C-8a}, Theorem \ref{ThmC1} applies and gives \bref{C-9c}.
\QED

A bootstrap argument, slightly more complicated than that used for
the interior bounds Theorem \ref{ThmIR6}, leads to higher,
$H^{1+k}$, regularity. Rather than stating
complicated conditions for general $k$, we describe the details
only for the case $k=1$ ($u\in H^2$). The coefficient regularity
conditions are most likely not optimal.
\begin{Theorem}
\label{ThmC2} In the setting of Theorem \ref{ThmC1}
let $u\in H^1(\YI)$ be the solution and suppose
the following additional regularity conditions are satisfied,
\begin{equation}
\label{C-10}
a^j\;, b\in 
W^{1,\infty}(\YI)\;,
\end{equation}
\begin{equation}
\gamma\in W^{2,n^*}(Y\times[0,1])\;,
\end{equation}
\begin{equation}
\label{C-11}
f\in H^{1}(\YI),\quad \sigma\in H^{3/2}_*(Y)\;,
\end{equation}
\begin{equation}
\label{C-12}
[A,K](1-P): (1-P)H^{1/2}_*(Y)\to PH^{1/2}_*(Y) \   \textrm{is bounded.}
\end{equation}
Then there exists $\delta''\le\delta/2$ such that $u\in H^{2}(Y\times I'')$, where $I''=[0,\delta'']$, and there is a constant $c_6$ depending on the
coefficient bounds (\ref{C-10})-(\ref{C-12}), such that
\begin{equation}
\label{C-13}
 \Norm{u}{H^{2}(Y\times I'')}
\le c_6 ( \Norm{f}{H^{1}(\YI)} + \Norm{\sigma} {H^{3/2}_*(Y)}
         + \Norm{u}{L^{2}(\YI)})\;.
\end{equation}
\end{Theorem}

\begin{remark} For the APS and chiral boundary conditions
\bref{APS}, \bref{GHHP}, $A$ commutes with $K$ and thus
\bref{C-12} is trivially satisfied.
\end{remark}


\Proof
The idea is to show that $Au$ satisfies a similar boundary value
problem. For convenience, let $\tilde{H}^1_0(\YI)$, $I=[0,\delta]$,
denote the $H^1$ completion of the $C^\infty$ functions of compact
support in $Y\times[0,\delta)$, where $\delta$ is the constant of
Theorem \ref{ThmC1}.
In particular, functions in $\tilde{H}^1_0(\YI)$ have vanishing
trace on $Y\times\{\delta\}$.  For any $v,\psi\in
C_c^\infty(Y\times[0,\delta))$, we have the identity
\begin{eqnarray}
\nonumber
\int_{\YI}\langle Av,\tilde{L}^\dagger\psi\rangle \, dx\,dv_Y
&=&
\int_{\YI}\left(\langle [\tilde{L},A]v,\psi\rangle
         + \langle \tilde{L}v,A\psi\rangle\right)\, dx\,dv_Y
\\
\label{C-15}
&&{}+
\oint_{Y}\langle v_0,A\psi_0\rangle \,dv_Y\;.
\end{eqnarray}
Since $A$ is formally self-adjoint with respect to $dv_Y$, this
formula follows by direct calculation.  The terms with
$\tilde{L}^\dagger\psi,\tilde{L}v$ are well-defined for $v,\psi\in
\tilde{H}^1_0(\YI)$.  Since $v_0,\psi_0\in H^{1/2}_*(Y)$ by Lemma \ref{LemB0},
the boundary integral extends also by writing it as
\begin{equation}
\label{C-15b}
\oint_{Y}\langle v_0,A\psi_0\rangle \, dv_Y =
\oint_{Y}\langle Jv_0,J^{-1}A\psi_0\rangle \,dv_Y\;.
\end{equation}
Here $J=(1+|A|)^{1/2}$ is defined as $Ju = \sum_{\alpha\in\Lambda}
u_\alpha (1+|\lambda_\alpha|)^{1/2} \psi_\alpha$, with $u=\sum_{\alpha\in\Lambda}
u_\alpha  \psi_\alpha$; $J^{-1}$ is defined similarly.  This
shows that
\[
\left|\oint_{Y}\langle v_0,A\psi_0\rangle \,dv_Y\right|
\le
c\Norm{v_0}{H^{1/2}_*(Y)} \Norm{\psi_0}{H^{1/2}_*(Y)},
\]
for all $v_0,\psi_0\in H^{1/2}_*(Y)$, where the constant $c$ is
determined by $A$.

More care is required to control the commutator $[\tilde{L},A]$,
which takes the form
\[
[\tilde{L},A]v  =
  \alpha^{ij}\partial_{ij}^2 v
+ \alpha_1^i\partial_i v + \alpha_2 v
\]
where
\begin{eqnarray*}
\alpha^{ij} &=& [\nu a^i,\tilde{a}^j],
\\
\alpha_1^i  &=& \nu a^j\partial_j\tilde{a}^i -
            \tilde{a}^j \partial_j(\nu a^i)
           +[\nu a^i,\tilde{a}_0],
\\
\alpha_2   &=& \nu a^i\partial_i\tilde{a}_0
  -\tilde{a}^i\partial_i(\nu b) \;.
\end{eqnarray*}
Observe that $\alpha^{ij}\rY=0$.  The coefficient conditions
\bref{C-10} and $v\in H^1$ ensure that $\alpha_1^i\partial_iv$ is in
$L^2$ and hence may be combined with the source term $A\tilde{f}$.
Likewise \bref{C-10} ensures $\alpha_2 v$ is bounded in $L^2$.

Let $\lambda_0\in\bR\backslash\Lambda$, so $A-\lambda_0$ has
trivial kernel and satisfies elliptic estimates
$\Norm{v}{H^{s+1}}\le c \Norm{(A-\lambda_0)v}{H^s}$ for
$s=0,1$ at least. Since $A$ is self-adjoint and elliptic, the
cokernel of $A-\lambda_0$ is also trivial, so
$A-\lambda_0:H^{s+1}\to H^s$ is invertible. Now decompose
$\alpha^{ij}\partial_{ij}^2 v = B_2Av + B_3 v$ where
\[
  B_2 = \alpha^{ij}\partial_{ij}^2(A-\lambda_0)^{-1}
\]
is bounded from $H^1\to L^2$, and
\[B_3 = -\lambda_0\alpha^{ij}\partial_{ij}^2(A-\lambda_0)^{-1}
\]
is also bounded from $H^1\to L^2$.  Note that by perhaps
decreasing $\delta$ we may ensure that $B_2$ and $B_2^\dagger$
satisfy a smallness condition similar to \bref{B3:bnd}.

Direct calculation (noting that $[A,P]=0$) establishes the
boundary formula
\begin{eqnarray}
\nonumber
\oint_{Y}\langle (1+K)(1-P)v_0,A\psi_0\rangle \,dv_Y
&=&
\oint_{Y}\langle Av_0, (1-P+K^\dagger P)\psi_0\rangle \,dv_Y
\\
\label{C-16}
&&{}
+ \oint_{Y}\langle [A,K](1-P)v_0,\psi_0\rangle \,dv_Y\;,
\end{eqnarray}
for all $v_0,\psi_0\in H^{1/2}_*(Y)$, since  $K$ satisfies
\bref{C-12} and \bref{B:K1} by assumption.

Now $u\in \tilde{H}^1_0(\YI)$ satisfies $\tilde{L}u=\tilde{f}$ and
$u_0=\sigma+(1+K)(1-P)u_0$.  Substituting $u$ for $v$ in \bref{C-15}
and using these relations, shows that $u$ satisfies
\begin{eqnarray}
\nonumber
\int_{\YI}\langle Au,(\tilde{L}-B_2)^\dagger\psi\rangle\, dx\,dv_Y
&=&
\int_{\YI} \langle A\tilde{f} + \alpha_1^i\partial_i u
            + \alpha_2 u,\psi\rangle\,dx\,dv_Y
\\
\nonumber
&&{}
 +\oint_Y \langle A\sigma+[A,K](1-P)u_0,\psi_0\rangle\,dv_Y
\\
\label{C-17}
&&{}+
\oint_{Y}\langle Au_0,(1-P+K^\dagger P)\psi_0\rangle\, dv_Y\;.
\end{eqnarray}
In particular, if $\psi_0\in\ker(1-P+K^\dagger P)\cap H^{1/2}_*(Y)
$ then $w=Au \in L^2(\YI)$ satisfies
\begin{eqnarray}
\nonumber
\int_{\YI}\langle w,(\tilde{L}-B_2)^\dagger\psi\rangle\, dx\,dv_Y
&=&
\int_{\YI} \langle \tilde{f}_1,\psi\rangle\,dx\,dv_Y
\\
&&{}
 +\oint_Y \langle \sigma_1,\psi_0\rangle\,dv_Y,
\label{C-18}
\end{eqnarray}
for all $\psi\in H^1(\YI)$ such that $\psi\in\ker(1-P+K^\dagger P)
$, where
\begin{eqnarray*}
\tilde{f}_1 &=& A\tilde{f} + \alpha_1^i\partial_i u + \alpha_2 u \in L^2,
\\
\sigma_1 &=& A\sigma+[A,K](1-P)u_0 \in H^{1/2}_*(Y)\;.
\end{eqnarray*}

In other words, $w\in L^2(\YI)$ is a weak solution of the problem
\begin{eqnarray*}
(\tilde{L}-B_2)w &=& \tilde{f}_1,
\\
Pw_0 &=& \sigma_1 + K(1-P)w_0\;.
\end{eqnarray*}
By shrinking the boundary layer we may assume $\Norm{B_2}{H^1\to L^2}$
and $\Norm{B_2^\dagger}{H^1\to L^2}$ are sufficiently small that the
conditions of Theorem \ref{ThmC1} are met, so $w=Au\in
H^1(Y\times[0,\delta])$.  The equation now gives $\partial_x u =
\tilde{f}-Au-Bu\in H^1$ and thus $u\in H^2(Y\times[0,\delta])$.
\QED

\section{Fredholm properties on compact manifolds}
\label{secDI}
The interior and boundary estimates of \S~\ref{secC} lead to
solvability (Fredholm) results, by standard arguments. The main
interest lies in identifying the cokernel, and we give a simple
necessary and sufficient condition for solvability, in Theorem
\ref{ThmDI3}. This section treats only compact manifolds, leaving the
more difficult case of non-compact manifolds to the following section.
Because more detailed descriptions are given in \S~\ref{secD}, some of
the arguments are only briefly summarised here.


Throughout this section we assume the coefficients $a^j\;,
j=1,\dots,n$ and $b$ of $\cL$ satisfy the conditions of
\S~\ref{secC}, namely \bref{IR-2},\bref{IR-3}, \bref{C-2},
\bref{C-4d}, \bref{C-4f}, \bref{C-6}; $\cL^\dagger$ is given by
\bref{IR-5}, and the boundary operators $K, K^\dag$ satisfy
(\ref{B:K1},\ref{B:K2}), where $P=P_{\Lambda}+P_{\hat{\Lambda}}$
is a positive spectrum projection of $A$ \bref{C-7}, and
$\cK\;,\cK^\dagger$ are defined by (\ref{C-8d},\ref{C-8f}).


Recall the Sobolev space $H^1(M)$ of sections of $E$ over $M$ is
 defined by the norm \bref{IR-14b}
\begin{equation}
\label{DI-1}
\Norm{u}{H^1(M)}^2 = \int_M (|\nabla u|^2 +  |u|^2)\,dv_M\;,
\end{equation}
where lengths are measured using the metric $\langle\ ,\ \rangle$ on
$E$ and a fixed smooth background metric $\go$ on $TM$, and the
connection $\nabla$ satisfies (\ref{IR-14c},\ref{IR-14d}).  Note again
that $\nabla$ need not be compatible with the metric on $E$, and the
space $H^1(M)$ is independent of the choice of $\nabla$.

The following basic elliptic estimate extends  \bref{IR-14}
of Theorem \ref{ThmIR5} to manifolds with boundary, using the boundary
neighbourhood estimate \bref{C-9b} of Theorem \ref{ThmC1}.

\begin{Proposition}
\label{PropDI1}
There is a constant $C>0$ depending on $a^j\;,b,\Gamma$ and $\cK$ such
that for all $u\in H^1(M)$,
\begin{equation}
\label{DI-3}
\Norm{u}{H^1(M)} \le C(\Norm{\cL u}{L^2(M)} + \Norm{\cK u_0}{H^{1/2}_*(Y)}
    +\Norm{u}{L^2(M)})\;.
\end{equation}
\end{Proposition}

\Proof
The argument used in Theorem \ref{ThmIR5} to prove the interior
estimate \bref{IR-14} may be applied using Theorem \ref{ThmC1},
estimate \bref{C-9b}, to estimate $\Norm{u_\alpha}{H^1(U_\alpha)}$ over
boundary neighbourhoods $U_\alpha$.  The remaining details are unchanged.
\QED

\begin{Theorem}
\label{ThmDI2}
The linear operator
\begin{equation}
\label{DI-4}
   (\cL,\cK):H^1(M)\to L^2(M)\times PH^{1/2}_*(Y)
\end{equation}
is \emph{semi-Fredholm}\/ (i.e. has finite dimensional kernel and
closed range).
\end{Theorem}

\Proof Suppose $\{u_k\}_1^\infty$ is a sequence in $\ker(\cL,\cK)$,
normalised by $\Norm{u_k}{H^1(M)}=1$.  To show the kernel is finite
dimensional, it suffices to show there is a subsequence converging in
$H^1(M)$.  By Rellich's lemma there is a subsequence (which we also denote
$u_k$) which converges strongly in $L^2(M)$, to $\bar{u}\in L^2(M)$ say.
The elliptic estimate applied to the differences $u_j-u_k$ shows that the
sequence is Cauchy in $H^1(M)$ and thus converges strongly to $\bar{u}\in
H^1(M)$.  Since \bref{DI-4} is bounded, it follows that
$\bar{u}\in\ker(\cL,\cK)$, so the unit ball in the kernel is compact and
hence the kernel is finite dimensional.

To show the range is closed, let $\mathring{H}^1(M)$ be the finite codimension
subspace of $H^1(M)$ defined by the condition
\[
   \int_M (\langle \nabla u,\nabla\phi\rangle + \langle
   u,\phi\rangle)\,dv_M
  = 0\quad \forall\ \phi\in \ker(\cL,\cK).
\]
A Morrey-type argument by contradiction using \bref{DI-3} shows there is a
constant $C>0$ such that for all $u\in \mathring{H}^1(M)$,
\begin{equation}
\label{DI-5}
C^{-1} \int_M |u|^2\,dv_M \le  \int_M |\cL u|^2 \,dv_M
+ \oint_Y|J\cK u_0|^2\,dv_Y\;,
\end{equation}
where $J = (1+|A|)^{1/2}$.  Now suppose $\{u_k\}_1^\infty \subset H^1(M)$
is such that $\cL u_k =f_k \to f\in L^2(M)$ and $\cK (u_k)_0 = s_k\to
\sigma\in PH^{1/2}_*(Y)$.  (Note that by the definition \bref{C-8d} of
$\cK$, the range of $\cK$ is a subspace of $PH^{1/2}_*(Y)$).  Since the
kernel is finite dimensional we may normalise $u_k\in \mathring{H}^1(M)$,
and then
\bref{DI-5} and \bref{DI-3} show that $\{u_k\}_1^\infty$ is bounded in
$H^1(M)$.  It then follows as above that there is a subsequence converging
strongly in $H^1(M)$ to $\bar{u}$, and that $\cL\bar{u} =
\lim_{k\to\infty}f_k = f$ and $\cK\bar{u}_0=\lim_{k\to\infty}s_k =\sigma$,
so the range of $(\cL,\cK)$ is closed.  \QED

The general boundary value problem
\begin{equation}
\label{DI-5b}
\left\{ \begin{array}{rclr}
        \cL u   & = & f & \ \ \textrm{in } M\\
        \cK u_0 & = & \sigma & \ \ \textrm{on } Y
    \end{array}
\right.
\end{equation}
is solvable for $u\in H^1(M)$ provided $(f,\sigma)$ satisfies the condition
\bref{DI-6} of the following main result.

\begin{Theorem}
\label{ThmDI3}
$(f,\sigma)\in L^2(Y)\times PH^{1/2}_*(Y)$ lies in the range of
$(\cL,\cK)$ (that is, \bref{DI-5b} admits a solution $u\in H^1(M)$),
 if and only if
\begin{equation}
\label{DI-6}
\int_M \langle f,\phi\rangle\,dv_M +
\oint_Y \langle \sigma,\nu\phi_0\rangle\,dv_Y =0
\quad \forall\ \phi\in\ker(\cL^\dag,\cK^\dag)\;.
\end{equation}
\end{Theorem}

\Proof
If $u\in H^1(M)$ satisfies $\cL u=f$ and $\cK u_0=\sigma$ then $u$ is
also a weak solution. Condition \bref{DI-6} then follows directly from
the definition \ref{defC1} of weak solution, hence \bref{DI-6}
is a necessary condition for solvability.

To establish the converse,  consider first the case $\sigma=0$.
Thus we suppose $f\in L^2(M)$ satisfies $\int_M \langle
f,\phi\rangle\,dv_M =0$ for all $ \phi\in\ker(\cL^\dag,\cK^\dag)$, and
we must find $u\in H^1(M)$ satisfying $\cL u=f$, $\cK u_0=0$.

By Lemma \ref{LemB0} the trace map $r_Y:u\mapsto u_0$ is bounded,
hence
\begin{eqnarray}
\nonumber
\mathring{H}^1_{\cK} &:=& \{u\in H^1(M): \cK u_0 =0, \textrm{  and}
\\
\label{DI-7}
&& 
   \int_M (\langle \nabla u,\nabla\phi\rangle + \langle
   u,\phi\rangle)\,dv_M
  = 0\quad \forall\ \phi\in \ker(\cL,\cK)  \}
\end{eqnarray}
is a closed  subspace of $H^1(M)$.
The argument of Theorem \ref{ThmDI2} (ii) shows there is a constant
$C$ such that
\begin{equation}
\label{DI-8}
\int_M (|\nabla u|^2+|u|^2)\,dv_M \le C\int_M|\cL u|^2\,dv_M
\end{equation}
for all $u\in \mathring{H}^1_{\cK}(M)$.  In particular, $\int_M|\cL
u|^2\,dv_M$ is strictly coercive on $\mathring{H}^1_{\cK}$, so the
Lax-Milgram lemma gives $u\in \mathring{H}^1_{\cK}$ satisfying
\[
\int_M \langle f,\cL\phi\rangle\,dv_M =
\int_M \langle \cL u,\cL\phi\rangle\,dv_M
\]
for all $\phi\in \mathring{H}^1_{\cK}$.  This equality also holds if
$\phi\in\ker(\cL,\cK)$, so $\Psi=\cL u-f$ satisfies
\begin{equation}
\label{DI-9}
\int_M \langle \Psi,\cL\phi\rangle\,dv_M =0
\quad \forall \ \phi\in H^1(M)\;,\ \cK\phi_0=0\;.
\end{equation}
Lemma \ref{LemB4a} and the identity
\begin{equation}
\label{DI-9b}
\int_M \langle \Psi,\cL\phi\rangle\,dv_M =
\int_M \langle \cL^\dag\Psi,\phi\rangle\,dv_M
-\oint_Y \langle \nu\Psi_0,\phi_0\rangle\,dv_M
\end{equation}
show that \bref{DI-9} is the weak form of the
adjoint problem
\begin{equation}
\label{DI-10}
\cL^\dag \Psi =0,\quad \cK^\dag \Psi_0=0\;.
\end{equation}

By \bref{IR-5}, $\cL^\dag$ is elliptic with boundary representation
\[
\cL^\dag
= -\nu(\partial_x+\hat{A}+\hat{B}),
\]
 where $\hat{A}=-\nu^{-1}A\nu$
since $A^\dag=A$, and $\hat{B}=-\nu^{-1}B^\dag\nu$.  By
(\ref{C-5},\ref{C-6}) the leading terms in $\hat{A}$ are
$\tilde{a}^i\partial_i$ so $\hat{A}$ is elliptic on $Y$, and
self-adjoint by \bref{C-4f}.  Since
$\hat{A}(\nu^{-1}\phi_\alpha)=-\lambda_\alpha\nu^{-1}\phi_\alpha$ if
$A\phi_\alpha = \lambda_\alpha\phi_\alpha$, we see that $\hat{A}$
satisfies the spectral conditions, and $\mathrm{spec}\,\hat{A}
=-\mathrm{spec}\,A$.  (Note that in the usual case of Dirac operators,
$\hat{A}=A$ and the spectrum is symmetric).  Now $\hat{P} :=
1-\nu^{-1}P\nu$ is a positive eigenspace projector for $\hat{A}$, with eigenvalues
$-\lambda_\alpha$ for $\alpha\in \Lambda^-\cup
(\Lambda^0\backslash\hat{\Lambda})$,  and the boundary
operator satisfies
\begin{equation}
\label{DI-11}
\cK^\dag \Psi_0 = (\hat{P}+\nu^{-1}K^\dag\nu(1-\hat{P}))\Psi_0\;.
\end{equation}
Since $\hat{K}=-\nu^{-1}K^\dagger\nu$ maps negative eigenvectors (of
$\hat{A}$) to positive eigenvectors, it follows that
$\cK^\dag\Psi_0=0$ is an elliptic boundary condition for $\cL^\dag$.
The boundedness conditions (\ref{B:K1},\ref{B:K2}) for $\hat{K}$
follow from the corresponding conditions for $K$.

Since $(\cL^\dag,\cK^\dag)$ is elliptic and satisfies the conditions
for Theorem \ref{ThmC1}, we conclude that $\Psi\in H^1(M)$ and
$\Psi$ satisfies the strong form \bref{DI-10}.

Since  $\Psi\in\ker(\cL^\dag,\cK^\dag)$,  assumption \bref{DI-6} with
$\sigma=0$ gives
\[
    \int_M\langle f,\Psi\rangle \,dv_M=0.
\]
By construction $\cK u_0=0$, so we may use $u$ as a test function
in the weak form \bref{DI-9} of the equation satisfied by $\Psi$,
giving
\[
    \int_M\langle \cL u,\Psi\rangle \,dv_M=0.
\]
It follows from $\Psi=\cL u-f$  that $\Psi=0$ and thus
$u$ is the required solution.

Now consider the case $\sigma\ne0$.  By Lemma \ref{LemB0} there is an
extension $v=e_Y(\sigma)\in H^1(M)$ supported in a neighbourhood of
$Y$ such that $v_0=\sigma$, $\Norm{v}{H^1(M)}\le
2\Norm{\sigma}{H^{1/2}_*}$.  Let $\tilde{f}=f-\cL v$ and consider the
equation
\begin{equation}
\label{DI-12}
\cL\tilde{u}=\tilde{f},\quad \cK\tilde{u}_0=0\;.
\end{equation}
The previous case shows
there is a solution provided $\tilde{f}$ satisfies
\[
    \int_M\langle \tilde{f},\psi \rangle \,dv_M=0\quad \forall\
\psi\in \ker(\cL^\dag,\cK^\dag)\subset H^1(M).
\]
Now \bref{DI-9b} shows that for all $\psi\in\ker(\cL^\dag,\cK^\dag)$,
\begin{eqnarray*}
   \int_M\langle \tilde{f},\psi \rangle \,dv_M
&=&    \int_M\langle f,\psi \rangle \,dv_M
    -    \int_M\langle v,\cL^\dag\psi \rangle \,dv_M
    +    \oint_Y\langle v_0,\nu\psi_0 \rangle \,dv_Y
\\
&=&    \int_M\langle f,\psi \rangle \,dv_M
    +     \oint_Y\langle \sigma,\nu\psi_0 \rangle \,dv_Y\;.
\end{eqnarray*}
Thus if \bref{DI-6} is
satisfied then there exists a solution $\tilde{u}$ of \bref{DI-12},
and then $u=\tilde{u}+v$ is the required full solution.  This
establishes sufficiency for the condition \bref{DI-6}.

\QED
We note two important consequences of Theorems \ref{ThmDI2}, \ref{ThmDI3}.
\begin{Corollary}
\label{CorDI4}
\bref{DI-5b} admits a solution for all $(f,\sigma)\in L^2(M)\times
PH^{1/2}_*(Y)$ if and only if $\ker(\cL^\dag,\cK^\dag)=\{0\}$.
\end{Corollary}

\begin{Corollary}
\label{CorDI5}
 $ (\cL,\cK):H^1(M)\to L^2(M)\times PH^{1/2}_*(Y)$  is Fredholm.
\end{Corollary}
\Proof
The argument of Theorem \ref{ThmDI3} shows that $(\cL^\dagger,\cK^\dagger)$
is elliptic and thus has finite dimensional kernel by Theorem \ref{ThmDI2}.
Now \bref{DI-6} shows that the range of $(\cL,\cK)$ has finite codimension.
\QED

\section{Fredholm properties on complete noncompact manifolds}
  \label{secD}

In this section we establish conditions under which the Fredholm and
existence results of the previous section for the operator
$(\cL,\cK):H^1(M)\to L^2(M)\times H^{1/2}_*(Y)$, may be extended to
non-compact manifolds.  This includes in particular, a generalisation of
the solvability criterion \bref{DI-6} of Theorem \ref{ThmDI3}.  Results of
this type may be applied to establish positive mass results in general
relativity, for example.


The non-compactness of $M$ causes some difficulties not found in the
compact case.  A classical result \cite{Wolf73,LawsonMichelsohn89}
shows that a Dirac operator $\cD$ on a non-compact manifold is
essentially self-adjoint on $L^2(M)$.  However, this elegant result is
useless for our purposes, since it implies only that
$\{(\phi,\cD\phi):\phi\in\textrm{dom}\,\cD\subset L^2(M)\}$ is closed
in the graph topology on $L^2(M)\times L^2(M)$.  This is weaker than
the closed range property, which is necessary for useful solvability
criteria.  In fact, because $L^2(M)$ often does not encompass natural
decay rates of solutions, the self-adjoint closure may not have closed
range.  In such cases the Dirac operator defined on $L^2(M)$ will not
be semi-Fredholm.  This is shown explicitly in the following example.

Consider the self-adjoint
closure $\overline{\cD}:\textrm{dom}\, \overline{\cD}\subset
L^2(\bbR^3) \to L^2(\bbR^3)$ of the constant coefficient Dirac operator
$\cD = \gamma^i\partial_i$ and let $f=\cD u$,
$u=(1-\chi)|x|^{-1}\psi$, where $\chi(r)$ is a smooth compactly
supported function identically one around $0$ and $\psi$ is a constant
spinor on $\bR^3$.  Clearly $f\in L^2(\bbR^3)$ but $u\not\in
L^2(\bbR^3)$, so in particular, $u
\not\in\textrm{dom}\,\overline{\cD}$.  However, $f$ still lies in the
closure of the range of $\overline{\cD}$, since $\cD(\chi_Ru) =
\chi_Rf+D\chi_Ru\to f$ in $L^2(\bbR^3)$, where $\chi_R(x)=\chi(x/R)$,
but $\chi_Ru$ can not converge in $L^2(\bbR^3)$.  Clearly
$(u,f)\not\in \textrm{graph}\,\overline{\cD}$ since $u\not\in
L^2(\bbR^3)$, and it can be shown (using the corresponding
Schr\"odinger-Lichnerowicz identity) that there is no $\bar{u}\in
L^2(\bbR^3)$ satisfying $\cD\bar{u}=f$.  Thus the self-adjoint closure
$\overline{\cD}$ does not have closed range.

In order to obtain an operator with closed range, it is thus necessary
to enlarge the domain, which raises the question of determining the
appropriate decay rate.  We sidestep this problem by using the $L^2$
size of the covariant derivative as a norm.  To obtain sufficient
control on the $L^2_\loc$ behaviour, we then must postulate a \emph{
  weighted Poincar\'e inequality} \bref{D:WP}.  The existence of such
inequalities can be established for the applications of most interest
in general relativity; see Proposition~\ref{PropD1} and \S\ref{secF}.

The elliptic estimate \bref{DI-3} plays a central role in the analysis
over a compact manifold, but its noncompact analogue cannot be obtained
directly by similar localisation arguments.  However, in cases of geometric
interest an identity of Schr\"odinger-Lichnerowicz form (generalising
\bref{A-2}) is available, and can be used to construct suitable
\emph{global} estimates.

The weighted Poincar\'e and Schr\"odinger-Lichnerowicz estimates are
the two additional ingredients needed for establishing solvability and
Fredholm properties on a non-compact manifold.

For ease of further reference, let us summarize the hypotheses
which will be made throughout this section:
\begin{hypotheses}\label{Hh}
$M$ is a non-compact manifold with compact boundary $Y$,
which is complete with respect to a
$C^\infty$ background metric $\go$. The case $Y=\emptyset$ is
admitted. The operator $\cL = a^j\partial_j +b$ satisfies the
global uniform ellipticity and boundedness condition
\begin{equation}
\label{D-1}
\eta^2 |V|^2
\le \go_{jk}\langle a^j(x)V,a^k(x)V\rangle
\le \eta^{-2} |V|^2\;,
\end{equation}
for some $\eta>0$, for all $V\in E_x$ and all $x\in M$.  The
coefficients of $\cL$ satisfy the interior regularity conditions
\bref{IR-2}, and the boundary regularity and structure conditions
of \S\ref{secC}, namely \bref{C-2}, \bref{C-4d}, \bref{C-4f},
\bref{C-6}.  Let $A$ be the boundary operator and $P$ its
associated positive spectrum projection, as in \S\ref{secC}. The
boundary operator $K:(1-P)L^2(Y)\to PL^2(Y)$ satisfies
(\ref{B:K1},\ref{B:K2}), and  $\cK\;,\cK^\dagger$ are defined in
(\ref{C-8d},\ref{C-8f}). The connection \bel{coneq}\nabla=\partial
-\Gamma\ee satisfies (\ref{IR-14c},\ref{IR-14d}) and we note again
that $\nabla$ need not be compatible with the metric on $E$ ---
this is important in some applications. \end{hypotheses}
 We may
express $\cL$ in terms of $\nabla$ by
\[
\cL = a^j\nabla_j +(b+a^j\Gamma_j) = a^j\nabla_j + \beta,
\]
where $\beta \in L^{n^*}_\loc(M)$. Additional, rather weak, decay
conditions will be imposed on $\beta$ \bref{D-5new}, on the negative
part of the curvature endomorphism $\rho$ \bref{D-5}, and on
$\Gamma^S=\half(\Gamma+{}^t\Gamma)$ in \S\ref{secF}.


\begin{Definition}
\label{defD1} The covariant derivative $\nabla$ on $E$ over $M$
admits a \emph{weighted Poincar\'e inequality} if there is a
weight function $w\in L^1_\loc(M)$  with $\mbox{\rm ess}\inf_\Omega w >0$ for all
relatively compact $\Omega \Subset M$, such that for all $u\in
C_c^1(M)$ we have
\begin{equation}
\label{D:WP}
\int_M |u|^2\,w\,dv_M \le \int_M |\nabla u|^2\, dv_M\;.
\end{equation}
\end{Definition}

Here the length $|\nabla u|^2$ is measured by the metric on $E$ and
the background Riemannian metric $\go$ on $M$, and $dv_M$ is the
volume measure of $\go$.  
It is clear that the weight function $w$ can be chosen to be smooth.

The semi-norm
\begin{equation}
\label{D-2a}
\Norm{u}{\hb}^2 = \int_M |\nabla u|^2 \,dv_M
\end{equation}
on $C_c^\infty(M)$ may be completed to form the space
\begin{equation}
\label{D-2b}
\hb := \Norm{\cdot}{\hb}\textrm{-completion of }C_c^\infty\Gamma(E)\;,
\end{equation}
which consists of equivalence classes of $\hb$-convergent
sequences in $C_c^\infty(M)$.  The weighted Poincar\'e inequality
\bref{D:WP} ensures that an $\hb$-convergent sequence converges
locally in $L^2$, so the equivalence classes may be identified
with cross-sections in the usual Lebesgue sense: with
cross-sections having coefficient functions agreeing
$dv_M$--a.e.

If there is no weighted Poincar\'e inequality, then it may be that
$\hb$ can not be identified with a space of Lebesgue-measurable
cross-sections in this sense.  For example, the trivial spinor
bundle over $M=\bbT^2\times\bR$ with the flat connection
$\nabla_i=\partial_i$ admits a global parallel spinor $\nabla_i
\psi=0$ which is approximated in the $\hb$ seminorm by $\psi_k =
\chi(x/k)\psi$ for $\chi\in C_c^\infty(\bR)$, $\chi=1$ on
$[-1,1]$.  Now $\int_M |\nabla \psi_k|^2\,dv_M \to 0$, but
$\lim_{k\to\infty}\psi_k = \psi\ne0$, so the $\hb$-equivalence
class $[0]$ contains $\psi\ne0$ everywhere.  In other words,
\bref{D-2a} does not define a norm on spinors in this example.
This shows, {\em inter alia}, that \bref{D:WP} will not hold in
all cases.

More generally, a weighted Poincar\'e inequality fails for manifolds
of the form $N\times\bR$, where $N$ is compact and itself admits a
parallel spinor.  It follows from the proof of Theorem \ref{ThmF2}
below that in such cases the orthogonal complement in $\hb$ of the
subspace of all parallel spinors will admit a weighted Poincar\'e
inequality.  Note also that the presence of a weighted Poincar\'e
inequality \bref{D:WP} does not imply there are no global parallel
spinors --- $\bR^3$ provides a simple counterexample.

However, weighted Poincar\'e inequalities can be demonstrated in many
cases of interest. In the next section we will prove:
\begin{Proposition}
\label{PropD1}
A covariant derivative $\nabla$ on $E$ admits a weighted Poincar\'e
inequality if any one of the following conditions holds:
\begin{enumerate}
\item
there is a relatively compact domain $\Omega\subset M$ and a constant $c>0$
such that
\begin{equation}
\label{D-4}
\int_\Omega |u|^2\,dv_M \le c\int_M|\nabla u|^2\,dv_M
\end{equation}
for all $u\in C_c^\infty(M)$;
\item there are no nontrivial globally parallel sections ($\nabla u=0
  \Rightarrow u=0$);
\item $M$ has a {weakly asymptotically flat end} $\tM $ (see
  Definition \ref{defF1}), with $\dim M\ge 3$;
\item $M$ has a weakly asymptotically hyperboloidal end (see
  Definition \ref{defF1.&}), with $\dim M \ge 2$.
\end{enumerate}
\end{Proposition}

When $M$ is non-compact, the global G\aa rding inequality  (generalizing
\bref{DI-3}) cannot be constructed from local estimates.  Motivated by some
classical and fundamental identities, we instead introduce the following
definition.

\begin{Definition}
\label{defD2}
The operator pair $(\cL,\cK)$ admits a
\emph{Schr\"odinger-Lichnerowicz estimate} if there is $C>0$ and a
non-negative function $\rho$ such that\footnote{ The function $\rho$
  here should not be confused with the energy density arising in
  general relativity: in Section~\ref{Section:pmt} $\rho$ will be
  zero.}  for all $u\in C_c^1(M)$,
\begin{equation}
\label{D:SL}
C^{-1}\int_M |\nabla u|^2 \,dv_M \le \int_M (|\cL u|^2 + \rho|u|^2)\,dv_M
 + \oint_Y |J\cK u_0|^2\,dv_Y\;,
\end{equation}
where  $J=(1+|A|)^{1/2}$.
\end{Definition}

\begin{Lemma}
\label{Lemd0}
Suppose that the {Schr\"odinger-Lichnerowicz estimate} \bref{D:SL}
holds for all $u\in C^1_c(M)$ with $\rho$ and $\beta=\cL-a^j\nabla_j $
satisfying $\rho\in L^{n^*/2}_\loc$, $\beta\in L^{n^*}_\loc$, and
\begin{equation}
\label{D-5new}
\lim_{R\to\infty} \sup_{M\backslash M_R} \frac{\rho+|\beta|^2}{w} <  \infty\;,
\end{equation}
where $\{M_R\}_{R\to\infty}$ is an exhaustion of $M$.  Then
$\cL:\hb\to L^2(M)$ is bounded and
\bref{D:SL} holds for all $u\in\hb$.
\end{Lemma}
\proof
It will suffice to show that the individual terms of the
right-hand-side of \bref{D:SL} are bounded by
$\|u\|_{\hb}^2$.
  Now
\begin{equation}
\label{D-5b}
\int_{M} |\cL u|^2
\le C \int_{M} |\nabla u|^2 + 2 \int_{M} |\beta|^2 |u|^2 \;,
\end{equation}
and we use \eq{D-5new} and \eq{D:WP} to estimate
\begin{eqnarray*}
\int _{M\setminus M_R} (\rho+ |\beta|^2) |u|^2 &\le&
\sup_{M\backslash M_R} \frac{\rho+|\beta|^2}{w}
\int_{M\setminus M_R}  |u|^2 w
\\ &\le & \sup_{M\setminus M_R} \frac{\rho+|\beta|^2}{w}
\int_{M\setminus M_R}  |\nabla u|^2
\\ &\le&  C \int _{M}  |\nabla u|^2 \;,
\end{eqnarray*}
for some $R<\infty$.  Let $\chi_R\in C^\infty_c(M)$ be a cut-off
function with support contained in $M_{2R}$, $\chi_R=1$ on $M_R$. Then
\begin{eqnarray*}
\int_{M_{R}} (\rho+|\beta|^2) |u|^2 &\le&
\int_{M_{2R}} (\rho+|\beta|^2) |\chi_R u|^2
\\ &
\le &\left(
  \|\rho\|_{L^{n^*/2}(M_{2R})}+\|\beta\|^2_{L^{n^*}(M_{2R})}\right)
\|\chi_R u\|^2_{L^{n^*}(M_{2R})} \;.
\end{eqnarray*}
Applying the Sobolev inequality for $\nabla$ on the compact set
$M_{2R}$ and the weighted Poincar\'e inequality show that the last
term is controlled by $\int_M |\nabla u|^2$.  Finally, the $K$-bound
\eq{B:K1} and the restriction Lemma \ref{LemB0} show that the boundary
term is also controlled by $\int_M |\nabla u|^2$.  \qed

Schr\"odinger-Lichnerowicz \emph{identities} hold for many common
examples, and can easily be adapted to produce estimates of the form
\bref{D:SL}. We will not attempt to give general conditions which imply
such inequalities --- it is simpler to ask only that \bref{D:SL} be
established separately in any particular case of interest.

For example, consider the classical Dirac operator $\cD$ of the
metric $g$ as in \S\ref{secA}, on a non-compact spin manifold $M$.
Combining \bref{A:Lich} and \bref{A:Lichb} gives
\begin{equation}
\label{D:Lich}
\int_M |\nabla\psi|^2\,dv_M = \int_M
(|\cD\psi|^2 - \tfrac{1}{4}R(g)|\psi|^2)\,dv_M
+ \oint_Y\langle \psi_0,(\cD_Y+\half H_Y)\psi_0\rangle\, dv_Y\;,
\end{equation}
for any $C^1_c$  spinor field on $M$.  Suppose the boundary operator is $\cK =
P_+$, the orthogonal projection onto the positive spectrum
eigenspinors of $\cD_Y$. If the boundary mean curvature $H_Y$
satisfies $H_Y\le \sqrt{16\pi/\textrm{Area}(Y)}$, then the argument in
\S\ref{secA} shows that the boundary term in \bref{D:Lich} is not
greater than
\[
\oint_Y \langle
P_+\psi_0,\cD_Y P_+\psi_0\rangle \,dv_Y\le \Norm{\cK
\psi_0}{H^{1/2}_*(Y)}^2\;,
\]
and \bref{D:SL} follows immediately, with
\begin{equation}
\label{imprhoeq}
\rho=\max(0,-\tfrac{1}{4}R(g))\;.
\end{equation}
Since $\beta=0$ in this example, the inequality holds for all
 $u\in\hb$ provided $\rho$ satisfies \bref{D-5new}.
 For general mean curvatures
$H_Y\in L^\infty(Y)$, note again that
\[
\oint_Y \langle \psi_0,\cD_Y \psi_0 \rangle \,dv_Y \le \Norm{P_+
\psi_0}{H^{1/2}_*(Y)}^2 \;.
\]
If $H_Y\in L^\infty(Y)$ then
$\oint_Y H_Y|\psi_0|^2 \le \Norm{H_Y}{L^\infty(Y)}
\Norm{J\psi_0}{H^{1/2}_*(Y)}^2$. Using a fractional Sobolev
inequality, the control on $H_Y$ may be weakened to $H_Y\in
L^p(Y)$, $p=n-1$ for $n\ge3$ and $p>1$ for $n=2$.  Lemma
\ref{LemB0} shows that $\Norm{J\psi_0}{H^{1/2}_*(Y)} \le
c \Norm{\tilde{\psi}_0}{H^1(Y\times I')}$, where $\tilde{\psi}=\chi
\psi$ and $\chi=\chi(x)$ is a cutoff function supported in
$I'=[0,\delta/2]$, as in the proof of Theorem \ref{ThmC1}.  Now
Corollary \ref{CorC1} shows that
\[
C^{-1}\Norm{\tilde{\psi}}{H^1_*(Y\times I')}^2
\le  \int_{Y\times I} (|\cD\psi|^2 +|\psi|^2) \,dv_M +
\oint_Y |JP_+\psi_0|^2\,dv_Y,
\]
which provides the required Schr\"odinger-Lichnerowicz estimate
\bref{D:SL}.

In applications, a Schr\"odinger-Lichnerowicz estimate is usually
obtained in the special case of homogeneous boundary data ($\cK
u_0=0$).  The above trick shows that the estimate in the homogeneous
case implies the general estimate \bref{D:SL}:

\begin{Lemma} \label{LemD1a} Under the hypotheses of Lemma
 \ref{Lemd0}, suppose there is $\bar{C}>0$ such that for
  all $u\in\hb$ with $\cK u_0=0$ we have
\begin{equation} \label{D-4b}
\bar{C}^{-1}\int_M |\nabla u|^2\,dv_M
 \le
\int_M(|\cL u|^2+\bar{\rho}|u|^2)\,dv_M\;,
\end{equation}
for some $\bar \rho$. Then there is $C>0$ such that \bref{D:SL} holds
for all $u\in \hb$.
\end{Lemma}

\Proof
Suppose $u\in\hb$ and let $\tilde{u}=u-\chi u$, where $\chi=\chi(x) \in
C^\infty(M) $  is a cutoff function supported in $Y\times I'$ as in the
proof of Theorem \ref{ThmC1}.  Then $\cK \tilde{u}_0=0$ so \bref{D-4b}
applies to $\tilde{u}$, giving
\begin{eqnarray*}
\int_M|\nabla u|^2\,dv_M
&\le&
2\int_M(|\nabla\tilde{u}|^2+|\nabla(\chi u)|^2)\,dv_M
\\
&\le&
C\int_M(|\cL \tilde{u}|^2+\bar{\rho}|\tilde{u}|^2)\,dv_M
 +2\int_{Y\times I'} |\nabla(\chi u)|^2\,dv_M
\\
&\le&
C\int_M(|\cL u|^2+(\bar{\rho}+|d\chi|^2)|u|^2)\,dv_M
 +2\int_{Y\times I'} |\nabla(\chi u)|^2\,dv_M\;.
\end{eqnarray*}
Now it follows easily from Corollary \ref{CorC1} that
\[
C^{-1} \int_{Y\times I'} |\nabla(\chi u)|^2\,dv_M
\le
\int_{Y\times I} (|\cL u|^2 +|u|^2) \,dv_M +
\oint_Y |J\cK u_0|^2\,dv_Y,
\]
which gives the required inequality.
\qed

\begin{Theorem}
\label{ThmD2}
Under the hypotheses \ref{Hh}, suppose $(M,\nabla,\cL,\cK)$ admits a
weighted Poincar\'e inequality \bref{D:WP} and a
Schr\"odinger-Lichnerowicz inequality \bref{D:SL} with $\rho$ and
$\beta$ satisfying the conditions of Lemma~\ref{Lemd0}. If $\rho\in
L^{p}_\loc(M)$ for some $p>n^*/2$,
 and if
\begin{equation}
 \label{D-5}
 \lim_{R\to\infty} \sup_{M\backslash M_R} \frac{\rho}{w} =  0\;,
\end{equation}
 where $\{M_R\}_{R\to\infty}$, is any exhaustion of $M$, then
\begin{equation}
\label{D-6}
(\cL,\cK):\hb \to L^2(M)\times H^{1/2}_*(Y)
\end{equation}
is semi-Fredholm.
\end{Theorem}

\Proof Lemma~\ref{Lemd0} gives $\cL u\in L^2(M)$ for $u\in\hb$.  We
first show the unit ball in the kernel is compact.  Let
$\{u_k\}_{k=1}^\infty$ be a sequence in the kernel of $(\cL,\cK)$,
normalised by $\Norm{u_k}{\hb}=1$. Weak compactness of bounded sets in
$\hb$ shows there is $\bar{u}\in\hb$ and a subsequence, which we also
denote by $u_k$, such that $u_k \rightharpoonup \bar{u}\in\hb$ and
$\Norm{\bar{u}}{\hb}\le \lim\inf \Norm{u_k}{\hb}=1$.

Since \bref{D-5} is independent of the choice of exhaustion, we
may suppose for definiteness that $M_R=\{x\in M: d(x) < R\}$ where
$d(x)$ is the smoothed distance function from some fixed base
point.  Let $\chi\in C^\infty_c(\bbR)$ satisfy $\chi(x)=1$ for
$x\le1$, $\chi(x)=0$ for $x\ge 2$ and $0\le\chi(x)\le1$,
$|\chi'(x)|\le2$ for all $x$.  Then the functions
$\chi_R(x)=\chi(d(x)/R) $ form support functions for the
exhaustion $M_R$ which satisfy $\supp \chi_R \subset M_{2R}$,
$\chi_R =1 $ on $M_R$ and $|d\chi_R| \le 2$. Using the weighted
Poincar\'e inequality we have
\begin{eqnarray*}
\int_M |\nabla(\chi_R u_k)|^2\,dv_M &\le&
  2\int_{M_{2R}\backslash M_R} |d\chi_R|^2 |u_k|^2\,dv_M
+ 2\int_{M_{2R}} |\nabla u_k|^2\,dv_M
\\
&\le&
  2(1+2\sup_{M_{2R}\backslash M_R}w^{-1}) \, \int_M |\nabla u_k|^2\,dv_M\;,
\end{eqnarray*}
which shows that for any $R>1$ the sequence $\chi_R u_k$ is
bounded in $H^1(M_{2R})$.  Since $\chi_Ru_k\rightharpoonup
\chi_R\bar{u}$ in $H^1(M_{2R})$, the
Rellich lemma
implies $\chi_Ru_k\to \chi_R\bar{u}$ strongly  in $L^q(M_{2R})$ for any
$q<\twostar=2n/(n-2)$ and any $R>1$.

Applying \bref{D:SL} to any difference $u_j-u_k$ gives
\begin{eqnarray}
\nonumber
\int_M |\nabla(u_j-u_k)|^2\,dv_M &\le&
  \int_M \rho\,  |u_j-u_k|^2\,dv_M
\\
\nonumber
&\le&  \Norm{\rho}{L^p(M_R)} \Norm{u_j-u_k}{L^q(M_R)}^2
\\
\label{D-7}
&&{}
  + \sup_{M\backslash M_R} \frac{\rho}{w}\ \int_M |u_j-u_k|^2\,w\,dv_M
\end{eqnarray}
where, since $p>n^*/2$, we have $q=2p/(p-1)<\twostar$.
Now \bref{D:WP} and $\Norm{u_k}{\hb}=1$ combine to show that
\[
\int_M |u_j-u_k|^2\,w\,dv_M \le 4,
\]
so by \bref{D-5}, for any $\epsilon>0$ there is $R=R(\epsilon)$ such
that the second term of \bref{D-7} is less than $\epsilon/2$ for all
$j,k$.  Since $u_k$ converges in $L^q(M_R)$  there is $N=N(\epsilon,R)$
such that the first term is less than $\epsilon/2$ for all $j,k\ge N$.
This shows $u_k$ is a Cauchy sequence, hence strongly convergent to
$\bar{u}$, in $\hb$.

As noted above, $\Norm{\cL u}{L^2(M)}\le C\Norm{u}{\hb}$ and thus
\begin{eqnarray*}
\int_M |\cL \bar{u}|^2 \,dv_M & = & \int_M |\cL (\bar{u}-u_k)|^2 \,dv_M
\\
&\le& C \int_M |\nabla(\bar{u}-u_k)|^2\,dv_M
\\
&\to &0 \ \ \textrm{as }k\to\infty\;,
\end{eqnarray*}
which shows that $\cL\bar{u}=0$.  Similarly, since
$\cK:H^{1/2}_*(Y)\to H^{1/2}_*(Y)$ is bounded, for any $u\in \hb$ we have
\begin{eqnarray*}
\oint_Y|J\cK u_0|^2\,dv_Y &\le& c\|\cK u_0\|_{H^{1/2}_*(Y)} \ \le\ ck
\|u_0\|_{H^{1/2}_*(Y)}
\\
&\le& C\int_M |\nabla u|^2\,dv_M \;,
\end{eqnarray*}
by \bref{B:K1} and the trace lemma \ref{LemB0}.  Choosing $u=\bar{u}-u_k$
gives
\[
\oint_Y|J\cK \bar{u}_0|^2\,dv_Y \le c\int_M |\nabla (\bar{u}-u_k)|^2\,dv_M
= o(1),
\]
which shows also that $\cK\bar{u}_0=0$.  Thus $\bar{u}\in
\ker(\cL,\cK)$ and the kernel is finite dimensional.


To show the closed range property, observe that by \bref{D-5} and
\bref{D:WP},  the elliptic estimate \bref{D:SL} may be
strengthened to
\begin{equation}
\label{D-8}
C^{-1} \int_M\left(|\nabla u|^2+|u|^2w\right)\,dv_M
\le
\int_M|\cL u|^2 dv_M
+ \int_{\Omega}\rho|u|^2\,dv_M
+ \oint_Y |J\cK u_0|^2\,dv_Y\;,
\end{equation}
for some relatively compact domain $\Omega\Subset M$.  Now we claim
there is a constant $C>0$ such that
\begin{equation}
\label{D-10}
 \int_\Omega   \rho|u|^2\,dv_M
\le
C\left(
\int_M|\cL u|^2 dv_M
+ \oint_Y |J\cK u_0|^2\,dv_Y
\right)\;,
\end{equation}
for all $u\in\hb$ such that
\begin{equation}
\label{D-9}
\int_M\langle\nabla u,\nabla\phi\rangle\,dv_M =0
\quad\forall \ \phi\in\ker(\cL,\cK)\;.
\end{equation}
Suppose
\bref{D-10} fails, so there is a sequence $u_k\in\hb$, $k=1,2,\ldots$,
such that \bref{D-9} holds for each $u_k$, and
\[
\int_\Omega\rho|u_k|^2\,dv_M=1,\quad
\int_M|\cL u_k|^2 dv_M
+ \oint_Y |J\cK (u_k)_0|^2\,dv_Y
\le 1/k.
\]
The sequence is bounded in $\hb$ by \bref{D-8}, so by passing to a
subsequence we may assume $u_k$ converges weakly to $\bar{u}\in\hb$
and strongly in $L^q(\Omega)$, $q=2p/(p-1)<\hat{2}$ as
before.
Applying \bref{D-8} to $u_j-u_k$ shows the sequence is Cauchy and thus
converges strongly in $\hb$.  It follows that
\[
\int_M|\cL \bar{u}|^2 dv_M
+ \oint_Y |J\cK\bar{u}_0|^2\,dv_Y = 0,
\]
so $\bar{u}\in\ker(\cL,\cK)$.  Strong convergence shows that
\bref{D-9} is also satisfied by $\bar{u}$, so testing \bref{D-9} for
$\bar{u}$ with $\phi=\bar{u}$ shows that $\bar{u}=0$. However, strong
convergence in $L^q(\Omega)$ shows that $\int_\Omega
\rho|\bar{u}|^2dv_M =1$, which is a contradiction and establishes
the claim \bref{D-10}.

Combining \bref{D-10} with \bref{D-8} gives
\begin{equation}
\label{D-11}
\int_M\left(|\nabla u|^2+|u|^2w\right)dv_M
\le
C\left(
\int_M|\cL u|^2 dv_M
+ \oint_Y |J\cK u_0|^2\,dv_Y
\right)
\end{equation}
for all $u\in\hb$ satisfying \bref{D-9}.  Now suppose $u_k\in\hb$
is a sequence such that $\cL u_k=f_k\to f\in L^2(M)$ and
$\cK(u_k)_0=s_k\to s\in H^{1/2}_*(Y)$.  These convergence
properties are retained if we replace $u_k$ by $u_k+y_k$ for any
convergent sequence $y_k\in\ker(\cL,\cK)$, so we may assume the $u_k$
all satisfy \bref{D-9}. In particular, applying \bref{D-11} to
$u_j-u_k$ shows that $u_k$ is Cauchy in $\hb$ and converges to
$\bar{u}$ satisfying $\cL \bar{u}=f$, $\cK \bar{u}_0=s$. This
shows $(\cL,\cK)$ has closed range. \QED

By Definition \ref{defC1}, $u$ is a weak solution of
\begin{equation}
\label{D-12}
\cL u=f,\quad \cK u_0=\sigma\;,
\end{equation} for $f\in L^2(M)$, $\sigma \in PH^{1/2}_*(Y)$,
if $u\in L^2_\loc(M)$ and
\begin{equation}
\label{D-12b}
\int_M \langle u,\cL^\dagger \phi\rangle\,dv_M =
\int_M \langle f,\phi\rangle\,dv_M +
\oint_Y \langle \sigma,\nu\phi_0\rangle\,dv_Y\;,
\end{equation}
for all $\phi\in H^1_c(M)$ such that $\cK^\dag \phi_0=0$.   Similarly, the
argument of Theorem \ref{ThmDI3} shows that the weak form of the adjoint
problem
\begin{equation}
\label{D-12c}
\cL^\dagger u=g,\quad \cK^\dagger u_0=\tau\;,
\end{equation}
for $g\in L^2(M)$, $\tau \in \hat{P}H^{1/2}_*(Y)$,
$\hat{P}=1-\nu^{-1}P\nu$,  is that $u\in L^2_\loc(M)$ and
\begin{equation}
\label{D-12d}
\int_M \langle u,\cL \phi\rangle\,dv_M =
\int_M \langle g,\phi\rangle\,dv_M -
\oint_Y \langle \tau,\nu^{-1}\phi_0\rangle\,dv_Y
\end{equation}
for all $\phi\in H^1_c(M)$ such that $\cK \phi_0=0$.

We now extend the solvability criterion (Fredholm alternative) of
Theorem \ref{ThmDI3} to the non-compact case.

\begin{Theorem}
\label{ThmD3}
Under the conditions of
Theorem~\ref{ThmD2}, suppose the formal adjoint
$(\cL^\dagger,\cK^\dagger)$ also satisfies a
Schr\"odinger-Lichnerowicz estimate \bref{D:SL} with the same
covariant derivative $\nabla$ and with a curvature term $\hat \rho$
satisfying \eq{D-5}.  Then the system \bref{D-12} with $(f,\sigma)\in
L^2(M)\times PH^{1/2}(Y)$ has a solution $u\in\hb$ if and only if
$(f,\sigma)$ satisfies
\begin{equation}
\label{D-14}
\int_M\langle f,\phi\rangle\,dv_M +
\oint_Y\langle\sigma,\nu\phi_0\rangle\,dv_Y
=0\;,
\end{equation}
for all $\phi\in \hb\cap L^2(M)$ satisfying $\cL^\dagger \phi=0$,
$\cK^\dagger\phi_0=0$. In particular, the system \bref{D-12} is
solvable for all $(f,\sigma)\in L^2(M)\times PH^{1/2}(Y)$ if and
only if there are no $0\ne \Psi\in \hb\cap L^2(M)$ satisfying
$\cL^\dagger \Psi=0$, $\cK^\dagger\Psi_0=0$.
\end{Theorem}

\begin{remark} We emphasise that in Theorem~\ref{ThmD3}
it is not necessary to impose conditions
on $\hat \rho$ other than \eq{D-5}, and  no conditions on the map
$\hat \beta:=\cL^\dag-{}^ta^i\nabla_i$ are needed.
\end{remark}

\begin{remark} See Theorem~\ref{pet3} for an example where
$\cL^\dagger\ne\cL$, with $\cL^\dagger$ satisfying two
Lichnerowicz-Schr\"odinger identities with respect to two
\emph{different} connections.
\end{remark}


\Proof The necessity of \bref{D-14} follows immediately from the weak form
\bref{D-12b}.  To show sufficiency, the argument of Theorem \ref{ThmDI3}
applies to reduce to the case $\sigma=0$, which we now consider.

Let $\hb_\cK=\{u\in\hb : \cK
u_0=0\}$.  The elliptic estimate \bref{D:SL} gives
\[
\int_M |\nabla u|^2\,dv_M \le C\int_M(|\cL u|^2
+\rho|u|^2)\,dv_M\;,\quad\forall\ u\in\hb_\cK\;.
\]
The arguments used to show \bref{D-8} and \bref{D-10}  apply and give
\begin{equation}
\int_M |\nabla u|^2\,dv_M \le C\int_M  |\cL u|^2 \,dv_M\quad\forall\
u\in\mathring\hb_\cK\;,
\end{equation}
where we define
\begin{equation}
\mathring\hb_\cK := \{u\in\hb_\cK : \int_M\langle\nabla u,\nabla\phi\rangle
\,dv_M = 0 \ \forall\ \phi\in\ker(\cL,\cK)\}\;.
\end{equation}
Thus the bilinear form $u\mapsto\int_M|\cL u|^2\,dv_M$ is strictly coercive
on the Hilbert space $\mathring\hb_\cK$, and for each $f\in L^2(M)$ the map
$\phi\to\int_M\langle f,\cL\phi\rangle\,dv_M$ is bounded on
$\mathring\hb_\cK$.  The Lax-Milgram lemma shows there is $u\in
\mathring\hb_\cK$ satisfying
\[
\int_M\langle\cL u,\cL\phi\rangle\,dv_M
=
\int_M\langle f,\cL\phi\rangle\,dv_M\quad\forall\ \phi\in \mathring\hb_\cK\;.
\]
Thus setting $\Psi=\cL u-f$ we have
\begin{equation}
\label{D-15}
\int_M\langle \Psi,\cL\phi\rangle\,dv_M=0\quad\forall\ \phi\in \hb_\cK\;,
\end{equation}
since $\phi\in\ker(\cL,\cK)$ will also satisfy the relation
\bref{D-15}. Lemma \ref{Lemd0} shows that $\Psi\in L^2(M)$ and from
Definition \ref{defC1} and \bref{D-15} we see that $\Psi$ is
a weak solution of
\[
\cL^\dagger \Psi =0,\quad \cK^\dagger\Psi =0\;.
\]
If there are no such non-trivial $\Psi$ then $\cL u=f$, and $u$ is the
required solution.  The arguments of Theorem \ref{ThmDI3} show that
$(\cL^\dagger,\cK^\dagger)$ is elliptic and Theorem \ref{ThmC1}
applies to show $\Psi\in H^1_\loc(M)$.  Let $M_R$ be the exhaustion of
$M$ constructed in Theorem \ref{ThmD2}, with associated cutoff
functions $\chi_R\in C^\infty_c(M)$, and let $\Psi_k=\chi_k\Psi\in
H^1_c(M)\subset \hb$.  The assumed Schr\"odinger-Lichnerowicz estimate
\bref{D:SL} for $(\cL^\dagger,\cK^\dagger)$ gives (with $\cL^\dagger$
curvature term $\hat{\rho}$)
\begin{equation}
\label{D-16}
\int_M |\nabla(
\Psi_k-\Psi_l)|^2\,dv_M \le
C \int_M \left(|\cL^\dagger(\Psi_k-\Psi_l)|^2 +
\hat{\rho}|\Psi_k-\Psi_l|^2\right)\,dv_M\;.
\end{equation}
Since $\cL^\dagger\Psi=0$ we have
\[
\int_M |\cL^\dagger(\Psi_k-\Psi_l)|^2\,dv_M \le
c \int_M \left(|d\chi_k|^2|\Psi|^2 + |d\chi_l|^2|\Psi|^2\right) \,dv_M \to 0
\;,\]
because $\Psi\in L^2(M)$, $|d\chi_k|\le 2$ and $\supp d\chi_k\subset
M_{2k}\backslash M_k$.   Now
\[
\int_M
\hat{\rho}|\Psi_k-\Psi_l|^2\,dv_M \le \epsilon
\int_M
|\Psi_k-\Psi_l|^2\,w\,dv_M
\]
by the condition \bref{D-5} on $\hat \rho$, for sufficiently large
$k,l$.  By the weighted Poincar\'e inequality \bref{D:WP}, this is in
turn bounded by $\epsilon$ times the left side of \bref{D-16} and may
therefore be discarded in \eq{D-16} by choosing $\epsilon$
sufficiently small.  It follows that $\Psi_k$ is a Cauchy sequence in
$\hb$, so $\Psi\in\hb\cap L^2$ and thus $\cL^\dagger\Psi=0$,
$\cK^\dagger\Psi_0=0$.  If there is no such $\Psi\ne 0$ then $\cL u =
f$, and $u$ is the required solution.  More generally we have $\cL
u=f+\Psi$, $u\in \mathring{\hb}_\cK$, and since $\int_M\langle\cL
u,\Psi\rangle\,dv_M=0$ by \eq{D-15}, the condition \bref{D-14} (with
$\sigma=0$ and $\phi=\Psi$) shows that $\Psi=0$ and we have solved
$\cL u=f$, as required.  \QED



\section{Weighted Poincar\'e Inequalities}
\label{secF}
%
%
 Before proceeding with the analysis, define the \emph{symmetric part} $\Gamma^S$
of the connection $\nabla$  by the formula
\begin{equation}
\label{symGamdef}
\langle\phi,\Gamma^S(X)\psi\rangle := \frac
12\Big(X\langle\phi,\psi\rangle -
\langle\phi,\nabla_X\psi\rangle-\langle\nabla_X\phi,\psi\rangle\Big)
\;,
\end{equation}
for all smooth sections $\phi,\psi$ of $E$ and all smooth vector
fields $X$: One easily checks that \eq{symGamdef} defines a linear map
$\Gamma^S(X)$ from fibers of $E$ to fibers of $E$, symmetric with
respect to the scalar product $\langle\cdot,\cdot\rangle$, with the
map $X\to\Gamma^S(X)$ being linear as well. Clearly, $\nabla$ is
compatible with $\langle\cdot,\cdot\rangle$ if and only if $\Gamma^S$
vanishes.  If $\Gamma$ is defined by \eq{coneq}, then
$$\Gamma^S=\half(\Gamma+{}^t\Gamma)\;.$$
 We establish Proposition~\ref{PropD1} via a special
case, based on  an argument of Geroch--Perng
\cite{GerochPerng94}:

\begin{Lemma}
\label{LemF1}
Let $\Omega,\tilde{\Omega}$ be any two relatively
compact domains in $M$, and assume that
\bel{F-1.0}\Gamma^S \in L^{n^*}_\loc(M)\;.\ee
  There is a constant $\epsilon>0$ such
that for all sections $u\in H^1_{\loc}(M)$ of $E$ we have
\begin{equation} \label{F-1}
\epsilon \int_{\tilde{\Omega}}|u|^2\,dv_M \le \int_\Omega|u|^2\,dv_M
+\int_M|\nabla u|^2\,dv_M\;.
\end{equation}
\end{Lemma}

\Proof
Let $q$ be any point of $\tilde{\Omega}$, fix $p\in\Omega$ and let
$r_p$ be small enough that the $\go$-geodesic ball $B(p,r_p)$ of radius
$r_p$ and centred at $p$, lies within $\Omega$.  Let $X$ be a
$C^\infty$ compactly supported vector field, such that the associated
flow $\phi_t$ satisfies $\phi_1(B(p,r_p))\supset B(q,r_q)$ for some
$r_q>0$.  (Since $M$ is $C^\infty$ and connected, it is always
possible to construct such an $X$.)  Let $\Omega_t=\phi_t(B(p,r_p))$.

By direct calculation and H\"older's inequality we have, for any $u\in
H^1_{\loc}(M)$,
\begin{eqnarray*}
\frac{d}{dt}\int_{\Omega_t}|u|^2 dv_M
&=& \int_{\Omega_t}\left( 2\langle u,(\nabla_X+\Gamma^S_X)u\rangle +|u|^2
\mathrm{div}_{\go} X\right) \,dv_M
\\
&\le& C\left(\int_{\Omega_t}\left(|u|^2 +|\nabla u|^2\right)\,dv_{\go}
+\Norm{\Gamma^S}{L^{n^*/2}(\Omega_t)} \Norm{u}{L^{\hat{2}}(\Omega_t)}^2\right)\;,
\end{eqnarray*}
where $C$ depends on $\Norm{X}{L^\infty}$,
$\Norm{\mathrm{div}_{\go}X}{L^\infty}$. By the Sobolev inequality in
the coordinate ball $\Omega_t$ for functions,
$\Norm{f}{L^{\twostar}(\Omega_t)}\le C(\Norm{\partial
  f}{L^2(\Omega_t)} + \Norm{f}{L^2(\Omega_t)})$.  Applying this to
$f=|u|$ gives $\Norm{u}{L^{\twostar}(\Omega_t)}\le C(\Norm{D
  u}{L^2(\Omega_t)} + \Norm{u}{L^2(\Omega_t)})$, where $D$ is any
metric-compatible connection. Since $\Gamma^S\in L^{n^*}$ may be
written as $\Gamma_1+\Gamma_2$, $\Gamma_1\in L^\infty$,
$\Norm{\Gamma_2}{L^{n^*}}\le\epsilon$, the Sobolev inequality gives
\[
\Norm{u}{L^{\twostar}(\Omega_t)}\le C\left(\Norm{\nabla u}{L^2(\Omega_t)} +
\Norm{u}{L^2(\Omega_t)}\right)\;,
\]
for some constant $C$ depending on $\Gamma$. Defining
$F(t)=\int_{\Omega_t}|u|^2dv_{\go}$, we have
\[
\frac{d}{dt}F(t) \le CF(t)+ C\int_M|\nabla u|^2dv_M\;,
\]
and Gronwall's lemma gives $F(1)\le e^C(F(0)+\int_M|\nabla u|^2dv_M)$.
Thus there is $\epsilon>0$ such that
\[
\epsilon\int_{B(q,r_q)}|u|^2\,dv_M \le \int_\Omega |u|^2\, dv_M +
\int_M |\nabla u|^2\,dv_M\;.
\]
Since $\tilde{\Omega}$ has compact closure, it is covered by finitely
many such balls $B(q,r_q)$ and \bref{F-1} follows.  \QED

\begin{Corollary}
\label{CorF2}
Under condition~\eq{F-1.0}, if there is a domain $\Omega\subset M$ and
a constant $\epsilon>0$ such that
\begin{equation}
\label{F-2}
\epsilon\int_\Omega |u|^2\,dv_M \le \int_M |\nabla u|^2\,dv_M
\end{equation}
for all $u\in C^1_c(M)$, then $M$ admits a weighted Poincar\'e
inequality \bref{D:WP}.
\end{Corollary}

\Proof
By paracompactness and Lemma \ref{LemF1}, there is a countable locally
finite covering of $M$ by domains $\Omega_k$ and constants
$1\ge\epsilon_k>0$, $k\in\bbZ^+$, such that for each $k$,
\[
\epsilon_k\int_{\Omega_k} |u|^2\,dv_M \le \int_\Omega
|u|^2\,dv_M+\int_M |\nabla u|^2\,dv_M\;.
\]
This is in turn bounded uniformly by  \bref{F-2}, so the function
\begin{equation}
\label{F-3}
w(x) = \sum_{k:x\in\Omega_k}\frac{2^{-k}\epsilon\epsilon_k}{1+\epsilon}
\end{equation}
is bounded, strictly positive, and satisfies
\[
\int_M |u|^2\,w\,dv_M \le \int_M |\nabla u|^2\,dv_M\;,
\]
which is the required weighted Poincar\'e inequality.
\QED
\medskip

This establishes part (i) of Proposition~\ref{PropD1}, and we next turn
to the proof of part (ii).

\begin{Theorem}
\label{ThmF2}
Suppose that $M$ has a locally finite cover such that
\bel{F-3a2}\nabla_i = \partial_i - \Gamma_i\;, \qquad \mbox{with } \
\Gamma_i \in L^{n^*}_\loc\;. \ee If there are no global
$\nabla$-parallel sections of the bundle $E$, then $M$ admits a
weighted Poincar\'e inequality. Equivalently, if $M$ does not admit a
weighted Poincar\'e inequality then $M$ admits a global
$\nabla$-parallel section.
\end{Theorem}

\Proof
 Assume $M$ does not admit a weighted Poincar\'e inequality, so by
Corollary \ref{CorF2}, for each domain $\Omega\Subset M$ and each constant
$\epsilon>0$, there is $u\in H^1_{\loc}(M)$ such that \bref{F-2} fails.  In
particular, fixing $\Omega$, for each $k>0$ there is $u_k\in H^1_{\loc}(M)$
such that
\begin{equation}
  \label{F-3a}
\int_\Omega |u_k|^2\, dv_M =1,\quad \int_M |\nabla u_k|^2\,dv_M \le k^{-1}\;.
\end{equation}
It follows that $\nabla u_k \to 0$ strongly in $L^2(M)$. Under
\eq{F-3a2} Rellich's lemma holds, so  there is a subsequence
converging strongly to $u \in L^2(\Omega)$.  Then $\nabla u =0$ and
$u\ne 0$ in $\Omega$.

 Now let $M_j$, $j=1,2,\dots$ be the
exhaustion of $M$ from Theorem \ref{ThmD2}, and let $u_j\in
H^1(M_j)$ be the corresponding parallel spinors, constructed in
the preceding paragraph.  Since $u_j\ne0$ there is $M_j'\Subset
M_j$ such that $\int_{M_j'}|u_j|^2 \ne0$.  Lemma~\ref{LemF1}
applied with $M_j$ replacing $M$ shows there is $\eta_j>0$ such
that for all $v\in H^1_{\loc}(M_j)$,
\[
\eta_j \int_{M_j'}|v|^2\,dv_M
\le \int_{M_0}|v|^2\,dv_M  + \int_{M_j}|\nabla v|^2\,dv_M\;.
\]
In particular this implies $\int_{M_1}|u_j|^2\,dv_M\ne 0$ and we
may impose the normalisation
$\int_{M_1}|u_j|^2\,dv_M=1$. By Rellich's lemma there is
$\bar{u}_1\in H^1(M_1)$ and a subsequence, also denoted by $u_j$,
such that $u_j\to \bar{u}_1$ in $H^1(M_1)$ and
$\int_{M_1}|\bar{u}_1|^2=1$, $\nabla\bar{u}_1=0$.

Again by Lemma~\ref{LemF1}, for each $k\ge1$ there is $\epsilon_k>0$ such
that
\[
\epsilon_k \int_{M_k}|v|^2\,dv_M
\le \int_{M_1}|v|^2\,dv_M  + \int_{M_{k+1}}|\nabla v|^2\,dv_M \;,
\quad\forall\ v\in H^1_{\loc}(M_{k+1})\;.
\]
Setting $v=u_i-u_j$, $i,j>k$, shows that the sequence $u_j$ is Cauchy in
 $L^2(M_k)$ and therefore converges strongly in $L^2(M_k)$
for all $k\ge1$ to some nontrivial $\bar{u}\in
L^2_{\loc}(M)$, and $\nabla\bar{u}=0$.
\QED

Another application of Corollary \ref{CorF2} leads to Proposition
\ref{PropD1} part~3, for asymptotically flat manifolds.  In fact the proof
works for a much broader class of manifolds:
\begin{Definition}
\label{defF1} A \emph{weakly asymptotically flat end} $\tM \subset
M$ of a Riemannian manifold $M$ with metric $g$ is a connected
component of $M\backslash K$ for some compact set $K$, such that
$\tM \simeq \bR^n\backslash B(0,1)$ and there is a constant
$\eta>0$ such that
\[
\eta \,\delta_{ij}\xi^i \xi^j \le g_{ij}(x)\xi^i \xi^j \le
\eta^{-1}\delta_{ij}\xi^i \xi^j\;;
\]
for all $x\in \bR^n\backslash B(0,1)$ and all vectors $\xi\in \bR^n$.
\end{Definition}

\begin{Theorem}
\label{ThmF3}
Suppose $(M,g)$ is a (connected) Riemannian manifold of dimension $n\ge 3$,
$g\in C^0(M)$, and $M$ has a weakly asymptotically flat end $\tM $.
Suppose also the connection $\nabla_i=\partial_i-\Gamma_i$ on $E$ satisfies
$\Gamma\in L^{n^*}_\loc(M)$ and the decay conditions
\begin{equation}
\label{F-4a}
 \Norm{r^{-1}\Gamma^S}{L^{n/2}(\tM )}
+  \Norm{\Gamma^S}{L^{n}(\tM )}
< \infty\;,
\end{equation}
where $\Gamma^S$ is the symmetric,
scalar product incompatible, component of  $\nabla$ defined by \Eq{symGamdef}.
Then $M$ admits a weighted Poincar\'e inequality.
\end{Theorem}
\begin{remark} The restriction $\dim M\ge 3$ is rather
harmless as far as the applications to the positive mass theorems are
concerned, since the notion of asymptotic flatness for two dimensional
manifolds, relevant to general relativistic applications, has to be
defined in a completely different way. An adequate
analogue of mass here when $\dim M =2 $ is provided by the Shiohama
theorem~\cite{Shiohama}.\end{remark}

\begin{remark} The decay condition \bref{F-4a} is independent of the choice of
flat background metric $\mathring{g}_{ij}=\delta_{ij}$: \Eq{symGamdef}
shows that $\Gamma^S$ is a tensor. By comparison with the $g$-distance
function from any chosen point $p$, the function $r$ is equivalent to
this distance function, which implies the result.
\end{remark}
\begin{remark} The proof below establishes the inequality \eq{F-2} for
  spinors supported in $\Omega:=\R^3\setminus B(0,R)$ for some $R$
  without assuming that $\Gamma\in L^{n^*}_\loc(M)$.
\end{remark}

\Proof Let $r=(\sum (x^i)^2)^{1/2}\in C^\infty(\tM )$ and
$\chi=\chi(r)\in C^1_c(\tM )$ satisfy, for some $R_0>1$ and $k\ge10$,
\[
  \chi(r) = \frac{\log(r/R_0)}{\log k},\quad 2 R_0 \le r \le (k-1) R_0
\]
and $\chi(r)=1$ for $r>kR_0$, $\chi(r)=0$ for $r\le R_0$.  Then
$|\chi'(r)| \le 2/(r\log k)$, so for any section $u\in C_c^1(M)$
\begin{equation}
\label{F-4}
  \int_M|\nabla (\chi u)|^2 \,dv_M \le 2\int_M |\nabla u|^2\,dv_M
  +\frac{4}{(\log k)^2} \int_{R_0\le r\le kR_0}\frac{1}{r^2}|u|^2 \,dv_M\;.
\end{equation}
Now $\Delta_0(r^{2-n})=0$ for $r\ge1$ in $\bR^n$, $n\ge 3$, so for any
$v\in C_c^1(\bR^n\backslash B(0,R_0))$ we have
\begin{eqnarray*}
0 &=& -\int_{\bR^n} \partial_i(\partial_i(r^{2-n})\,|v|^2r^{n-2})\,dx
\\
&=&(n-2)^2 \int_{\bR^n}r^{-2}|v|^2\,dx +
(n-2)\int_{\bR^n}r^{-1}2\langle v,(\nabla_r+\Gamma^S_r)v\rangle\,dx\;,
\end{eqnarray*}
where $\Gamma^S_r = r^{-1}x^i\Gamma^S_i$ and
 lengths are measured by $\mathring{g}$ and the metric on $E$.
Using H\"older's inequality we obtain
\[
\frac{(n-2)^2}{4} \int_{\bR^n}r^{-2}|v|^2\,dx
\le
\int_{\bR^n}|\nabla v|^2\,dx
+ (n-2)\int_{\bR^n} r^{-1}|v|^2 |\Gamma^S_r|\,dx\;.
\]
The Sobolev inequality in $\bR^n$, $n\ge 3$,
\[
\left(\int_{\bR^n} |v|^{\hat{2}}\,dx\right)^{1-2/n}
\le C_S \int_{\bR^n} |Dv|^2\,dx\;,
\]
where $D=\nabla+\Gamma^S$ is the metric-compatible connection, gives the estimate
\begin{eqnarray*}
\int_{\bR^n} |Dv|^2\,dx &\le& 2\int_{\bR^n} (|\nabla
v|^2+|v|^2|\Gamma^S|^2) \,dx
\\
&\le&
 2\int_{\bR^n} |\nabla v|^2 \,dx
+2\,C_S\,\Norm{\Gamma^S}{L^n(\bR^n\backslash B(0,R_0))} \int_{\bR^n} |Dv|^2 \,dx\\
&\le&
4 \int_{\bR^n} |\nabla v|^2\;,
\end{eqnarray*}
provided $2C_S\Norm{\Gamma^S}{L^n(\bR^n\backslash B(0,R_0))}\le\half$.
Now \bref{F-4a} implies there is $R_0<\infty$ such that this condition
will be satisfied, so for any $v\in C_c^1(\bR^n\backslash B(0,R_0))$
we have
\begin{eqnarray*}
\int_{\bR^n} r^{-1}|\Gamma^S_r|\,|v|^2\,dx
&\le&
\Norm{r^{-1}\Gamma^S}{L^{n/2}(\bR^n)}\,C_S\,\int_{\bR^n}|Dv|^2\,dx
\\
&\le&
4C_S\Norm{r^{-1}\Gamma^S}{L^{n/2}(\bR^n)} \int_{\bR^n}|\nabla v|^2\,dx\;.
\end{eqnarray*}
Hence there is  $\epsilon>0$ such that for all $v\in
C^1_c(\tM \cap\{r>R_0\})$,
\begin{equation}
\label{F-5}
\epsilon\int_{\tM } r^{-2}|v|^2\,dv_M
\le
\int_{\tM } |\nabla v|^2\,dv_M\;.
\end{equation}
Combining \bref{F-5} with $v=\chi u$ and \bref{F-4} gives
\begin{eqnarray*}
\int_{\{r>kR_0\}} r^{-2}|u|^2\,dv_M &\le &
\int_{\tM }  r^{-2}|\chi u|^2\,dv_M
\\
&\le&
C \int_{\tM } |\nabla(\chi u)|^2\,dv_M
\\
&\le&
C \int_{\tM } |\nabla u|^2\,dv_M
+ \frac{C}{(\log k)^2}\int_{\tM }  r^{-2}| u|^2\,dv_M \;,
\end{eqnarray*}
where now $|\nabla u|^2=g^{ij}\langle \nabla_iu,\nabla_ju\rangle$.
If $k$ is chosen so that $C/(\log k)^2 \le\half$ then the last term
may be absorbed into the left hand side, giving
\begin{equation}
\label{F-6}
 \int_{r\ge kR_0}  r^{-2}|u|^2\,dv_M
\le C \int_{M} |\nabla u|^2\,dv_M\;.
\end{equation}
Lemma \ref{LemF1} now applies and gives the required weighted Poincar\'e
inequality.
\QED

In order to prove part 4.~of Proposition~\ref{PropD1} the
following Definition is needed:

\begin{Definition}
\label{defF1.&} A \emph{weakly hyperboloidal end} $\tM \subset M$
is a connected component of $M\backslash K$ for some compact set
$K$, such that $\tM \simeq (0,x_0)\times \cN$, where $(\cN,\nmet)$
is a (boundaryless) compact Riemannian manifold with continuous
metric $\nmet$, with $g|_{\tilde M}$ being uniformly equivalent to
$$ \mathring{g}\equiv x^{-2} (dx^2 + \nmet)\;.$$ Here $x$ is the
coordinate running along the $(0,x_0)$ factor of $(0,x_0)\times
\cN$.
\end{Definition}

We have the following hyperboloidal counterpart of Theorem~\ref{ThmF3}:
\begin{Theorem}
\label{ThmF3.1} Suppose $(M,g)$ is a (connected) Riemannian
manifold of dimension $n\ge 2$, $g\in C^0(M)$, and $M$ has a
weakly hyperboloidal end $\tM $. Suppose also the connection
$\nabla_i=\partial_i-\Gamma_i$ on $E$ satisfies $\Gamma\in
L^{n^*}_\loc(M)$ and the decay condition
\begin{equation}
\label{F-4b} \limsup_{x\to 0} |x\Gamma^S_x| < \frac{n-1}{2}
\end{equation}
in $\tM $, where $\Gamma^S_x$ is the symmetric part
of $\nabla_{\partial _x}$, with norm understood as that of an endomorphism of
 fibres of $E$. Then $M$ admits a weighted
Poincar\'e inequality.
\end{Theorem}

\proof This is essentially McKean's inequality \cite{McKean}; we
follow the proof in \cite{ChSob}. Let, first, $f$ be a function in
$ C^1([0,x_0]\times \cN)$ with $f=0$ at $\{x=0\}$; we have
\begin{eqnarray}
\nonumber
f^2(x,v) & = & 2\int_0^x f(s,v)\frac{\partial f(x,v)}{\partial x} ds
\\ \label{hwpi1}
&\le & \frac{n-1}2\int_0^x \frac{f^2(s,v)}{s}ds + \frac{2}{n-1}\int_0^x s
\left(\frac{\partial f}{\partial x}(s,v)\right)^2 ds \;.
\end{eqnarray}
Here we use the symbol $v$ to label points in $ \cN$.  Integrating on
$ [0,x_0]\times \cN $, 
a change of the order of integration in $x$ and $s$
together with some obvious manipulations gives
\begin{eqnarray}
\nonumber \int_{[0,x_0]\times \cN}f^2\;x^{-n}dx\,d\mu_h &\le &
\frac{4}{(n-1)^2}\int_{[0,x_0]\times \cN}
\left(x\frac{\partial f}{\partial x}\right)^2 \;x^{-n}dx\,d\mu_h\\
&\le & \frac{4}{(n-1)^2}\int_{[0,x_0]\times \cN}
\mathring{g}(df,df) \;x^{-n}dx\,d\mu_h \;.\label{wpif}
\end{eqnarray}
 This is the desired inequality on $\tilde M$
with metric $\mathring g$ for functions, with weight function $w=
(n-1)^2/4$. The result for general weakly asymptotically
hyperboloidal metrics and for functions follows immediately from
the above, using uniform equivalence of $g$ with $\mathring{g}$ on
the asymptotic region, and using Lemma~\ref{LemF1}.

Let, finally, $v$ be a smooth compactly supported section of a
Riemannian bundle with not-necessarily-compatible
connection $\nabla$. Let $\phi$ be any
smooth compactly supported function equal to $1$ on the support of
$v$, set
$$f_\epsilon = \phi \sqrt{\epsilon + \langle v,v\rangle}\;.
$$
We have
\begin{eqnarray*}
\left  |\frac{\partial f_\epsilon}{\partial x}\right|^2 & = &
|d_x\phi|^2 (\epsilon + \langle v,v\rangle) +
  \phi^2 \frac{ \langle v,(\nabla_x +\Gamma^S_x)v\rangle\langle
    v,(\nabla_x+\Gamma^S_x)v\rangle}{\epsilon + \langle v,v\rangle}
\\ & \le & \epsilon |d\phi|^2 + \phi^2 |(\nabla_x + \Gamma^S_x) v|^2\;.
\end{eqnarray*}
The first line of \eq{wpif} yields
\begin{eqnarray*}
 \int_M f_\epsilon^2 \;x^{-n}dx\,d\mu_h& = & \int_M \phi^2 (\epsilon + \langle
 v,v\rangle) \;x^{-n}dx\,d\mu_h
\\ & \le & \frac{4}{(n-1)^2}\int_M x^2\left(\epsilon |d\phi|^2 + \phi^2
|(\nabla_x+\Gamma^S_x)
  v|^2\right)  \;x^{-n}dx\,d\mu_h\;.
\end{eqnarray*}
Passing with $\epsilon$ to zero gives
\begin{eqnarray*}
 \int_M \langle v, v \rangle  \;x^{-n} dx\,d\mu_h&  \le & \frac{4}{(n-1)^2}\int_M
 x^2
|(\nabla_x+\Gamma^S_x)
  v|^2 \;x^{-n}dx\,d\mu_h\
  \\
  &  \le & \frac{4}{(n-1)^2}\int_M
\left((1+\frac 1 \delta)|\nabla
  v|_{\mathring{g}}^2 + x^2(1+\delta)|\Gamma^S_x
  v|^2dx\right)\;x^{-n}dx\,d\mu_h\;,
\end{eqnarray*}
for any $\delta>0$, and if condition~\eq{F-4b} holds the last term
can be carried over to the left hand side, leading to
\begin{eqnarray*}
 {C^{-1}}\int_M \langle v, v \rangle  dv_M\ & \le & \int_M   |\nabla
  v|_{{g}}^2 dv_M\;.
\end{eqnarray*}
Lemma~\ref{LemF1} gives then the desired inequality, with a weight
function $w$ equal to $1/C$ in the asymptotic region. \qed

\section{Examples and Applications}
\label{secE}
%
%
The structure and regularity conditions may be readily verified in
situations of interest, which we illustrate by considering the Dirac
operator examples of \S\ref{secA}.

Suppose $M$ is a Riemannian spin manifold.  Fix local coordinates
$(x^\mu)$ and a local orthonormal framing $e_i=e_i^\mu\partial_\mu$ of
the tangent bundle $TM$, and let $\phi_I$, $I=1,\dots,\dim S$ denote
an associated spinor frame, determined by some choice of
representation $c:C\ell_n\to \mathrm{End}(S)$.  The Dirac
operator $\cD$ defined by \bref{A:gpsi}-\bref{A:dirac} is
\begin{equation}
\label{E-1}
\cD = e_i^\mu\gamma^i\partial_\mu -
\tfrac{1}{4}\omega_{ij}(\partial_\mu)e_k^\mu \gamma^k\gamma^i\gamma^j\;,
\end{equation}
where the skew-symmetric matrices $\gamma^i=c(e_i)\in
\mathrm{End}(S)$ are constant in the local spinor frame and satisfy
the Clifford relation
\begin{equation}
\gamma^i\gamma^j + \gamma^j\gamma^i = -2\delta^{ij}\;.
\end{equation}
Clearly $\cD$ has the form
$a^\mu\partial_\mu + b$ where
\begin{equation}
\label{E-3}
\begin{array}{rcl}
a^\mu &=& \displaystyle{\sum_{i=1}^n} e_i^\mu \gamma^i\;, \\[2pt]
b &=& {}-\tfrac{1}{4}
 \displaystyle{\sum_{i,j,k,\mu=1}^n} \omega_{ij}(\partial_\mu)e_k^\mu
    \gamma^k\gamma^i\gamma^j\;.
\end{array}
\end{equation}
If $g_{\mu\nu}\in W^{1,n^*}_\loc \cap C^0_\loc$ then by the
Gram-Schmidt construction,
the local orthonormal frame may be chosen so that the
coefficients also satisfy $e_i^\mu\in W^{1,n^*}_\loc \cap
C^0_\loc$,  and then $a^\mu,b$ also satisfy \bref{IR-2},
{\em cf.} Proposition~\ref{pspinnew}, Appendix~\ref{AppAGC}.
The Clifford identity shows $|\xi_\mu a^\mu
V|^2=g^{\mu\nu}\xi_\mu\xi_\nu|V|^2$ pointwise, which implies
\bref{IR-3} so Theorem \ref{ThmIR5} may be applied to establish the
$H^1_\loc$ interior regularity of $L^2$ weak solutions.

For boundary regularity we assume there is a diffeomorphism of
$Y\times I$ with a neighbourhood of $Y\times \{0\}=\partial M$ in $M$,
such that for any chart of $Y$, the associated chart of adapted
coordinates $(y^A,x)$, $x\ge0$ satisfies
\begin{equation}
\label{E-4}
\begin{array}{rcl}
g_{xA}(y,0) &=& 0 , \ \ A=1,\ldots,n-1\;, \\
g_{xx}(y,0) &=& 1\;, \\
\end{array}
\end{equation}
and
\begin{equation}
\label{E-5}
g_{\mu\nu}\rY \in W^{1,(n-1)^*}(Y)\cap C^0(Y) \;.
\end{equation}
(Recall that the symbol ``$|_Y$'' stands for ``$|_{Y\times
\{0\}}$''.)

For smooth metrics this follows easily using Gaussian coordinates
about $Y$. Geodesic uniqueness may not be available in the more
general case of a $W^{k+1,p}$ manifold with metric $g_{\mu\nu}\in
W^{k,p}$ for $k>n/p$; in this case the existence of boundary
coordinates satisfying \bref{E-4} is guaranteed by
Proposition~\ref{PAGC.1}.

Then \bref{C-4d} follows from the regularity conditions on
$g_{\mu\nu}$, \bref{C-2} follows from \bref{E-5}, and \bref{C-6}
follows from \bref{E-4}, the Clifford relations and the skew-symmetry
of the $\gamma^i$. In terms of an adapted frame $e_i$, where $e_n\rY =
-\partial_{x^n}$ is the outer normal at $Y=\partial M$, the boundary
operator \bref{C-7} may be taken as
\begin{equation}
  \label{E-6}
  A=\cDb = -\sum_{i=1}^{n-1} e^\mu_i \gamma^n\gamma^i\partial_\mu
  +\tfrac{1}{4} \sum_{i,j,k=1}^{n-1}\omega_{ij}(\partial_\mu) e_k^\mu
  \gamma^n\gamma^k\gamma^i\gamma^j\;,
\end{equation}
by an appropriate choice of $\tilde{b}_0$ in \bref{C-7}.  Again note
that other choices of $A$ are possible, such as $\cDb + F$ for any
function or symmetric endomorphism $F$.

The boundary condition \bref{A:herzlich}
\begin{equation}
  \label{E-7}
  P_+\psi = \sigma, \qquad\sigma\in P_+H^{1/2}_*(Y)\;,
\end{equation}
where $P_+$ is the orthogonal projection to the positive eigenspace of
$A$, corresponds via \bref{C-8o} to $K=0$, which clearly satisfies the
conditions of Theorem \ref{ThmC1} for $K,K^\dagger$.  From Theorem
\ref{ThmC1} it follows that the weak-strong property and the elliptic
estimate \bref{C-9b} hold at the boundary for the boundary condition
\bref{E-7}.

Boundary operators of the type \bref{GHHP} were used in
\cite{GHHP83,LudvigsenVickers82,DouganMason91}, for example.  More
generally suppose there is an endomorphism $\epsilon:S\to S$ acting on
sections of $E$, which satisfies
\begin{equation}
\label{E-8}
\epsilon^2=1,\quad {^t}\epsilon = \epsilon,\quad \epsilon A + A\epsilon=0\;.
\end{equation}
Note that if there is a splitting $E=E_0\oplus F_0$ such that $A$ has
the form \bref{SC-10} then $\epsilon =
\left[\begin{array}{cc}I&0\\0&-I\end{array}\right]$ will satisfy
\bref{E-8}.  Conversely, if \bref{E-8} holds then setting $E_0,F_0$
equal to the $\pm1$ eigenspaces of $\epsilon$ shows that $A$ may be
written in the supersymmetric form \bref{SC-10}.

If $A\psi=\lambda\psi$ then $A(\epsilon\psi)=-\epsilon A \psi =
-\lambda\epsilon\psi$, so the spectrum of $A$ is symmetric and
$\epsilon$ interchanges the positive and negative eigenspaces.  Let
$P_+$ ($P_-$) be orthogonal projection to the positive (negative)
eigenspaces, and consider the eigenspace splitting $L^2(Y)=H_+\oplus
H_- \oplus H_0^+ \oplus H_0^-$, where $H_0=\ker A$ and $H_0= H_0^+
\oplus H_0^-$ is the decomposition into $\pm1$ eigenspaces of
$\epsilon$.  There is an isometric isomorphism $e:H_+\to H_-$ such
that $\epsilon$ has the block decomposition
\[ \epsilon =
\left[\begin{array}{cccc}&e^{-1}&&\\e&&&\\&&1&\\&&&-1\end{array}\right]
\textrm{,  acting on }=H_+\oplus H_- \oplus  H_0^+ \oplus
H_0^-.
\]
It follows that the action of $\cBK_+ = \half(1+\epsilon)$ is given by
\begin{equation}
\cBK_+\psi =  \half(1+\epsilon)
\left[\begin{array}{c}p\\q\\r^+\\r^-\end{array}\right]
=
\left[\begin{array}{c}
\half(p+e^{-1}q)\\
\half(ep+q)\\
r^+\\
0
\end{array}\right]\;.
\end{equation}
Hence the boundary condition $\cBK_+\psi=\sigma$ where
$\sigma=\epsilon\sigma = {}^t[\sigma_+,\ e\sigma_+,\ \sigma_0^+,\ 0]$, is
equivalent to the component conditions $\half(p+e^{-1}q) = \sigma_+$
and $r^+=\sigma_0^+$.  This may be expressed in the form \bref{C-8o} if
we define $P,K$ by
\begin{equation}
P=P_++P_0^+ =
\left[\begin{array}{cccc}1&&&\\&0&&\\&&1&\\&&&0\end{array}\right],
\quad
K =
\left[\begin{array}{cccc}&e^{-1}&&\\0&&&\\&&0&\\&&&0\end{array}\right]\;,
\end{equation}
since then
\begin{equation}
P\psi-K(1-P)\psi =
\left[\begin{array}{c}p\\0\\r^+\\0\end{array}\right]
+
\left[\begin{array}{c}e^{-1}q\\0\\0\\0\end{array}\right]
=
\left[\begin{array}{c}2\sigma_+\\0\\\sigma_0^+\\0\end{array}\right]
\end{equation}
is equivalent to $\cBK_+\psi=\sigma$ as above.  Clearly $K$ is bounded
on both $L^2(Y)$ and $H^{1/2}_*(Y)$ as required by the regularity
theorems of \S\ref{secB}.  Similarly we find that the boundary
condition $\cBK_-\psi=\sigma$ is equivalent to $\half(p-e^{-1}q) =
\sigma_-$ and $r^-=\sigma_0^-$, so an equivalent spectral projection
condition may be constructed using $\hat{P}=P_++P_0^+$ and
$\hat{K}=-K$, so
\begin{equation}
\hat{P}\psi-\hat{K}(1-\hat{P})\psi =
\left[\begin{array}{c}p-e^{-1}q\\0\\0\\r^-\end{array}\right]
=
\left[\begin{array}{c}2\sigma_+\\0\\0\\\sigma_0^-\end{array}\right]\;.
\end{equation}

\section{Positive mass theorems}
\label{Section:pmt}

Our motivation for the present work arose from positive energy
theorems, and we shall present here some such theorems which follow
from our work elsewhere in the paper. The main point is to give a
complete proof of Herzlich's inequality, \emph{cf.}\/
Theorem~\ref{Tpet3} below. In addition, our remaining results in this
section improve the previous ones
\cite{SchoenYau81,Witten81,ParkerTaubes82,%
ChBlesHouches,Chremark,BM,ChHerzlich,Zhang:hpet,Wang,%
Herzlich97a,Herzlich:mass,Bartnik:mass} in two respects: 1) the
low differentiability of the metric; 2) we do not require $M$ to
have a compact interior. This second aspect of our results is
critical for some applications of the positive mass theorem to
black holes uniqueness theory
\cite{bunting:masood,Ruback,Chstatic,Chstaticelvac}. Some of the
arguments already presented in other sections will be repeated in
the proofs below, whenever useful for the clarity of presentation.

A set $(M,g,K)$ will be called initial data for Einstein equation if
$(M,g)$ is a three\footnote{The results here generalize without any
  difficulties to spin manifolds of higher dimensions, so that the
  restriction $n= 3$ is only made for simplicity of presentation.}
dimensional Riemannian manifold, and $K$ is a symmetric tensor on $M$;
this is a slight abuse of terminology as we are not requiring any
constraints to be satisfied. Only complete $(M,g)$'s will be
considered, with boundary either compact or empty.  Given such a
triple we set
\begin{deqarr}
 & \mu := R - |K|^2 + (\mbox{\rm tr} K)^2\;, & \label{constraint1}\\
& \nu^j := 2 D_i (K^{ij} - \mbox{\rm tr} K g^{ij}) \;. &
\arrlabel{constraints}\label{constraints2}
\end{deqarr}
In a field theoretic framework one should provide some further initial
data for the complete system of equations in the model under
consideration; $\mu$ and $\nu$ correspond then to the energy and
momentum densities of the matter fields while \Eqsone{constraints}
become constraint equations. Now, we are not assuming anything about
matter fields, so $\mu$ and $\nu$ should be thought of as derived from
the data $(M,g,K)$, as in \eq{constraints}.

Throughout this section we shall be working in a space of spinors
which carries a representation $\cgamma $ of the Clifford algebra
associated with the $(n+1)$-dimensional Lorentzian metric
\begin{eqnarray}
&
\gamma := {} -{}\theta^0\otimes \theta^0 + \sum_{i=1}^n \theta^i
\otimes \theta^i
={}-{}\theta^0\otimes \theta^0 + g\;,\label{Lormet} &
\end{eqnarray}
with $\theta^i$ denoting a co-frame dual to a frame $e_i$. We assume
that the $\cgg (e_i)$'s are antisymmetric, and that $\cgamma (e_0)$ is
symmetric. The symbol $D$ will be used to denote the canonical spinor
connection on $(M,g)$ defined by Equation~\eq{A:gpsi} (and denoted by
$\nabla$ there). The connection \bel{4dc} \nabla_j:=D_j + \frac 12
K^i{}_{j}\cgg (e_i)\cgamma (e_0)\ee will be called the
\emph{space-time} spin connection on $M$; $\nabla$ is sometimes
referred to as the \emph{Sen} connection. We note that $\nabla$ is
\emph{not} compatible with the positive definite metric
$\langle\cdot,\cdot\rangle$: \bel{nmc} \partial_j
\langle\phi,\psi\rangle = \langle \nabla_j\phi,\psi\rangle + \langle
\phi,\nabla_j\psi\rangle - \langle\phi,K^i{}_j c(e_i)c(e_0)\psi\rangle
\;. \ee
However, $\nabla$ is compatible with the Lorentz-invariant (and hence
indefinite) inner product $(\phi,\psi) :=
\langle\phi,\cgg(e_0)\psi\rangle$.
The identity associated with $\nabla$, which replaces the identity
\eq{A:Lich} of Schr\"odinger-Lichnerowicz, takes the form
\begin{equation}
  \int_{\Omega}\left(|\nabla\psi|^{2} +
    {\textstyle{\frac{1}{4}}}\left(\mu|\psi|^{2} +
    \nu^i\langle \psi, \cgg (e_ie_0)\psi \rangle
    \right)-|\mfd\psi|^{2}\right) =
  \oint_{\partial \Omega}\langle\psi,\cgg(ne^{A})\nabla_{A}\psi\rangle\;,
\label{A:Lichspt}
\end{equation}
where \bel{opdef} \mfd:= g^{ij}\cgg (e_i) \nabla_j\;. \ee For
sufficiently differentiable $(g,K)$'s, as will be made precise
below, the identity~\eq{A:Lichspt} holds in the following
circumstances:
\begin{enumerate}
\item If $\psi$ is a $H^1_{\loc}$ spinor field, we may take as
  $\Omega$ a domain in $M$ with compact closure and differentiable
  boundary $\partial \Omega$;
\item If $\psi$ is a compactly supported $C^1$ spinor field, then
  \Eq{A:Lichspt} holds with $\Omega=M$, and $\partial \Omega=\partial
  M$; in particular no boundary term is present if $M$ has no
  boundary;
\item Suppose $(M,\nabla)$ admits a weighted Poincar\'e inequality,
  assume $\partial M=\emptyset$, and let $\HH$ be the space defined in
  \eq{D-2b}. We then have the following:
\end{enumerate}
\begin{Lemma}\label{Lcont}
  Suppose that $\partial M=\emptyset$.  The function
\[
C^1_c(M)\ni \psi \to  G(\psi):= \int_{M}\left(\mu|\psi|^{2} +
    \nu^i\langle \psi, \cgg (e_i)\cgamma (e_0)\psi \rangle \right)
\]
extends by continuity to a continuous function on $\HH$, still denoted
by the same symbol.
\end{Lemma}

\proof For $\psi,\chi\in  C^1_c(M)$ the identity \eq{A:Lichspt2} gives
\begin{eqnarray*}
 G(\psi)-G(\chi) & = &  4\int_{M}\left(-\langle \nabla\psi,\nabla
   \psi\rangle+\langle \nabla\chi,\nabla \chi\rangle +\langle \mfd
   \psi,\mfd  \psi\rangle-\langle \mfd \chi,\mfd  \chi\rangle\right)
\\ & = &  4\int_{M}\Big(\langle\nabla\chi- \nabla\psi,\nabla
   \psi+\nabla \chi \rangle
+\langle \mfd \chi+\mfd \psi,\mfd  \psi-\mfd \chi\rangle
\\
&\le& C\|\psi+\chi\|_{\HH}\|\psi-\chi\|_{\HH}\;;
\end{eqnarray*}
in the last step we have used the fact that  the $|\cgg (e_i)|$ are
uniformly bounded.
\qed
Lemma~\ref{Lcont} implies that the left-hand-side of \eq{A:Lichspt} is
continuous on $\HH$ for any measurable $\Omega\subset M$. Since
compactly supported $C^1$ fields are by definition dense in $\HH$, one
easily checks, using continuity, that
\eq{A:Lichspt} holds with $\Omega=M$ and with vanishing
right-hand-side for all $\psi\in\HH$.
(Nonempty compact boundaries $\partial M$ will be considered
shortly.)

\newcommand{\Mext}{M_{\ext}}

We shall say that $M_{\ext}\subset M$ is \emph{an asymptotically
flat end} if $M_{\ext}$ is diffeomorphic to $\R^3 \setminus
B(0,R)$ for some $R$, with \begin{deqarr}&
r^{-1}(g_{ij}-\delta_{ij})\;,\
\partial_kg_{ij}\;,\ K_{ij}\in L^2(M_{\ext})\;,
&\label{afc1.0}\arrlabel{afc} \\&
\partial_kg_{ij}=o(r^{-3/2})\;,\quad K_{ij}=o(r^{-3/2})\;,&\label{afc1.1}
\\& \mu\;,\ \nu\in L^1(M_{\ext})\; ,&\label{afc2}
\end{deqarr} compare Definition~\ref{defF1}.
Those conditions guarantee that the ADM four-momentum of the data
set is finite and well defined, as follows from what is said in
\cite{ChErice} (compare \cite{Bartnik:mass}): \eq{afc1.0} and
 \eq{afc2} guarantee convergence of the mass and momentum integrals, while
\eq{afc1.1} guarantees geometric invariance.

One of the ingredients of Witten-type proofs of  positive energy
theorems is the introduction of appropriate boundary conditions on
the spinor field in the asymptotic regions.
In the asymptotically flat case this is straightforward: one chooses a
$g$-orthonormal triad such that, in the coordinate system of \bref{afc},
\[
dx^k(e_i - \partial_i) \to_{r\to\infty} 0\;, \ \partial_i\left(
  dx^k( e_j) \right)\in L^2(\Mext)\;;
\]
this is easily achieved by a
Gram-Schmidt orthonormalisation of the frame $\{\partial_i\}$.
Then a spin frame on $\Mext$ is introduced, such that the $\cgamma
(e_i)$'s are represented by constant matrices, as in
Section~\ref{secE}. The boundary condition then is that the spinor
field $\psi$, which will be required to solve the generalized Dirac
equation,
\bel{Dieq} \mfd \psi=g^{ij}\cgg (e_i) \nabla_j\psi=0\;,\ee
asymptotes, as $r$ tends to infinity, to a spinor $\psi_\infty$ which,
for $r\ge R$ for some $R$, has constant entries in the spin frame
above. It is convenient to choose $\psi_\infty$ so that $\psi_\infty$
is smooth, and supported in $\Mext$.
The procedure is somewhat more
delicate in the asymptotically hyperboloidal setting; an elegant
geometric framework for such constructions has been provided in
\cite{Bourguignon92,AndDahl}.

Consider the identity \eq{A:Lichspt} with a spinor field
$\psi=\psi_\infty + \chi$, with $\chi$ differentiable and
compactly supported, while $\partial \Omega=S_R$, a coordinate
sphere of radius $R$ in the exterior region, with $R$ large enough
that $\chi$ vanishes there. A classical calculation along the
lines of \cite{Bartnik:mass} shows that the boundary term in
\eq{A:Lichspt} is then proportional to
\[
4\pi p_\alpha\langle\psi_\infty\;,\cgamma
({e^\alpha}e^0)\psi_\infty\rangle  + o(1) \;,
\]
where $p^\alpha$ is the ADM four-momentum of $(\Mext,g)$, with
$o(1)\to 0$ as $R$ tends to infinity. Passing to this limit we
thus have
\begin{equation}
  \int_{M}\left(|\nabla\psi|^{2} +
    {\textstyle{\frac{1}{4}}}\left(\mu|\psi|^{2} +
    \nu^i\langle \psi, \cgg (e_i)\cgamma (e_0)\psi \rangle \right)-|\mfd \psi|^{2}\right) =
 4\pi p_\alpha\langle\psi_\infty\;,\cgamma ({e^\alpha}e^0)\psi_\infty\rangle\;,
\label{A:Lichspt2}
\end{equation}
still for $C^1$ compactly supported $\chi$'s. But the
left-hand-side of \eq{A:Lichspt2} is continuous on $\HH$,
which is shown by a calculation similar to that in Lemma~\ref{Lcont}:
 Let $F(\chi)$ denote the left-hand-side of \Eq{A:Lichspt2} with
$\psi=\psi_\infty+\chi$ there, let $\chi_i\in \HH$ converge in $\HH$
to $\chi\in \HH$, so we have
\begin{eqnarray*}
F(\chi)-F(\chi_i) & = & \|\chi\|^2_\HH-\|\chi_i\|^2_\HH \\
& & +2 \int_{M}
\langle\nabla^k\psi_\infty\;,\nabla_k(\chi-\chi_i)\rangle \\
&& -  2 \int_{M}
\langle{\mfd}\psi_\infty\;,{\mfd}(\chi-\chi_i)\rangle \\
&&  + \frac{1}{2}\int_{M}
\langle \psi_\infty\;, \left(\mu + \nu^k\cgg (e_k)\cgamma
  (e_0)\right)(\chi-\chi_i)\rangle\;.
\end{eqnarray*}
Because $\nabla\psi_\infty\in L^2(M)$, the first three terms above
converge to zero as $i\to\infty$.
The convergence of the final term can be
justified by applying the Cauchy-Schwarz inequality
$|\langle u, Q v\rangle|\le \sqrt{\langle u, Q u\rangle}\sqrt{\langle
  v, Q v\rangle} $ whenever $Q$ is positive:
\begin{eqnarray*}
\lefteqn{\left|\int_{M}
\langle \psi_\infty\;, \left(\mu + \nu^j\cgg (e_j)\cgamma
  (e_0)\right)(\chi-\chi_i)\rangle\right|} &&
\\&&\le
\left(\int_{M}\langle \psi_\infty\;, \left(\mu + \nu^j\cgg
    (e_j)\cgamma (e_0)\right)\psi_\infty\rangle\right)^{1/2}
\\&&\phantom{\le}\times
\left(\int_{M}\langle (\chi-\chi_i), \left(\mu + \nu^j\cgg
    (e_j)\cgamma (e_0)\right)(\chi-\chi_i)\rangle\right)^{1/2}
\\
&& \le
C(\|\mu\|_{L^1}+\|\nu\|_{L^1})
  \Big(\|\chi+\chi_i\|_\HH\|\chi-\chi_i\|_\HH\Big)^{1/2}\;;
\end{eqnarray*}
in the last step of the calculation of Lemma~\ref{Lcont} has been used.
Now, $F(\chi_i)=F(0)$, and density implies that
\eq{A:Lichspt2} remains true for any $\psi $ of the form
$\psi_\infty+\chi$, with $\chi\in \hb$.

We are ready now to prove the following version of the positive
energy theorem, the regularity conditions of which  have been
chosen as a compromise between those needed for solvability of the
Dirac equation \bref{Dieq} ({\em cf.}\/ Remark~\ppref{Rreg}),
those needed for a well defined notion of ADM mass, and those
needed for a Banach manifold structure for the set of solutions of
the general relativistic vacuum constraint equations.

\begin{Theorem}\label{pet1}
Let $(M,g,K)$ be initial data for the Einstein equations with $g\in
W^{2,2}_{\loc}$, $K\in W^{1,2}_{\loc}$, with $(M,g)$ complete
(without boundary). Suppose that $M$ contains an asymptotically
flat end and let $p_\alpha=(m,\vec p)$ be the associated ADM
four-momentum.\footnote{There is a signature-dependent ambiguity
in the relationship between $p^0$, $p_0$ and the mass $m$: in the
space-time signature $(-,+,+,+)$ used in this paper this sign is
determined by the fact that $p_0$, obtained by  Hamiltonian
methods, is usually positive in Lagrangean theories on Minkowski
space-time such as the Maxwell theory, while the mass $m$ is a
quantity which is expected to be positive.} If
\begin{equation}
 \mu \ge |\nu|_g\;,
\label{DEC}
\end{equation}
then
\bel{posmasineq} m \ge |\vec p|_{\delta}\;, \ee
with equality if and only if $m$ vanishes. Further, in that last
case there exists a non-trivial covariantly constant (with respect
to the space-time spin connection) spinor field on M.
\end{Theorem}

\begin{remark} Under the supplementary assumption of smoothness of $g$
and $K$, it has been shown in \cite{ChBeig1} that the existence of
a covariantly constant spinor implies that the initial data can be
isometrically embedded into Minkowski space-time, \emph{cf.} also
\cite{Yip}. We expect this result to remain true under the current
hypotheses, but we have not attempted to prove this.
\end{remark}

\Proof Suppose that for all $\psi_\infty$ we can establish
existence of $\chi=\chi[\psi_\infty]\in\HH$ such that
$\psi_\infty+\chi$ satisfies the Dirac equation \eq{Dieq}.
\Eq{A:Lichspt2} would then show that the quadratic form
\[
\psi_\infty \to 4\pi p_\alpha\langle\psi_\infty\;,\cgamma
 ({e^\alpha}e^0)\psi_\infty\rangle
\]
 is non-negative, and the
Theorem  follows by a standard calculation.
The existence of $\chi$ will be a
consequence of Theorem~\ref{ThmD3},
provided that the relevant
hypotheses are met.
We have\newcommand{\trg}{\textrm{tr}_g}
$$
\mfd = c(e^i)D_i-\frac 12 \trg K \cgamma (e_0)\quad
\Longrightarrow\quad
 \mfd ^\dagger= \mfd \;.
$$
\Eq{A:Lichspt2} with $\psi_\infty=0$ shows that the
Schr\"odinger-Lichnerowicz estimate of Definition~\ref{defD2},
with $Y=\emptyset$ and  $\rho=0$, holds both for $\cL:=\mfd $
and its formal adjoint $\mfd ^\dagger=\mfd $. Next, we note that
the symmetric part $\Gamma^S$ of the connection~\eq{4dc} is
$$\Gamma^S=\frac 12 K^i{}_{j}\cgg (e_i)\cgamma (e_0)\otimes dx^j$$
which does not vanish for non-zero $K$'s, but satisfies nevertheless
the fall-off condition \eq{F-4a} by \eq{afc1.1}. It follows from
Theorem~\ref{ThmF3} that the weighted Poincar\'e inequality holds.
The regularity conditions on the metric imply that the requirements of
Hypothesis \ref{Hh} with $\mathring{g}=g$ are met: for trivial
bundles, or for smooth initial data, this is a straightforward
calculation, compare Remark~\ref{Rreg}; for non-trivial bundles
Proposition~\ref{pspinnew}, Appendix~\ref{AppAGC}, has to be
invoked.
 The map
$\beta$ of \Eq{D-5b} is zero, as is the curvature term $\rho$ in the
Schr\"odinger-Lichnerowicz inequality \bref{D:SL}, by the energy
condition $\mu\ge|\nu|_g$.   From what has been said it follows that
spinor fields in $\HH$ which are also in the kernel of $\mfd $ are
covariantly constant; they then have constant length, and are not in
$L^2$ if they are non-zero.  Theorem~\ref{ThmD3} now shows that for
any $\psi_\infty$ there exists a solution $\chi\in\HH$ of the equation
$$ \mfd \chi = - \mfd \psi_\infty\;,$$ and the existence of the
desired $\psi$ follows. \QED

\newcommand{\pM}{\partial M}%
Let us now turn our attention to manifolds with boundary. We shall
say that a boundary $\pM$ is \emph{future-trapped} if
\bel{trapped}\theta_+:= H+\sum_{A=2,3}K(e_A,e_A)\le 0.\ee Here
$H$ is the mean curvature of $\pM$ with respect to an
inner-pointing normal, while the $e_A$'s form an ON basis for
$T\pM$. A future-trapped boundary in the sense above is
future-trapped in the usual sense~\cite{HE} for a surface in
space-time. The following result generalises one by
Herzlich\footnote{The proof in~\cite{Herzlich:mass} is the
rigorous version of an argument proposed in~\cite{GHHP83}; it
also extends that argument, as  in~\cite{GHHP83} only
\emph{marginally trapped} boundaries
  are considered.}~\cite{Herzlich:mass}:

\begin{Theorem}
\label{Tpet2}
Under the remaining hypotheses of Theorem~\ref{pet1}, suppose instead
that $M$ has a differentiable, compact, future-trapped boundary
$\partial M$. Then the conclusions of Theorem~\ref{pet1} hold.
\end{Theorem}

\begin{NRemark} One expects that the equality case cannot occur in
  \bref{posmasineq}, and a possible argument could proceed as follows:
  First, the existence of a covariantly constant spinor implies
  existence of a non-spacelike, covariantly constant, Killing vector
  field in the associated space-time. Further, if the metric is $C^2$
  and $K$ is $C^1$, then the space-time metric fulfills the Einstein
  equations with a null fluid as a source~\cite[Appendix~B]{ChBeig1};
  this conclusion is expected to hold under the weaker
  differentiability conditions considered here. By reduction of the
  field equations, this should imply smoothness of the metric.
  (Alternatively, one could assume at the outset that $g$ is $C^3$ and
  $K$ is $C^2$, in which case the argument presented in the current
  remark settles the issue).  Topological censorship results
  \cite{ChWald} applied to the Killing development \cite{ChBeig1} of
  the initial data show that the boundary is then the union of a
  finite number of spheres. Arguing as in the proof of Theorem~4.6
  of~\cite{ChHerzlich}, the restriction to the boundary of the
  covariantly constant spinor would be harmonic, which is impossible
  by the Hijazi-B\"ar inequality \eq{Hbest}~\cite{Hijazi91,Bar92}.
\end{NRemark}

\begin{remark} \emph{Past-trapped} boundaries are defined by changing the
sign of $K$ in \eq{trapped}; as the remaining hypotheses of
Theorem~\ref{Tpet2} are invariant under this change of sign, an
identical result holds for compact past-trapped $\partial M$'s.
\end{remark}

 \proof The proof follows closely that of
Theorem~\ref{pet1}, the main difference being the need to impose
suitable boundary conditions. Indeed, when $\pM$
is non-empty \Eq{A:Lichspt2} becomes
\begin{eqnarray}
\lefteqn{
  \int_{M}\left(|\nabla\psi|^{2} +
    {\textstyle{\frac{1}{4}}}\left(\mu|\psi|^{2} +
    \nu^i\langle \psi, \cgg (e_i)\cgamma (e_0)\psi \rangle \right)
    -|\mfd \psi|^{2}\right)
    }&&
 \nonumber \\ &&=
 4\pi p_\alpha\langle\psi_\infty\;,\cgamma ({e^\alpha}e^0)\psi_\infty\rangle+
\oint_{\pM}\langle\psi,\cgg (n)\cgg (e^{A})\nabla_{A}\psi\rangle
\nonumber
\\ &&=
4\pi p_\alpha\langle\psi_\infty\;,\cgamma
  ({e^\alpha}e^0)\psi_\infty\rangle
 + \nonumber \\ &&   \oint_{\pM} \langle\psi,\mfdb\psi +
\half\Big( H
-\sum_A(K_{AA}\cgg (n)-K_{A1}\cgg (e_A)\Big)\cgamma (e_0)\psi\rangle\;,
\label{A:Lichspt4}
\end{eqnarray}
for, say, continuously differentiable $\psi$'s of the form
$\psi=\psi_\infty+\chi$, with $\chi$ compactly supported, and
$\psi_\infty$ as in the proof of Theorem~\ref{pet1}. Further, $\mfdb $
is the Dirac boundary operator defined by \Eq{A:bdirac}, and
\Eq{A:Lichb} has been used. Finally, $e_i$ is an ON frame on $\pM$
with $n\equiv e_1$ normal to $\partial M$. Following \cite{GHHP} we
impose the boundary condition \bel{GHHPbc} \cK_- := \half
(1-\epsilon)\psi = 0 \quad
\textrm{on } \pM\;, \ee where $\epsilon:=-\cgg (n)\cgamma (e_0)$.
We then have
\begin{eqnarray*}
\langle \psi,\cgg (e_A)\cgamma (e_0)\psi\rangle & = & -\langle
\psi,\cgg (e_A)\cgamma (e_0)\cgg (n)\cgamma (e_0)\psi\rangle \\ & =
& -\langle
\psi,\cgamma (e_0)\cgg (n)\cgg (e_A)\cgamma (e_0)\psi\rangle \\ & =
& -\langle
\cgg (n)^t\cgamma (e_0)^t\psi,\cgg (e_A)\cgamma (e_0)\psi\rangle \\
& = & \langle
\cgg (n)\cgamma (e_0)\psi,\cgg (e_A)\cgamma (e_0)\psi\rangle=-\langle
\psi,\cgg (e_A)\cgamma (e_0)\psi\rangle\;,
\end{eqnarray*}
which shows that the last term in the last line of \Eq{A:Lichspt4}
vanishes. Then
$$\mfdb  \epsilon = -\epsilon \mfdb $$ and $\epsilon^t=\epsilon$, so
\begin{eqnarray*}
\langle \psi,\mfdb \psi\rangle & = & \langle \psi,\mfdb \epsilon\psi\rangle \\
& = & -\langle \psi,\epsilon\mfdb \psi\rangle \\ & = & -\langle
\epsilon\psi,\mfdb \psi\rangle =-\langle \psi,\mfdb \psi\rangle\;,
\end{eqnarray*}
which shows that the first term in the last line of
\Eq{A:Lichspt4} vanishes. Next,
$$
\langle\psi,\Big( H
-\sum_AK_{AA}\cgg (n)\cgamma (e_0)\Big)\psi\rangle =
\langle\psi,\Big( H +\sum_AK_{AA}\Big)\psi\rangle =
\theta_+\langle\psi,\psi\rangle\;,
$$ which shows that the sum of the second and third term in the
last line of \Eq{A:Lichspt4} gives a non-positive contribution
when $\pM$ is trapped. When $\mu\ge |\nu|_g$, $\psi_\infty=0$, and
\eq{GHHPbc} holds, from \Eq{A:Lichspt4} we obtain
\begin{eqnarray}\nonumber 
& \displaystyle\int_{M}\left(|\nabla\psi|^{2} +
    {\textstyle{\frac{1}{4}}}\left(\mu|\psi|^{2} +
    \nu^i\langle \psi, \cgg (e_i)\cgamma (e_0)\psi \rangle \right)\right)-
    \half\oint_{\pM}\theta_+\langle\psi,\psi\rangle = \int_M|\mfd \psi|^{2}
 \;,
 &     \label{A:Lichspt5}
\end{eqnarray}
so conditions \bref{DEC}, \bref{trapped} give
\begin{equation}
  \int_{M}|\nabla\psi|^{2} \le \int_{M}|\mfd \psi|^{2}\;,
\label{A:Lichspt7}
\end{equation}
for all $\psi\in C^1_c(M)$ which satisfy \bref{GHHPbc}.
Define
\begin{equation}
  \label{HK}
\HH_{\cK_-}  := \{ \psi\in \HH, \ \cK_-\psi= 0 \textrm{  on } \pM\}\;,
\end{equation}
where $\HH$ is defined in \bref{D-2b}. If we  let
 \bel{gpri} G'(\psi):=\displaystyle\int_{M}\left(
  {\textstyle{\frac{1}{4}}}\left(\mu|\psi|^{2} + \nu^i\langle \psi,
    \cgg (e_i)\cgamma (e_0)\psi \rangle \right)\right)-
\half\oint_{\pM}\theta_+\langle\psi,\psi\rangle\;,\ee then the
calculation of the proof of Lemma~\ref{Lcont} shows that
$G'(\psi)$ can be extended by continuity to a continuous function
on
 $\HH_{\cK_-}$.
The boundary integral in \eq{gpri} is
 continuous, so the volume integral is also continuous on
 $\HH_{\cK_-}$. which implies that \Eqsone{A:Lichspt5}  holds for all
$\psi\in\HH_{\cK_-}$. Lemma~\ref{LemD1a} establishes the
Schr\"odinger-Lichnerowicz estimate  \bref{D:SL} with $\rho=0$ for
$(\mfd ,\cK_-)$.

As explained in Section~\ref{secE}, the boundary value problem
determined by \Eq{GHHPbc} belongs to the family of problems
considered in Theorem~\ref{ThmD3}. If both $M$ and $\pM$ are
simultaneously parallelizable, and if the metric $g$ is a product
near $\partial M$, then the regularity conditions of
Theorem~\ref{ThmD3} are met by hypothesis; the general case is
handled by Propositions~\ref{pspinnew} and \ref{PAGC.1},
Appendix~\ref{AppAGC}. Repeating now the arguments of the
proof of Theorem~\ref{pet1} gives the non-negativity of $p_0$.
\qed

It is expected that the positivity statement of
Theorem~\ref{Tpet2} can be strengthened to the so-called Penrose
inequality when trapped boundaries occur. This question remains
wide open, except in the special\footnote{Similarly to
Theorem~\ref{Tpet3}, for the results in
\cite{HI1,Bray:preparation2} it actually suffices that $R\ge 0
\Leftrightarrow \mu\ge -|K|^2_g+(\textrm{tr}_gK)^2$, and that
$\textrm{tr}_h K$ vanishes on $\partial M$, where $\textrm{tr}_h$
is the trace of the restriction of $K$ to $\pM$.} case $K_{ij}=0$
\cite{HI1,Bray:preparation2}. An interesting related inequality
has been, essentially, proved by Herzlich~\cite{Herzlich97a};
however, the arguments of that last reference do not include a
sufficient justification of existence of the required spinor
field, except in the rather special case of a smooth metric which
is a product near the boundary, as analyzed by Bunke~\cite{Bunke}.
Here we fill this gap and establish the following:

\begin{Theorem}
\label{Tpet3}
  Let $(M,g)$ be a complete Riemannian manifold
with $g\in W^{2,2}_{\loc}$, and suppose that $M$ has a  boundary
$\partial M$ diffeomorphic to $S^2$, with non-positive inwards
pointing mean curvature. Suppose that the curvature scalar $R(g)$ of
the metric $g$ is non-negative, and that $M$ contains an
asymptotically flat end with mass $m$.
If $\sigma$ is the
dimensionless quantity defined as $$\sigma:=
\sqrt{\frac{\textrm{Area}(\pM)}{\pi}} \inf_{f\in C^\infty_c(\bM),
f\not \equiv 0}\frac{\|df\|^2_{L^2(M)}}{\|f\|^2_{L^2(\pM)}}\;,$$
then $$m\ge \frac \sigma{1+\sigma}\sqrt{\frac{\textrm{Area}(\pM)}{4\pi}} \;.$$
Moreover, if the metric is smooth, then equality is achieved if and
only if $(M,g)$ can be isometrically embedded in the
Schwarzschild space-time with mass
$\sqrt{{\textrm{Area}(\pM)}/{16\pi}}$.
\end{Theorem}

\begin{remark} If $M$ is the union of a compact set with a finite number
of asymptotically flat ends, then $\sigma>0$.
\end{remark}

\medskip

\proof The details of the argument follow closely those of the
proof of Theorem~\ref{Tpet2}, the pointwise boundary conditions
\eq{GHHPbc} being replaced by the spectral boundary conditions
\eq{A:bMle0} with $\cal K$ given by \eq{APS}; compare the
discussion of Section~\ref{secA}, as well as that in the paragraph
following \Eq{E-7}. The main elements missing in the arguments
of~\cite{Herzlich97a} are provided by the boundary regularity
results of Section~\ref{secB}; those are the key to the proof of
Theorem~\ref{ThmD3}. The reader is referred to
\cite{Herzlich:mass} and \cite[p.~679]{Chstatic} for the analysis
of the equality case, \emph{cf.\/} also \cite{Beigconformal}.
\qed

Following \cite{GHHP83,GibbonsHull}, let us pass now to
inequalities with an electric charge contribution. A set
$(M,g,K,E,B)$ will be called \emph{initial data for the
Einstein-Maxwell equations} if $(M,g)$ is a three dimensional
Riemannian manifold, $K$ is a symmetric tensor on $M$, while $E$
and $B$ are vector fields on $M$; as before, this is a serious
abuse of terminology, as we are not requiring any constraint
equations to be satisfied. Given such a triple we set
\begin{eqnarray}
\nonumber & \divE :=D_iE^i\;, & \\ \nonumber& \divB :=D_iB^i\;, & \\
\nonumber&\mu  := R - |K|_g^2 + (\mbox{\rm tr} K)^2 - 2|E|_g^2 -2
|B|_g^2\;, & \\ & \nu_i := 2D_j (K^{j}{}_i - \mbox{\rm tr} K
\delta^{j}_i) +4\epsilon_{ijk}E^j B^k\;. \label{EMconstraints}
\end{eqnarray}
Here $D$ is the Levi-Civita connection associated with the metric
$g$. In a general relativistic context, $\divE $ is the electric
charge density, $\divB $ is the magnetic charge density (usually
zero, whether electro-vacuum or not), $\mu $ is the energy
density remaining after subtracting the electro-magnetic
contribution, and $\nu_i$ is the left-over matter current; $\mu $,
$\divE $, $\divB $ and $J_i$ vanish when the Einstein-Maxwell
constraint equations hold.

\newcommand{\clg}[1]{\cgg (#1)}%

Let a new connection $\nabla$ be defined as\
\begin{equation}
  \label{4connect}
  \nabla_i:=D_i + {1\over 2 } K_{ij}\cgg ({e^j})\cgamma (e_0) - \frac 12
   \clg {E}\clg {e_i}\cgamma (e_0) - \frac 14
  \epsilon_{jk\ell}B^j\clg {e^k}\clg {e^\ell}\clg {e_i}\;;
\end{equation}
The connection $\nabla $ will be called the \emph{space-time
  Einstein-Maxwell spin connection}  on $M$.
$\nabla$ is again \emph{not} metric compatible, with symmetric
part $\Gamma^S$ given by \bel{EMsymp}\Gamma^S=\Big({1\over 2 }
K_{ij}\cgg ({e^j})\cgamma (e_0) - {E_i}\cgamma (e_0) - \frac 14
  \epsilon_{jk\ell}B^j\clg {e^k}\clg {e^\ell}\clg {e_i}\Big)\otimes
  \theta^i\;,\ee
  where, as before, $\theta^i$ is the co-frame dual to $e_i$.
  In this context the
  asymptotic flatness conditions have to be complemented by
  conditions on $E$ and $B$: we shall require
  \bel{Maxbc}
  E,B\in L^2(\Mext)\;,\quad \divE ,\divB \in L^1(\Mext)\;,\quad E=
  o(r^{-1})\;, \quad B= o(r^{-1})\;.
  \ee
  Following an argument proposed by Gibbons and Hull~\cite{GibbonsHull} we
  have:
\begin{Theorem}\label{pet3}
  Let $(M,g,K,E,B)$ be initial data for Einstein-Maxwell equations with $g\in
  W^{2,2}_{\loc}$, $K,E,B\in W^{1,2}_{\loc}$, with $(M,g)$ complete
  (without boundary). Suppose that $M$ contains an asymptotically
  flat end $\Mext$ with $E$ and $B$ satisfying the fall-off conditions
  \eq{Maxbc} there,
  and let $ p_\alpha=(m,\vec p)$ be the associated ADM
  four-momentum. Let $Q$ and $P$ be the total electric and magnetic
  charge of $M_\ext$,
$$ Q=\lim_{R\to\infty}\frac 1 {4\pi} \oint_{r=R} E^i dS_i\;, \qquad
P=\lim_{R\to\infty}\frac 1 {4\pi}\oint_{r=R} B^i dS_i\;.$$ If
\begin{equation}\label{hypineq} \mu  \ge \sqrt{|\nu|^2_g + |\divE
    |^2+|\divB |^2}\;,
 \end{equation}
then \bel{Maxwellineq}
 m \ge \sqrt{|\vec p|_\delta^2 + Q^2+P^2}\;,
\ee where $|\vec p|_\delta \equiv \sqrt{\sum(p^i)^2}$, with
equality if and only if there exists a spinor field on M which is
covariantly constant with respect to the  Einstein-Maxwell
space-time spin connection \eq{4connect}.
\end{Theorem}

\begin{remark} Under the hypothesis of smoothness of the metric,
Tod~\cite{Tod} has found the local form of the metrics which admit
covariantly constant spinors as above; however, no classification
of globally regular such space-times is known. It is expected that
the only singularity-free solutions here have vanishing Maxwell
field, or belong to the \emph{standard} Majumdar-Papapetrou family
(\emph{cf., e.g.,}\/~\cite{ChNad}). It would be of interest to
fill this gap.\end{remark}

\begin{remark}
 Charged matter
  might violate \eq{hypineq}; however, there might exist a constant
  $\alpha\in(0,1)$ such that
\begin{equation}\label{hypineq2}
\mu  \ge \sqrt{|\nu|^2_g + \alpha^2(|\divE |^2+|\divB |^2)}\;.
\end{equation}
Replacing in \eq{4connect} the fields $E$ and $B$ by $\alpha E$ and
$\alpha B$, an essentially identical argument leads to
\bel{Maxwellineq2} m \ge \sqrt{|\vec p|_\delta^2 +
  \alpha^2(Q^2+P^2)}\;.  \ee\end{remark}

 \proof
For the connection~\eq{4connect} the identity~\eq{A:Lichspt2}
becomes \cite{GibbonsHull}
\begin{eqnarray}
 \lefteqn{ \int_{M}\left.\Big( |\nabla\psi|^{2}  -|\mfd \psi|^{2}\right.}&&
  \nonumber \\
  &&+
    {\textstyle{\frac{1}{4}}}\langle \psi,\left(\mu +
    \nu^i \cgg (e_i)\cgamma (e_0) - \divE \cgamma (e_0) -
    \divB \clg {n}\clg {e_2}\clg {e_3} \right)\psi \rangle
    \Big)
  \nonumber \\ && =4\pi
 \langle\psi_\infty\;, \left[p_\alpha\cgamma ({e^\alpha})\cgamma (e^0)
    +   Q \cgamma (e_0)-
 P\clg {n}\clg {e_2}\clg {e_3} \right] \psi_\infty\rangle\;,
\label{A:Lichspt6}
\end{eqnarray}
again for  $\psi $ of the form $\psi_\infty+\chi$, with  $C^1$
compactly supported $\chi$'s. The Dirac operator
$$\mfd :=\clg{e^i}\nabla_i=\clg{e^i}D_i - \frac 12 \Big(\trg K -\clg E
\Big)\cgamma (e_0)
-\frac 14 \epsilon_{jk\ell}B^j\clg {e^k} \clg{e^\ell}$$ is not
formally self-adjoint, we have instead
$$\mfd ^\dagger=\clg{e^i}D_i - \frac 12 \Big(\trg K -\clg E \Big)\cgamma (e_0)
+\frac 14 \epsilon_{jk\ell}B^j\clg {e^k} \clg{e^\ell}\;.$$ This
shows that the adjoint of $\mfd $ coincides with $\mfd $ modulo
the replacement
$$B\to-B\;.$$ The
arguments follow now the previous ones, basing on the identity
\eq{A:Lichspt6}. We simply note that \eq{hypineq}
implies non-negativity of the quadratic form appearing in the second line
of \eq{A:Lichspt6}. Similarly, positivity of the quadratic form
defined by the third line of \Eq{A:Lichspt6} 
implies \eq{Maxwellineq}. Some comments are in order here, related
to the fact that $\cL$ is not formally self-adjoint when the
magnetic field does not vanish.  Since we are not assuming
interior compactness of $M$, $M$ could have other asymptotic
regions in which $B$ could grow in an uncontrollable way, so that
$\cL^\dag$ will not map $\hb$ into $L^2$. Now, $\cL^\dag$ differs
from $\cL$ by a change of the sign of $B$, which implies that
$\cL^\dag$ also satisfies a Schr\"odinger-Lichnerowicz identity
with a connection in which $B$ is replaced by $-B$. The arguments
already given show that the weak equation $\cL^\dag \phi=0$ has no
$L^2$ solutions, and Corollary~\ref{ThmD3} provides the desired
isomorphism property of $\cL$.  \qed

In the presence of boundaries we have:

\begin{Theorem}\label{pet4} Under the remaining hypotheses of
Theorem~\ref{pet3}, suppose instead that $M$ has a compact
future-trapped boundary $\partial M$. Then the conclusions of
Theorem~\ref{pet3} hold.
\end{Theorem}
\proof This is a repetition of the argument of the proof of
Theorem~\ref{Tpet2}; one imposes again the boundary condition
\eq{GHHPbc}, and we only need to check  that \eq{A:Lichspt7} still
holds. This is indeed the case, which is established as follows:
the electromagnetic field leads to a supplementary contribution
$$ \oint_{\pM} \langle \psi, \left[E^i\cgamma (e_0) - B^i
\cgg (e^1)\cgg (e^2)\cgg (e^3)\right]\psi\rangle n_i
$$
to the boundary integral \eq{A:Lichspt4}. When \eq{GHHPbc} holds
we have
\begin{eqnarray*}
\langle \psi,\cgamma (e_0) \psi \rangle & = & -\langle
\psi,\cgamma (e_0)\cgg (n) \cgamma (e_0)\psi \rangle
\\& = & \langle
\psi,\cgg (n) \psi \rangle
\\& = & -\langle \cgg (n)
\psi, \psi \rangle
\\ & = & -\langle \cgg (n)\cgamma (e_0)\cgamma (e_0)
\psi, \psi \rangle
\\ & = &\langle \cgamma (e_0)\cgg (n)\cgamma (e_0)
\psi, \psi \rangle
\\ & = &-\langle \cgamma (e_0)
\psi, \psi \rangle = - \langle \psi,\cgamma (e_0) \psi \rangle\;,
\end{eqnarray*}
hence $$\langle \psi,\cgamma (e_0) \psi \rangle=0.$$ Similar
manipulations show that \Eq{GHHPbc} implies
$$\langle \psi,\cgg (e^1)\cgg (e^2)\cgg (e^3) \psi \rangle=0\;,$$
and the result follows. \qed

 We finish this section by
noting that positive energy results follow by identical arguments
for asymptotically hyperboloidal manifolds
\cite{ChHerzlich,AndDahl,Zhang:hpet,Wang,GHHP,ReulaTod,DahlBanach,CJL};
here Theorem~\ref{ThmF3.1} should be used instead of
Theorem~\ref{ThmF3}. The definition of mass in that case is
considerably more delicate, we refer the reader to
\cite{ChHerzlich,CJL} for details.

\appendix

\section{Fields on manifolds of $W^{k+1,p}$ differentiability class}
\label{AppAGC}

Consider a smooth manifold $M$; on such a manifold one can define in a
geometrically invariant way tensor
fields which are of $C^\infty$ differentiability class, or of $C^k$
class, or of $W^{k,p}_\loc$ class. For example, one says that a tensor
field is of $W^{k,p}_\loc$ class if there exists a covering of $M$ by
coordinate patches such that the coordinate components of the tensor
in question are in $W^{k,p}_\loc$ in each of the coordinate patches.
Since the transition functions when going from one coordinate system
to another are smooth, this property will be true in any coordinate
system.

Let, now, $(M,g)$ be a smooth manifold with a \Riemn metric $g$ which
is of $W^{k,p}_\loc$ differentiability class. For various arguments it
is convenient to use local coordinate systems which are adapted to the
metric, such as geodesic coordinates, or harmonic coordinates.  In
this case the transition functions to the adapted coordinate system
will not belong to the original smooth atlas on $M$ in general. At
this point there are two strategies possible: either to enlarge the
atlas on $M$ to contain those new coordinate systems, or to ignore
this issue and try to analyze the problems that arise on an \emph{ad
  hoc} basis. For nearly all of this paper the \emph{ad hoc} approach,
working entirely within a $C^\infty$ structure on $M$, is quite
adequate.  However, the proofs of Theorems \ref{Tpet3} and \ref{pet4}
require the existence of approximately Gaussian coordinates near a
boundary or near a hypersurface of $M$ (\emph{cf.\/} Lemma
\ref{PAGC.1} below).  Direct construction of such coordinates with
respect to the $W^{k,p}_\loc$ metric $g$ produces a coordinate change
which is not $C^\infty$, which forces us to analyse the problems
involved when constructing systematically manifolds of
$W^{k+1,p}_\loc$ differentiability class.  For this reason we will
present such a construction here.  For technical reasons we shall
always assume that
\begin{equation}
  \label{kpn}
p\in [1,\infty]\ , \qquad k\in \bN\ , \qquad  kp > n\ ;
\end{equation}
these restrictions are more than sufficient for our purposes.
Generalising the condition $k\in\bN$ to $k\in \bR^+$ would require an
analogue of Lemma \ref{lcomp} for non-integral $k,\ell$, which seems
not to be available.  Condition \eq{kpn} and the Sobolev embedding
$W^{k+1,p}_\loc(\Omega)\subset C^1(\Omega)$ (for appropriately regular
open\ domains $\Omega\subset\bR^n$) mean that we will consider only
manifolds which are at least of $C^1$ differentiability class.

Consider, thus, a connected paracompact Hausdorff manifold $M$ of
$C^1$ differentiability class. We shall say that $M$ is of
$W^{k+1,p}_\loc$ differentiability class if $M$ has an atlas for which
all the transition functions are of $W^{k+1,p}_\loc$ differentiability
class. Unless indicated otherwise, the Lebesgue measure in local
coordinates is used.



A tensor field with components which are $C^\infty$  with respect to
some coordinate chart (belonging to a $C^\infty$ sub-atlas of the
$W^{k+1,p}_\loc$ atlas), will not generally have smooth components in
all $W^{k+1,p}_\loc$ charts.
In this situation a
$C^\infty$ tensor field ``comes equipped'' with a preferred atlas of
coordinate charts in which it has smooth coordinate components. This
is \emph{a priori\/} the case for tensor fields of \emph{any}
differentiability class on $\wkpl$ manifolds, and it is of interest to
single out those classes of tensor fields, the coordinate components of
which will be of a prescribed differentiability class in \emph{every}
coordinate system of the $\wkpl$ atlas on $M$. Differentiability
classes of this type will be referred to as \emph{invariantly defined}.
Our next result describes some such classes of tensor fields.  It is
convenient to introduce the following notation: let $x,y\in \bR$, we
shall write $x\gst y $ if the following holds:
\begin{equation}
  \label{eq:lst}
  x\gst y \Longleftrightarrow \cases{x\ge y \ ,&  if $y>0$\ , \cr
                                     x> y  \ , & if $y\le0$\  .}
\end{equation}
(We note that for $x\ge0$ the only value of $x$ at which ``$\gst$''
does not coincide with ``$\ge$'' is $x=0$.)  In this notation the
the Sobolev embedding theorem can be stated as:
\begin{equation}
  \label{eq:sobemb}
  W^{s,t}_\loc \subset W^{u,v}_\loc \quad  \Longleftrightarrow\quad u\le s
\ \mbox{  and  }\ \frac{1}{v} \gst \frac{1}{t}-\frac{s-u}{n}\ \;.
\end{equation}
\begin{Proposition}
  \label{pinvnew} Let $(M,g)$ be a $\wkpl$  manifold, $kp>n$, $p\in
  [1,\infty]$.
  \begin{enumerate}
\item Let $(\ell,q)$ be such that the Sobolev embedding
$$
\wkpl\subset W^{\ell,q}_\loc
$$
holds (which is equivalent to the condition
\begin{equation}
  \label{eq:p-1}
  \frac{1}{q}\gst \frac{1}{p}+ \frac{\ell - k - 1}{n}\ \;,
\end{equation}
with $\gst$ defined in \eq{eq:lst}). Then the space of
$W^{\ell,q}_\loc$ scalar fields on $M$ is invariantly defined.
  \item Let $(\ell,q)$ be such that the Sobolev embedding
    $$
    \wkpgl\subset W^{\ell,q}_\loc
    $$
    holds (which is equivalent to the condition
\begin{equation}
  \label{eq:p0}
  \frac{1}{q}\gst \frac{1}{p}+ \frac{\ell - k}{n}\ .)
\end{equation}
Then the space of $W^{\ell,q}_\loc$ tensor fields on $M$
    is invariantly defined.
  \end{enumerate}
\end{Proposition}

\proof Point 1 is a straightforward consequence of the following
Lemma:

\begin{Lemma}
  \label{lcomp}
  Let $\Omega, \cU \subset \bR^n$ and let $\psi: \Omega\to \cU $ be a
  $C^1$ diffeomorphism such that $\psi\in \wkpl(\Omega;\bR^n)$, $kp>n$.
  If $(\ell,q)$ is such that the Sobolev embedding
  $\wkpl\subset\wlqgl$ holds, \emph{cf.\/} Equations \eq{eq:p-1} and
  \eq{eq:lst}, then for all $F\in \wlqgl(\cU)$ we have
$$F\circ \psi \in \wlqgl(\Omega)\ .$$
\end{Lemma}

\begin{remark} In \cite{Bourdaud,Sickel} some partial results can be found
concerning sharpness of this result.
\end{remark}

\medskip

\proof
We have, for $0\le |\alpha|\le \ell(\le k+1)$,
$$
\partial^ \alpha (F\circ \psi) = \sum C(\alpha_1,\ldots,\alpha_m)
\partial^ {\alpha_1} \psi \cdots \partial^ {\alpha_m} \psi F^{(m)}\circ
\psi\ ,
$$
where the sum is taken over sets $(\alpha_1,\ldots,\alpha_m)$
satisfying $\alpha_1+\cdots+\alpha_m=\alpha$, with $|\alpha_i|\ge 1$.
For any compact $K\subset \Omega$ it follows that
\begin{eqnarray}
  \label{eq:c1}
&  \|\partial^ \alpha (F\circ \psi)\|_{L^q(K)}\le C  \sum
\|\partial^ {\alpha_1} \psi\|_{L^{s(\alpha_1)}(K)} \cdots \|\partial^
{\alpha_m}
\psi\|_{L^{s(\alpha_m)}(K)} \|F^{(m)}\|_{L^t(\psi(K))}\ ,
& \\ \label{c2}
& \mbox{$\displaystyle \frac{1}{s(\alpha_1)}+ \cdots +
  \frac{1}{s(\alpha_m)} + \frac{1}{t}\le \frac{1}{q}$}
 \ . &
\end{eqnarray}
Here we have used the generalized H\"older inequality, and the change
of variables theorem to pass from $\|F^{(m)}\circ \psi\|_{L^t(K)}$ to
$\|F^{(m)}\|_{L^t(\psi(K))}$.  By Sobolev's embedding we have
$F^{(m)}\in L^r(\psi(K))$ for all $r$ satisfying
\begin{eqnarray}
  \label{eq:c0}
  \frac{1}{r}\gst \frac{1}{q} + \frac{m-\ell}{n}\ \;.
\end{eqnarray}
Consider, first, those terms in \eq{eq:c1} for which the right hand
side of \eq{eq:c0} is positive (if any). Let $r$ be defined by
Equation \eq{eq:c0} with $\gst$ replaced by $=$. Set
\begin{equation}
  \label{eq:c0.1}
  \frac{1}{s(\alpha_i)} = \frac{|\alpha_i|-1}{n}\ ;
\end{equation}
since $kp>n$ we have
\begin{equation}
  \label{eq:c3}
  \frac{1}{s(\alpha_i)} > \frac{1}{p} - \frac{k}{n} +
  \frac{|\alpha_i|-1}{n}\ \;,
\end{equation}
and Sobolev's embedding theorem implies that
$\|\partial^{|\alpha_i|}\psi\|_{L^{s(\alpha_i)}(K) }$ is finite. With
this choice of $r$ and of the $s(\alpha_i)$'s we have
\begin{eqnarray*}
  \sum \frac{1}{s(\alpha_i)} + \frac{1}{r} & = & \frac{|\alpha|-m}{n }
  + \frac{1}{q} +\frac{m-\ell}{n }
\\ & = & \frac{|\alpha|-\ell}{n }
  + \frac{1}{q}  \le \frac{1}{q}\ \;,
\end{eqnarray*}
so that those terms will give a finite contribution to the right hand
side of \eq{eq:c1} by setting $t=r$.

Consider, next, those terms in \eq{eq:c1} for which the right hand
side of \eq{eq:c0} vanishes. If all the $\alpha_i$'s have length one
the term in question will give a finite contribution to the right hand
side of \eq{eq:c1} by setting $t=q$. If one of the $\alpha_i$'s, say
$\alpha_1$, has length large than $1$, for $i\ge 2$ we choose the
$\alpha_i$'s as in \eq{eq:c0.1}, while we set $1/s(\alpha_1)=
(|\alpha_1|-1)/n -\epsilon>0$, with $0<\epsilon$ so chosen that
\eq{eq:c3} still holds, $\epsilon< (kp - n ) / 2pn$. Choosing
$1/t=\epsilon$ will lead to a finite contribution in \eq{eq:c1}.

It remains to consider those terms in \eq{eq:c1} for which the right hand
side of \eq{eq:c0} is negative. In this case we set $t=\infty$, and
\begin{equation}
  \label{eq:c4}
  \frac{1}{s(\alpha_i)}= \frac{|\alpha_i|}{\ell q}\quad
  \Longrightarrow \quad  \sum \frac{1}{s(\alpha_i)} =
  \frac{|\alpha|}{\ell q }\le \frac{1}{q}\ \;.
\end{equation}
By Sobolev's embedding
$\|\partial^{\alpha_i}\psi\|_{L^{s(\alpha_i)}(K) }$ will be finite
when \eq{eq:c3} holds.  Now Equation \eq{eq:c3} with $s(\alpha_i)$
defined by \eq{eq:c4} is equivalent to
\begin{equation}
  \label{eq:c5}
 \frac{|\alpha_i|}{\ell }\Big( \frac{1}{q} - \frac{\ell}{n}\Big)>
 \frac{1}{p} - \frac{k-1}{n}\ \;.
\end{equation}
The right hand side of \eq{eq:c5} is negative. If the left hand side
is positive or vanishes there is nothing to check. If both sides are
negative the worst case is obtained with $|\alpha_i|=\ell$, and the
inequality holds when \eq{eq:p-1} is an inequality. The simple
analysis of the case of equality in \eq{eq:p-1} is left to the reader.
\qed

Before returning to the proof of Proposition~\ref{pinvnew} we need one
more Lemma:
\begin{Lemma}
  \label{lproduct}
  Let $0\le m\le \ell \le k$, $q,p\in [1,\infty]$, $kp>n$. Suppose
  that $(\ell,q)$ is such that the Sobolev embedding $\wkpgl\subset
  W^{\ell,q}_\loc$ holds, \emph{cf.\/} Equations \eq{eq:p0} and
  \eq{eq:lst}. Then the product map
$$
W^{k-m,p}_\loc \times W^{\ell,q}_\loc \ \ni \ (f,g) \ \longrightarrow
fg \in W^{\ell-m,q}_\loc
$$
is continuous.
\end{Lemma}

\proof For any $0\le |\alpha|\le \ell - m$ the Leibniz rule gives
$$ \partial ^{\alpha}(fg) =
\sum_{\alpha_1+\alpha_2=\alpha}C_{\alpha_1,\alpha_2}\,  \partial
^{\alpha_1}f \, \partial^{\alpha_2}g\ ,
$$
so that on any compact set $K\subset M$ the H\"older inequality
gives
\begin{eqnarray}
  \label{eq:p1}
 & \|\partial ^{\alpha}(fg)\|_{L^q(K)} \le C
\mbox{$\displaystyle \sum_{\alpha_1+\alpha_2=\alpha}$} \|\partial
^{\alpha_1}f\|_{L^{s(\alpha_1)}(K)}
\|\partial^{\alpha_2}g\|_{L^{s(\alpha_2)}(K)}\ ,
& \\ & \label{eq:p2}
\mbox{$\displaystyle \frac{1}{s(\alpha_1)}+\frac{1}{s(\alpha_2)}\le
  \frac{1}{q} $}
\ . &
\end{eqnarray}
Let
\begin{eqnarray*}
  a(\alpha_1)& = & \frac{1}{p}+ \frac{|\alpha_1|-k+m}{n}\ ,
\\ a(\alpha_2)& = & \frac{1}{q}+ \frac{|\alpha_2|-\ell}{n}\ \;.
\end{eqnarray*}
By Sobolev's embedding we will have $\partial^{\alpha_1}f \in
{L^{s(\alpha_1)}(K)}$, $\partial^{\alpha_2}g \in {L^{s(\alpha_2)}(K)}
$ when
\begin{equation}
  \label{eq:p3}
\frac{1}{s(\alpha_i)}\gst  a(\alpha_i)\ \;.
\end{equation}
We have the following cases:
\begin{itemize}
  \item If $a(\alpha_1)\le 0$ and $a(\alpha_2)\le 0$ we set
    $s(\alpha_i)= q/2$, and we obtain
    \begin{equation}
      \label{eq:prodin}
       \|\partial ^{\alpha_1}f\|_{L^{s(\alpha_1)}(K)}
\|\partial^{\alpha_2}g\|_{L^{s(\alpha_2)}(K)} \le C
\|f\|_{W^{k-m,p}(K)} \|g\|_{W^{\ell,q}(K)}
    \end{equation}
\item If  $a(\alpha_1)> 0$ and $a(\alpha_2)> 0$ we set $1/s(\alpha_i)=
  a(\alpha_i)$ so that \eq{eq:p3} holds, and we obtain
$$    \mbox{$\displaystyle \frac{1}{s(\alpha_1)}+
  \frac{1}{s(\alpha_2)}$}
=  \mbox{$\displaystyle\frac{1}{q}+
  \frac{1}{p}+ \frac{|\alpha|- k-\ell+m }{n}$}
\le \mbox{$\displaystyle\frac{1}{q}+
  \frac{1}{p}- \frac{ k}{n}< \frac{1}{q}$}
 \ ,
$$
since $kp>n$, so that \eq{eq:p2} holds. We note that \eq{eq:prodin} is
again satisfied.
\item If $a(\alpha_1)=0$ and $a(\alpha_2)>0$ we have
  $$
  \frac{|\alpha_2|}{n}= \frac{|\alpha|-|\alpha_1|}{n}=
  \frac{|\alpha|- k+m }{n}+ \frac{1}{p}\ ,
$$
so that
\begin{equation}
  \label{eq:pr-1}
  a(\alpha_2) = \frac{1}{q}+ \frac{1}{p}+ \frac{|\alpha|- k-\ell+m
  }{n} \le \frac{1}{q}+ \frac{1}{p}- \frac{ k}{n} < \frac{1}{q} \ \;.
\end{equation}
Let $\epsilon$ be  any number satisfying $0<\epsilon<(kp-n)/np$, set
$$
\frac{1}{s(\alpha_1)}= \frac{\epsilon}{2}\ , \quad
\frac{1}{s(\alpha_2)}= a(\alpha_2)+\frac{\epsilon}{2}\ .
$$
Decreasing $\epsilon$ if necessary we will have $s(\alpha_2)>0$.
Then \eq{eq:p2}, \eq{eq:p3} and \eq{eq:prodin} hold by the calculation
in Equation \eq{eq:pr-1}. A similar analysis takes care of the case
$a(\alpha_1)>0$ and $a(\alpha_2)=0$.
\item If $a(\alpha_1)<0$ and $a(\alpha_2)>0$ we set
  $s(\alpha_1)=\infty$ and $s(\alpha_2)=q$.
\item  If $a(\alpha_1)>0$ and $a(\alpha_2)<0$ we set
  $s(\alpha_1)=q$ and $s(\alpha_2)=\infty$; \eq{eq:p2} obviously
  holds, while
  $$
  \frac{1}{s(\alpha_1)}-a(\alpha_1)= \frac{1}{q}- \frac{1}{p}+
  \frac{ k-m -|\alpha_1|}{n} \ge \frac{1}{q}- \frac{1}{p}+ \frac{
    k-\ell}{n}$$
  which is non-negative by \eq{eq:p0}, hence
  \eq{eq:p3} and \eq{eq:prodin} hold again.
\end{itemize}
This establishes that $fg\in  W^{\ell-m,q}_\loc$. The continuity
follows immediately from the inequality
$$
\|fg\|_{W^{\ell-m,q}(K)}\le C \|f\|_{W^{k-m,p}(K)}\|g\|_{W^{\ell,q}(K)}
$$
which has been established during the proof.
\qed

We can pass now to the proof of point 2 of Proposition~\ref{pinvnew}.
Let two coordinate systems on $M$ be given related to each other by a
map $\psi\in \wkpl$, set $\chi\equiv \psi^{-1}$. Let
${t^{\alpha_1\ldots \alpha_k}}_{\beta_1\ldots \beta_s}$ and ${{\hat
    t}^{\mu_1\ldots \mu_k}}{}_{\nu_1\ldots \nu_s}$ be the coordinate
components of a tensor field $t$, with ${ t^{\alpha_1\ldots
    \alpha_k}}_{\beta_1\ldots \beta_s} \in \wlqg$. We have the
transformation rule
\begin{equation}
  \label{eq:tr}
  {{\hat t}^{\mu_1\ldots\mu_k}}{}_{\nu_1\ldots \nu_s}(x)=
{  t^{\alpha_1\ldots \alpha_k}}_{\beta_1\ldots \beta_s}(\psi(x))
\frac{\partial \psi^{\beta_1}}{\partial x^{\nu_1}}(x)
\cdots
\frac{\partial \psi^{\beta_s}}{\partial x^{\nu_s}}(x)
\frac{\partial \chi^{\mu_1}}{\partial x^{\alpha_1}}(\psi(x)) \cdots
\frac{\partial \chi^{\mu_k}}{\partial x^{\alpha_k}}(\psi(x))
\ .\end{equation}
Now the matrix $\frac{\partial \chi^{\mu}}{\partial x^{\alpha}}\circ\psi$ is
the inverse matrix to$ \frac{\partial \psi^{\beta}}{\partial
  x^{\nu}}$, so the components of the former are rational
function of those of the latter.  We recall the
Gagliardo--Moser--Nirenberg inequalities ({\emph{cf.,
    e.g.}\/~\cite[Corollaries~6.4.4 and 6.4.5]{Hormander97})
\begin{eqnarray}
  \label{eq:moser}&&
  \forall \ f,g \in \wlqg\cap L^\infty\qquad \|fg\|_{\wlqg}\le
  C_1(\|f\|_{L^\infty} \|g\|_{\wlqg}+  \|f\|_{\wlqg}\|g\|_{L^\infty})\ ,
\\ \label{comp} &&
  \forall \ f \in \wlqg\cap L^\infty\qquad \|F(f)\|_{\wlqg}\le
  C_2(\|f\|_{L^\infty})(1+
  \|f\|_{\wlqg})\phantom{\|f\|_{\wlqg}\|g\|_{L^\infty})\ ,}
\ ,
\end{eqnarray} for some $f$ and $g$ independent constant $C_1$, and
for some constant $C_2(\|f\|_{L^\infty})$ depending upon $f$ only
through its $L^\infty$ norm. Here one assumes that $F$ is a smooth
function of its argument which is allowed to take values in $\bR^N$,
and the integrals are taken over compact sets.  Equations
\eq{eq:moser}--\eq{comp} show that the $\frac{\partial
  \chi^{\mu}}{\partial x^{\alpha}}\circ \psi$'s are $\wkpgl$ functions
of their arguments. Lemma \ref{lproduct} implies that
$$\frac{\partial \psi^{\beta_1}}{\partial x^{\nu_1}}(x)
\cdots
\frac{\partial \psi^{\beta_s}}{\partial x^{\nu_s}}(x)
\frac{\partial \chi^{\mu_1}}{\partial x^{\alpha_1}}(\psi(x)) \cdots
\frac{\partial \chi^{\mu_k}}{\partial x^{\alpha_k}}(\psi(x))
\in \wkpgl
\ .
$$
It follows that the right hand side of Equation \eq{eq:tr} is of the
form
$$
t\circ \psi \, A\ ,\qquad t\in W^{\ell,q}_\loc\ , \ \psi \in
\wkpl\ , \ A \in \wkpgl\ .
$$
Lemma \eq{lcomp} implies that $t\circ \psi \in  W^{\ell,q}_\loc$, and
Lemma \ref{lproduct} with $m=0$ shows that the right hand side of
Equation \eq{eq:tr} is in $ W^{\ell,q}_\loc$, as desired.
\qed

We wish to extend the above discussion to spinor fields; this requires
the introduction of orthonormal frames, and hence of the metric.
Consider, then, a $\wkpl$ manifold $M$ with a strictly positive
definite symmetric two-covariant tensor field $g$. We shall say that
$(M,g)$ is a \emph{\Riemn $\wkpl$ manifold} if $M$ is a $\wkpl$
manifold and if $g$ is a \Riemn metric of $\wkpgl$ differentiability
class. This is an invariantly defined notion by Proposition
\ref{pinvnew}.  Before proceeding further we note the following:

\begin{Proposition}
  \label{Priem} Let $(M,g)$ be a $\wkpl$ manifold with a \Riemn
  metric of $\wkpgl$ differentiability class, $kp>n$, $p\in
  [1,\infty]$. Then the following hold
  \begin{enumerate}
  \item In any coordinate system in the $\wkpl$ atlas the Christoffel
    coefficients $\Gamma^i_{jk}$ satisfy
$$ \Gamma^i_{jk} \in W^{k-1,p}_\loc \ . $$
\item The Riemann tensor is of $ W^{k-2,p}_\loc$ differentiability class.
\item The curvature scalar $R\equiv g^{ij}{R^k}_{ikj}$ is of $
  W^{k-2,p}_\loc$ differentiability class.
\item \label{tensorderivative} Assume that $\ell\ge 1$ and suppose
  that $(\ell,q)$ is such that the Sobolev embedding $\wkpgl\subset
  W^{\ell,q}_\loc$ holds, \emph{cf.\/} Equations \eq{eq:p0} and
  \eq{eq:lst}. Let $t$ be a tensor field of $ W^{\ell,q}_\loc$
  differentiability class, then for any vector field $X\in
  W^{k,p}_\loc$ we have
$$
X^i\nabla_i t \in W^{\ell-1,q}_\loc\ .$$
  \end{enumerate}
\end{Proposition}

\proof
1. By definition of the $\Gamma^i_{jk}$'s we have
$$
\Gamma^i_{jk}= \frac{1}{2}g^{i\ell}\Big(\frac{\partial g_{\ell
    j}}{\partial x^k} + \frac{\partial g_{\ell
    k}}{\partial x^j} - \frac{\partial g_{j
    k}}{\partial x^\ell}\Big)\ .
$$
By \eq{eq:moser} we have $g^{i\ell}\in \wkpgl$, by definition the
derivatives of the metric are in $ W^{k-1,p}_\loc$, and lemma
\ref{lproduct} with $m=0$ and $(\ell,q)=(k-1,p)$ gives the result.

2. By definition of the curvature tensor we have
$$
{R^i}_{jk\ell}= \frac{\partial \Gamma^i_{j\ell} }{\partial x^k }-
\frac{\partial \Gamma^i_{jk} }{\partial x^\ell }+
\Gamma^i_{mk}\Gamma^m_{j\ell} -  \Gamma^i_{m\ell}\Gamma^m_{jk}\ ,
$$
Lemma \ref{lproduct} with $m=1$ and  $(\ell,q)=(k-1,p)$ shows that the
product terms above are in $ W^{k-2,p}_\loc$, and the result follows
from point 1.

3. This claim follows immediately from point 2,  from $g^{i\ell}\in
\wkpgl$ and from Lemma \ref{lproduct} with $m=0$ and $(\ell,q)=(k-2,p)$.

4. Let $t$ be a tensor field of  $ W^{\ell,q}_\loc$ differentiability
  class, in local coordinates we have
$$\nabla_i {  t^{\alpha_1\ldots \alpha_k}}_{\beta_1\ldots \beta_s} =
\partial_i {  t^{\alpha_1\ldots \alpha_k}}_{\beta_1\ldots \beta_s}+
\Gamma^{\alpha_1}_{i\sigma }{  t^{\sigma\ldots
    \alpha_k}}_{\beta_1\ldots \beta_s} +  \ldots
- \Gamma^{\sigma}_{i\beta_1 }{  t^{\alpha_1\ldots
    \alpha_k}}_{\sigma\ldots \beta_s} -  \ldots\ .
$$
The $\Gamma$'s are in $ W^{k-1,p}_\loc$ by point 1, thus the
product terms are in $W^{\ell-1,q}_\loc$ by Lemma \ref{lproduct} with
$m=1$.  The claim about $X^i\nabla_i t $ follows again from Lemma
\ref{lproduct}. \qed

Let $(M,g)$ be a $\wkpl$ {\Riemn manifold} and let $\ON$ be the bundle
of $g$--orthonormal frames on $M$. We can equip $\ON$ with a $\wkpgl$
structure by considering only those $g$--orthonormal sets of vector
fields which are all of $\wkpgl$ differentiability class. Let us start
by showing that the set of such (locally defined) frames is not empty.
On $\cO$, the domain of a coordinate system $(x^i)$, we
can construct a $g$--orthonormal frame $e_j= {e_j}^i\partial/\partial
x^i$ by performing a Gram--Schmidt orthonormalisation of the basis
$\{{\partial}/{\partial x^i}\}$. By construction the coordinate
coefficients ${e_j}^i$ of the vector fields $e_j$ are smooth functions
of $g_{ij}$ (at least on a neighborhood of the range of values taken
by $g_{ij}$), where the $g_{ij}$'s are the coordinate coefficients of
the metric $g$, $g=g_{ij}dx^i dx^j$.  Since $kp>n$,
\eq{eq:moser}--\eq{comp} applied to the ${e_j}^i$ considered as
functions of the $g_{ij}$ shows that the vector fields $e_j$ are
indeed of $\wkpgl$ differentiability class, as desired.

The following shows that the $\wkpgl$ structure of $\ON$ is an
invariantly defined property of a $\wkpl$ \Riemn manifold:

\begin{Proposition}
  \label{lon} Any two (globally or locally defined) $g$--orthonormal
  frames of $\wkpgl$ differentiability class are related to each other
  by a $O(n)$--rotation of $\wkpgl$ differentiability class.
\end{Proposition}

\proof Consider two locally defined $g$--orthonormal frames $e_i$
and $f_i$, $i=1,\ldots, n=\mathrm{dim}\, M$, of $\wkpgl$
differentiability class. In particular each of the $e_i$ and
$f_i$' is a vector field of $\wkpgl$ differentiability class,
which is invariantly defined by Proposition~\ref{pinvnew}, so that
it is sufficient to prove the result in any coordinate system in
the $\wkpl$ atlas on $M$. In such a coordinate system $\{x^i\}$ we
can write $e_j={e_j}^i\frac{\partial}{\partial x^i}$,
$f_j={f_j}^i\frac{\partial}{\partial x^i}$, for some functions
${e_j}^i,{f_j}^i\in\wkpgl$.  Since both frames are orthonormal
there exists an $O(n)$--valued function ${w_i}^j$ such that
\begin{equation}
  \label{eq:omega}
 e_i = {w_i}^j f_j \ \;.
\end{equation}
It follows that
$$ {w_i}^j = {e_i}^k {f^j}_k \ , $$
where ${f^j}_k$ is the matrix inverse to ${f_j}^k$. We have
${f^j}_k\in \wkpgl$ by \eq{comp}, thus ${w_i}^j \in \wkpgl$ by
\eq{eq:moser}, hence the result. \qed

 Now suppose that $M$
has a spin structure, namely a $\mathrm{Spin}$--principal bundle
$\tmcF$ which double-covers the principal bundle $\mcF$ of
$g$--orthonormal frames of $(M,g)$: \bel{spinpro} 0 \to \mathbb{Z}_2
\to \tmcf \stackrel{\pi}{\to} \mcF \to M\;.\ee We note that the
obstruction to the existence of such structures is purely topological,
(\emph{cf., e.g.}, \cite[Chapter II]{LawsonMichelsohn89}) and
therefore independent of the choice of the metric and
differentiability class.  A bundle of spinors $\mcV =
\tmcF\times_{T}V$ is a vector bundle associated to $\tmcF$ and a
representation
$$
T:\mbox{\rm Spin}\to \mbox{\rm End}(V)\;,
$$
for some finite-dimensional vector space $V$. A choice of
$g$--orthonormal frame $e=(e_i)$ of $\wkpgl$ differentiability class
defined on an open set $\cO\subset M$ determines a local section of
$\mcF$.  This lifts to a section of $\tmcF$, which in turn is
associated with a local orthonormal frame $\phi = (\phi_I)$ in $\mcV$.
Let $\cU$ be another open set with a $g$--orthonormal frame
$e'=(e'_i)$, so by Proposition~\ref{lon} there exists an $O(n)$ valued
map $w=({w_i}^j)$ of $\wkpgl$ differentiability class such
that the frames $e, e'$ are related by $e_i = {w_i}^je'_j$ on
$\cO\cap \cU$.  The map $w:\cO\cap \cU \to \mathrm{SO}$ lifts to
$\tilde{w}:\cO\cap \cU\to \mathrm{Spin}$.  This lift is not
unique, but the possible lifts differ only by a fixed nontrivial
element $z$ of the centre $Z(\mathrm{Spin}) \simeq \bZ_2$.  The
corresponding spin frames $\phi, \phi'$ are related by $\phi =
T({\tilde{w}})\phi'$ or $T(z{\tilde{w}})\phi'$.  Analyticity of
the local inverse $\pi^{-1}(\cdot)$ and the inequality \eq{comp} show
that $T({\tilde{w}}), T(z\tilde{w})\in \wkpgl$, and it follows
from \eq{eq:moser} that the spin frames $\phi$ on ${\cU}$ and $\phi'$
on $\cO$ are $\wkpgl$ compatible on $\cU\cap\cO$.  This establishes
the following result:

\begin{Proposition}
  \label{pspin} Let $(M,g)$ be a $\wkpl$ spin manifold with a \Riemn
  metric of $\wkpgl$ differentiability class. Then every spinor bundle
  carries a natural $\wkpgl$ differentiable structure.
\end{Proposition}


An argument similar to that of  Proposition~\ref{pinvnew} shows:

\begin{Proposition}
  \label{pspininv} Let $(M,g)$ be a $\wkpl$ spin manifold with a \Riemn
  metric of $\wkpgl$ differentiability class, $kp>n$, $p\in
  [1,\infty]$. Let $(\ell,q)$ be such that the Sobolev embedding
$$
\wkpgl\subset W^{\ell,q}_\loc
$$
holds, \emph{cf.\/} Equations \eq{eq:p0} and \eq{eq:lst}. Then the
space of $W^{\ell,q}_\loc$ spinor fields is invariantly defined.
\end{Proposition}

To proceed further,  we recall the
definition of the covariant derivative of a spinor field. Let $V$,
$e_i$, $\cO\subset M$ and $\gamma_i$ be as before, and let
$\phi=(\phi_{I})$ be the spinor frame corresponding to the orthonormal
frame $e=(e_i)$. This defines a preferred local spin frame, with
respect to which the Clifford action is represented by locally
constant matrices $\gamma_i$.  Let $\psi$ be a spinor field over
$\cO$; the spinor covariant derivative of $\phi$ is given in terms of
the orthonormal frame connection matrix
$\omega_{ij}(e_{k})=g(e_i,\nabla_{e_k}e_j)$ and the spinor frame
components $\psi=\psi^{I}\phi_{I}$ by \bref{A:gpsi},
and the Dirac operator of $\nabla $ on $S$ is defined by
\begin{equation}
  \mfd \psi = \gamma^i\nabla_{e_{i}}\psi\ .
\label{A:dirac1}
\end{equation}
\begin{Proposition}
  \label{pspinnew} Let $(M,g)$ be a $\wkpl$ manifold with a \Riemn
  metric of $\wkpgl$ differentiability class, $kp>n$, $p\in
  [1,\infty]$. Then the following hold
  \begin{enumerate}
  \item Let $e^i$ be any $g$--orthonormal frame of $\wkpgl$
    differentiability class, then the spin connection coefficients
    $\omega_k$ defined as $\nabla_{e_k}\psi = e_k(\psi) + \omega_k
    \psi$ satisfy
    $$
    \omega_k\in W^{k-1,p}_\loc \ . $$
  \item If $(\ell,q)$ is as in
    Proposition~\ref{pspininv} with $\ell\ge 1$, and if $X$ is a
    vector field of $\wkpg$ differentiability class, then $\nabla_X$
    maps continuously $W^{\ell,q}_\loc$ spinor fields to
    $W^{\ell-1,q}_\loc$ spinor fields:
$$
W^{\ell,q}_\loc \ni \phi \longrightarrow \nabla_X \phi \in
W^{\ell-1,q}_\loc\ .
$$
In particular the Dirac operator maps continuously
$W^{\ell,q}_\loc$ to $W^{\ell-1,q}_\loc$.
  \end{enumerate}
\end{Proposition}

\proof To prove point 1 choose a spin frame in which the
$\clg {e^i}$'s are point independent matrices.  Then
$$
\omega_k\equiv -\frac{1}{4}\clg {e^i} \clg {e^j}\omega_{ij}(e_{k})
\ ,
$$
with $\omega_{ij}(e_{k})\equiv g(e_i,\nabla_{e_k}e_j)$. The claim
that $\omega_{ij}(e_{k})\in W^{k-1,p}_\loc$ follows immediately from
point \ref{tensorderivative} of Proposition~\ref{Priem} and from Lemma
\ref{lproduct}. The result in any spin frame follows from the
transformation rule of the connection coefficients under changes of
frames and from Lemma \ref{lproduct}. The proof of point 2 follows
that of point \ref{tensorderivative} of Proposition~\ref{Priem} and
will be omitted.  \qed

Given a smooth metric in a neighbourhood of a compact
boundary $\partial M$, geodesics normal to the boundary determine a
diffeomorphism of $Y\times I$ with a neighbourhood of $\partial
M \simeq Y$, such that in adapted coordinates $v =
(y^A,x)$, $x \in [0, x_0), y^A \in \co_i$ we have
\begin{equation} g^{xx}\equiv 1, \quad g^{xA}\equiv 0.
\label{(AGC.1)}
\end{equation}
The diffeomorphism determines a \emph{tubular neighbourhood} of
$\partial M$ and the resulting coordinates are called \emph{Gaussian
coordinates}.
If the metric has only low
differentiability then uniqueness of the geodesic equation may fail,
and the existence of Gaussian coordinates becomes problematic.
However, for our
applications it is sufficient for (\ref{(AGC.1)}) to hold only
approximately near $\partial M$, in which case we may rely on
the following result.


%
\begin{Proposition}
  (Almost Gaussian tubular neighbourhood coordinates for $\partial M$).
  \label{PAGC.1} Let $k\in\bN$, $\ell\in\bN\cup\{0\}$, and suppose
  $(M,g)$ be a $W^{k+1,\ptwo}_\loc$ Riemannian manifold with metric
  $g\in W^{k,\ptwo}_\loc(M)$, $(k-\ell)\ptwo>n$.
Let $Y\subset M$ be a compact connected component of the boundary of
$M$ with $Y$ of $W^{k+1,\ptwo}_\loc$ differentiability class. There is
a neighbourhood $\cO$ of $Y\subset M$ and $x\in W^{k+1,p}(\cO)$, and a
diffeomorphism $\cO\simeq Y\times I$, $I=[0,x_0)$, which determines
coordinates $(v^i) = (y^A,x) \in \cO$ such that $Y\cap \cO =\{x=0\}$ and
\begin{eqnarray}
&
g(dv^i, dv^j) = g^{ij} \in  W^{k,\ptwo}(\cO )\ ,
\label{(AGC.0.1)} & \\ &
g(dx\;, dx) -1 = O(x^{\ell+\sigma}) , \label{(AGC.2)}
 & \\ &
g(dx\;, dy^A) =  O(x^{\ell+\sigma}) \;.,\label{(AGC.3)} &
\end{eqnarray}
for some $\sigma>0$.
\end{Proposition}
 \textbf{Remarks:} 1. A similar result for $C^{k,1}$ metrics, $k\ge
 1$, follows from \cite[Appendix~B]{AndChDiss}.

2. Similar results hold for \Riemn manifolds provided
$Y$ is non--characteristic, and when $Y$ is a hypersurface in $M$.

\proof\ If $k = \infty$ we can use Gauss coordinates near $Y$, and the
result follows.  Suppose thus that $k < \infty$, let $x$ be any
defining function for $Y$ and let $\cO_\alpha$ be any conditionally
compact coordinate neighborhood of $Y$, with $g^{ij} \equiv g (dv^i,
dv^j) \in W^{k,\ptwo}({\cO_\alpha})$.  Passing to a subset of
$\cO_\alpha$ if necessary without loss of generality we may assume
$\cO_\alpha \approx [0, x_0)\times \cutut_\alpha$,
$\cutut_\alpha\subset Y$.  Coordinate systems of this form will be
called cylindrical.

We construct a suitable $W^{k+1,p}$ coordinate change
$(\bar{y}^A,\bar{x})$ in $\cO_\alpha\times I$ by noting first
that $d\bar{x}=\frac{\partial \bar{x}}{\partial x} dx + \frac{\partial
  \bar{x}}{\partial y^A} dy^A$, $d\bar{y}^A=\frac{\partial
  \bar{y}^A}{\partial x\ } dx + \frac{\partial \bar{y}^A}{\partial
  y^B} dy^B$.  Thus if $\bar{x}$ is also a boundary coordinate, so
$\bar{x}(y,0)=0$ and $\frac{\partial\bar{x}}{\partial y^A}=0$ on $Y$,
then the metric coefficients satisfy
\begin{eqnarray*}
  g^{\bar{x}\bar{x}} &=& g(d\bar{x},d\bar{x})
  \\
  &=&
  \left(\pd{\bar{x}}{x}\right)^2 g^{xx}
  + 2 \pd{\bar{x}}{x}\pd{\bar{x}}{y^A}
  \,g^{xA} + \pd{\bar{x}}{y^A}\pd{\bar{x}}{y^B}\,g^{AB}
  \\
  &=&  \left(\pd{\bar{x}}{x}\right)^2 g^{xx}\quad \mathrm{on}\ Y=\{x=0\},
  \\
  g^{\bar{x}\bar{A}} &=& g(d\bar{x},d\bar{y}^A)
  \\
  &=&
  \pd{\bar{x}}{x}\left( \pd{\bar{y}^A}{x}\, g^{xx} +
   \pd{\bar{y}^A}{y^B} \,g^{xB}\right)
   \quad \mathrm{on}\ Y\;.
\end{eqnarray*}
 Since
$g^{xx}\in W^{k,p}$, the restriction $g^{xx}|_Y$ lies in the Besov
space $\Lambda^{p,p}_{k-1/p}(Y)$ (see \cite[\S VI.4.4]{Stein70}, or
\cite[Theorem VII.1]{JonssonWallin84}) and
there is an extension $\bar{x}=\bar{x}(y,x) \in W^{k+1,p}(\co_\alpha\times I)$
satisfying the conditions
\[
\bar{x}(y,0) = 0,\quad \frac{\partial \bar{x}}{\partial x}(y,0) =
(g^{xx})^{-1/2}(y,0)\;,
\]
for $y = (y^A)\in \cO_\alpha$  (\cite[\S VI.6]{Stein70},
\cite[Theorem VII.3]{JonssonWallin84}).  This implies
\[
g^{\bar{x}\bar{x}} = g(d\bar{x},d\bar{x}) \in W^{k,p}(\co_\alpha\times I)
\]
and $g^{\bar{x}\bar{x}}=1$ when $x=0$.  Similarly, there is
$f^A(y,x)\in W^{k+1,p}(\co_\alpha\times I)$ such that
\[
f^A(y,0) = 0, \quad \frac{\partial f^A}{\partial x}(y,0) =
{}- g^{xA}/g^{xx}(y,0) \;,
\]
so the coordinates $(\bar{y},\bar{x})$, $\bar{y}^A(y,x)=y^A+f^A(y,x)$
also satisfy
\[
g^{\bar{A}\bar{x}}(y,0) = g(d\bar{x},d\bar{y}^A) = 0
\]
on $Y$, since $d\bar{y}^A = dy^A - g^{xA}/g^{xx}\,dx$.



As shown in Proposition~\ref{pinvnew}, in the new coordinate system
$(\bar v^i)\equiv(\bar{y}^A,\bar{x})$ we still have $g({ d\bar v^i,
  d\bar v^j}) \in W^{k,\ptwo}(\cO_\alpha)$, so that by embedding theorems the
metric coefficients are $\sigma$--H\"older continuous on ${\cO_\alpha}$, for
some $\sigma>0$, and
\begin{equation}
{g^{\bar{x}\bar{x}}} - 1=O(x^\sigma),\quad
g^{\bar{x}\bar{A}}=O(x^\sigma)\;.
\label{(AGC.10.-2)}
\end{equation}
Let $\phi , f^A \in W^{k+1,\ptwo}(\cO_\alpha)$,
and consider the effect of the change
\begin{equation} \bar x = x + \phi(y,x) , \label{(AGC.10)}
\end{equation} \begin{equation} \bar{y}^{A} = y^A + f^A(y,x) ,
\label{(AGC.11)}
\end{equation}
where $\phi(y,0) =0$, $f^A(y,0)=0$.  Then
\begin{eqnarray}
&  g^{\bar{x}\bar{x}} = g^{xx} +
2g^{xi}{\partial\phi\over\partial v^i} +
g^{ij}{\partial\phi\over\partial v^{i} }{\partial\phi\over\partial
  v^{j} }\,,
\label{(AGC.12)}
& \\ &
 g^{\bar{x}\bar{A}} \equiv g(d\bar x, d\bar y^{A}) = g^{xA} + g^{xi}{\partial
  f^A\over \partial v^i} + g^{Ai}{\partial\phi\over\partial v^i} +
g^{ij}{\partial\phi\over\partial v^i}{\partial f^A\over\partial v^j}
.\label{(AGC.13)} &
\end{eqnarray}
Suppose that for some $\ell \geq 0$ we have
\begin{equation}
g^{xx}-1 = O(x^{\ell+\sigma}) ,\quad g^{xA} =
O(x^{\ell+\sigma})\,.\label{(AGC.14)}
\end{equation}
This holds for $\ell = 0$ by \eq{(AGC.10.-2)} and we establish the
general case by induction.  Again by restriction and extension results
\cite{Stein70,JonssonWallin84} there exist $\phi, f^A \in
W^{k+1,\ptwo}({\cal O}_\alpha)$ satisfying
\begin{eqnarray*}
{\partial^{\ell +1} \phi\over \partial x }\Big|_{(y,0)} &=& {}-{1\over
  2}{\partial^\ell g^{xx}\over \partial x^{\ell}}\Big|_{(y,0)}\;,
\\
{\partial^{\ell +1} f^A\over\partial x^{\ell +1} }\Big |_{(y,0)}& =& {}-
{\partial^\ell
g^{xA}\over \partial x^\ell}\Big|_{(y,0)} \;,
\end{eqnarray*}
while all the lower order $x$--derivatives of $\phi,f^A$ vanish at
$x=0$.  Passing to coordinates $(\bar{y},\bar{x})$ on a (possibly smaller)
cylindrical neighborhood $\cO_\alpha$, one finds from
(\ref{(AGC.12)})--(\ref{(AGC.13)}) that (\ref{(AGC.14)}) still holds
and moreover,
\begin{eqnarray*}
g^{\bar{x}\bar{x}}-1  &=& O(x^{\ell+1+\sigma}) ,
\\
g^{\bar{x}\bar{A}} &=& O(x^{\ell+1+\sigma}) ,
\end{eqnarray*}
Dropping bars one finds that (\ref{(AGC.14)}) holds with $\ell$
replaced by $\ell + 1$, and the induction step is complete.

Finally we show that the local charts can be combined to form a
tubular neighbourhood diffeomorphism.  It follows from \bref{(AGC.12)}
that if $x,\bar{x}$ both satisfy \bref{(AGC.10.-2)} and vanish on $Y$,
then $g(dx,d\bar{x})=1$.  In particular, by combining the functions
$\bar{x}_\alpha$ from each of the local coordinate charts $\cO_\alpha$
using a subordinate partition of unity $\phi_\alpha$, the
function $x=\Sigma_\alpha\phi_\alpha \bar{x}_\alpha$ satisfies
$x\in W^{k+1,p}(Y\times I)$ and $x=0$, $g^{xx}=1$ on $Y$.

In order to construct a diffeomorphism with $Y\times I$ we need to
construct a similar averaging of the $\bar{y}^A$ coordinate functions.
Fix a smooth embedding $\Phi:Y\to \bR^K$ and let
$\Pi_{\Phi(Y)}:\cN\subset \bR^K\to \Phi(Y)$ be the orthogonal
projection in $\bR^K$ from a tubular neighbourhood $\cN$ back to
$\Phi(Y)$.  Let $y_\alpha=(y^A_\alpha):\cU_\alpha\subset Y\to
\bR^{n-1}$ denote both the coordinates of a $C^\infty$ chart on $Y$,
and their natural extension to
$y_\alpha=(y^A_\alpha):\cO_\alpha=\cU_\alpha\times I\subset Y\times I\to
\bR^{n-1}$.  Let
$\bar{y}_\alpha=(\bar{y}^A_\alpha):\cO_\alpha\to\bR^{n-1}$ be the
functions constructed above, so there is a neighbourhood
$\tilde{\cO}_\alpha\subset \cO_\alpha$ containing $Y=Y\times\{0\}$
such that $y^{-1}_\alpha\circ \bar{y}_\alpha:\tilde\cO_\alpha \to Y$.
Choose a finite covering
$\tilde\cO_\alpha$ of $Y\times I$ with
subordinate partition of unity $\tilde\phi_\alpha$ and define
$\Psi:Y\times I\to Y$,
\[
\Psi(p) = \Phi^{-1}\circ \Pi_{\Phi(Y)} \left(\textstyle{\sum_\alpha}
  \tilde\phi_\alpha(p) \Phi(y^{-1}_\alpha\circ\bar{y}_\alpha(p)) \right) \;.
\]
Since $\bar{y}^A_\alpha(y,0)=y^A$, $\Psi|_{Y\times\{0\}} =Id$ and
$(\Psi,x)$ defines a diffeomorphism of $Y\times I$.
Now for any $C^\infty$ chart $y=(y^A)$ on $Y$,
$\tilde{y}^A := y^A\circ\Psi$ defines a chart on $Y\times I$ by
$p\mapsto(\tilde{y}^A(p),x(p))$, which satisfies $\tilde{y}^A(y,0)=y^A$.
Moreover, $d\tilde{y}^A(y,0) = d\bar{y}^A(y,0)$, so
$g^{x\tilde{A}}(y,0)=0$.  The condition $g^{xx}(y,0)=1$ is not
affected by changes in the $y$-coordinate, so $(\Psi,x)$ defines the
required tubular neighbourhood.
  \QED

%
 \bibliographystyle{amsplain}
 \bibliography{Energy,GRPDE,PDE,GR,RBmaster,Chrusciel,DiracPTC}
%


\end{document}